\documentclass[a4paper,twoside,reqno,11pt]{article}
\usepackage{textcomp}
\usepackage{amsmath}
\usepackage{mathrsfs}
\usepackage{amssymb}
\usepackage{amsthm}
\usepackage{enumerate}
\usepackage{mathdots}
\usepackage[left=1.5in, right=1.5in, top=1in, bottom=1in]{geometry}
\usepackage{verbatim}
\usepackage{mathabx}
\usepackage{hyperref}
\usepackage{mathtools}
\usepackage{setspace}
\usepackage[titletoc]{appendix}
\usepackage{titlesec}

\usepackage{tikz}
\usetikzlibrary{arrows}
\usepackage{tikz-cd}
\usepackage{float}

\numberwithin{equation}{subsection}
\numberwithin{figure}{section}
\newtheorem{theorem}[equation]{Theorem}
\newtheorem{corollary}[equation]{Corollary}
\newtheorem{lemma}[equation]{Lemma}
\newtheorem{proposition}[equation]{Proposition}
\newtheorem{conjecture}[equation]{Conjecture}
\theoremstyle{definition}
\newtheorem{definition}[equation]{Definition}

\newtheorem{remark}[equation]{Remark}

\newcommand{\la}{\langle} \newcommand{\ra}{\rangle}

\newcommand{\nA}{{\mathbb A}}

\newcommand{\nC}{{\mathbb C}}

\newcommand{\nQ}{{\mathbb Q}}

\newcommand{\cA}{{\mathcal A}}
\newcommand{\cB}{{\mathcal B}}

\newcommand{\cH}{{\mathcal H}}
\newcommand{\cI}{{\mathcal I}}

\newcommand{\cK}{{\mathcal K}}

\newcommand{\cR}{{\mathcal R}}

\newcommand{\cV}{{\mathcal V}}

\newcommand{\cY}{{\mathcal Y}}

\newcommand{\sF}{{\mathscr F}}
\newcommand{\sG}{{\mathscr G}}
\newcommand{\sH}{{\mathscr H}}

 \DeclareMathOperator{\diag}{diag}

 \DeclareMathOperator{\N}{N}

 \DeclareMathOperator{\rank}{rank}
 
\DeclareMathOperator{\spn}{span}

\DeclareMathOperator{\GL}{GL}

\DeclareMathOperator{\Lie}{Lie}

\DeclareMathOperator{\Res}{Res}
\DeclareMathOperator{\Ind}{Ind}
\DeclareMathOperator{\Hom}{Hom}

\DeclareMathOperator{\graph}{graph}
\DeclareMathOperator{\Sym}{Sym}

\allowdisplaybreaks
\begin{document}
\begin{center}
\begin{Large}
\textbf{
THE BURGER-SARNAK METHOD AND OPERATIONS ON THE UNITARY DUALS OF CLASSICAL GROUPS}\end{Large}\\
	\vspace{1cm}
\begin{large}ANDREW HENDRICKSON\\
NATIONAL UNIVERSITY OF SINGAPORE\end{large}\end{center}
\vspace{1cm}
\textbf{Abstract.} \quad The Burger-Sarnak method shows that the restriction of an automorphic representation of a reductive group to a reductive subgroup has automorphic support. Clozel has conjectured a qualitative refinement of this result, which were first verified and quantified in the $GL_n$ case by Venkatesh. We give a proof of this conjecture in the classical group cases.
\vspace{1cm}
\tableofcontents
\newpage
\section{Introduction}
\subsection{Background}
Let $G$ be a split reductive group over a number field $F$, and let $\nA$ be the ad\`{e}le ring of $F$. In \cite{Arthur2013}, Arthur has given a decomposition of $L^2(G(F)\backslash G(\nA))$ in terms of global $A$-parameters, which are equivalence classes of maps
$$\psi:L_F\times SL_2(\nC)\rightarrow \check{G},$$
where $L_F$ is the conjectural Langlands group and $\check{G}$ is the Langlands dual group of $G$, where $\psi|_{L_F}$ satisfies some boundedness condition and $\psi|_{SL_2(\nC)}$ is an algebraic homomorphism. The global $A$-parameters give rise to the local $A$-parameters 
$$\psi_v:L_{F_v}\times SL_2(\nC)\rightarrow \check{G}$$
for each place $v$, where $L_{F_v}$ is the Weil group $W_{F_v}$ in the archimedean case or the Weil-Deligne group $W_{F_v}\times SL_2(\nC)$ in the non-archimedean case. To each local $A$-parameter is associated a finite set $A_{\psi_v}$ of unitary representations of $G(F_v)$, which is called a local $A$-packet. The representations in the local $A$-packets, which we refer to as \textit{Arthur type} representations, are then used to build the global $A$-packets, from which a decomposition of $L^2(G(F)\backslash G(\nA))$ is obtained.

Let $H\subset G$ be a split reductive subgroup over $F$. For an Arthur type representation $\pi$ of $G(F_v)$ in some local $A$-packet, one thing we are interested in describing is the direct integral decomposition of $\pi|_{H(F_v)}$. Our starting point to this end is a theorem of Burger, Li, and Sarnark first stated in \cite{Burger1992} and \cite{Burger1991} and later proved in full generality by Clozel and Ullmo in \cite{Clozel2003}.
\begin{theorem}[Burger-Li-Sarnak]\label{BSthm}
	Let $S$ be a finite set of places of $F$, and set $F_S=\prod_{v\in S}F_v$. Denote by $\hat{G}_{S,\text{aut}}$ the support of $L^2(G(F)\backslash G(\nA))$ as a $G(F_S)$-representation in $\widehat{G(F_S)}$. For a split reductive subgroup $H\subset G$, we analogously define $\hat{H}_{S,\text{aut}}.$
	\begin{enumerate}[(1)]
		\item If $\pi\in \hat{H}_{S,\text{aut}}$, then the support of $\Ind_{H(F_S)}^{G(F_S)}\pi$ is contained in $\hat{G}_{S,\text{aut}}.$
		\item If $\pi\in \hat{G}_{S,\text{aut}}$, then the support of $\Res^{G(F_S)}_{H(F_S)}\pi$ is contained in $\hat{H}_{S,\text{aut}}.$
		\item If $\pi,\rho\in\hat{G}_{S,\text{aut}}$, then the support of $\pi\otimes\rho$ is contained in $\hat{G}_{S,\text{aut}}.$
	\end{enumerate}
\end{theorem} 
Clozel in \cite{Clozel2004} and \cite{Clozel2007} then established the following consequences of this result. For a representation $\pi$ of $G(F_v)$ which occurs in a local $A$-packet $A_{\psi_v}$ associated to a local $A$-parameter $\psi_v$, we say $\psi_v|_{SL_2(\nC)}$ is an $SL_2$-type of $\pi$. Clozel conjectures that if $\pi$ is unramified, then it has a unique $SL_2$-type $\cA$. From this conjecture, theorem \ref{BSthm}, and Arthur's description of $L^2(H(F)\backslash H(\nA))$, Clozel then shows that this implies all of the unramified representations occuring in $\pi|_{H(F_v)}$ have the same $SL_2$-type $\cB$. Moreover, the type $\cB$ depends only on $\cA$. This produces a map
\begin{equation}\{SL_2\text{-types of }G\}\xrightarrow{\Res^G_H}\{SL_2\text{-types of }H\}.\label{sl2map}\end{equation}
which is independent of the choice of place $v$. Similarly, inducing an Arthur type representation $\pi$ of $H(F_v)$ to $G(F_v)$ produces a map in the opposite direction, and the tensor product of two Arthur type $G(F_v)$-representations produces a binary operation on the set of $SL_2$-types of $G$.

In the case $G=GL_n$, the map $(\ref{sl2map})$ as well as the analogues for induction and tensor product were explicated by Venkatesh in \cite{Venkatesh2005}. In that case, an $SL_2$-type is an algebraic homomorphism $SL_2(\nC)\rightarrow \check{G}=GL_n(\nC)$, which is parameterized up to conjugacy by a partition of $n$. Thus the map $(\ref{sl2map})$ in the $GL_n$ restricted to $GL_m$ case can be expressed as a map 
\begin{equation}\label{partitionmap}\{\text{Partitions of }n\}\xrightarrow{\Res^{GL_n}_{GL_m}} \{\text{Partitions of }m\}.\end{equation}
To describe Venkatesh's result, we will use the following notation. Let $\lambda=(\lambda_i)$ be a sequence of integers such that $\sum_i\max\{0,\lambda_i\}\leq n.$ We denote by $\la \lambda\ra_n$ the partition of $n$ consisting of all positive $\lambda_i$'s together with enough $1$'s to sum to $n$. For example, if $\lambda=(-3,-1,0,3,5)$, then $\la \lambda\ra_{10}=\la 3,5,1,1\ra.$ If the size of $n$ is clear, we will omit it from the notation. We will commonly write such partitions in forms such as $\la \lambda_i-k\ra$, by which we mean the partition formed by subtracting $k$ from each $\lambda_i$ and then following the above process on the resulting sequence to obtain a partition of $n$. For example, the trivial representation of $GL_n$ has the $SL_2$-type $\la n\ra$, and tempered representations have a type of $\la 1,\dots,1\ra$ (it follows from part (1) of theorem \ref{BSthm} that all tempered representations occur in $A$-packets). Venkatesh's result is the following.
\begin{theorem}[Venkatesh]
	\begin{enumerate}[(1)]
		\item Induction from a Levi subgroup: Suppose $(m_i)$ is a partition of $n$, and let $\pi_i$ be a representation of $GL_{m_i}(F)$ of Arthur type with $SL_2$ type $\la \lambda_{i,j}\ra$ for each $i$. Let $M=\prod_i GL_{m_i}(F)$ be the corresponding Levi subgroup of $GL_n(F)$. Then $\Ind_M^{GL_n(F)}\otimes_i \pi_i$ has type $\la \lambda_{i,j}-(n-m_i)\ra$.
		\item Restriction to a Levi subgroup: Suppose $\pi$ is a representation of $GL_n(F)$ of Arthur type with type $\la \lambda_i\ra$, and suppose $GL_m(F)$ is embedded into $GL_n(F)$ via the standard embedding. Then $\pi|_{GL_m}$ has type $\la \lambda_i-(n-m)\ra.$
		\item Tensor product: Suppose $\pi_1$ and $\pi_2$ are representations of Arthur type with $SL_2$ types $\la \lambda_i\ra$ and $\la \tau_i\ra.$  Then $\pi_1\otimes \pi_2$ has type $\la \lambda_i+\tau_j-n\ra.$
	\end{enumerate}
\end{theorem}
Our goal is to extend this result to symplectic and orthogonal groups.

\subsection{The classical group case}
Let $G$ be $Sp_{2n}$, $SO_{2n+1}$, or $O_{2n}$. In this case,
$$\check{G}=\begin{cases}
SO_{2n+1}(\nC)&\text{if }G=Sp_{2n}\\
Sp_{2n}(\nC)&\text{if }G=SO_{2n+1}\\
O_{2n}(\nC)&\text{if }G=O_{2n}
\end{cases}.$$
An $SL_2$ type can once again be described by a partition. However, in this case, there are some restrictions on the partitions determined by the nilpotent orbits in $\check{G}$ \cite{Collingwood2017}. Namely,
\begin{itemize}
	\item An $SL_2$-type for $Sp_{2n}$ is a partition of $2n+1$ such that each even number occurs an even number of times,
	\item An $SL_2$-type for $SO_{2n+1}$ is a partition of $2n$ such that each odd number occur an even number of times,
	\item An $SL_2$-type for $O_{2n}$ is a partition of $2n$ such that each even number occurs an even number of times.
\end{itemize}
The even/odd numbers that are restricted in this way are often referred to as having ``bad parity."

As we mentioned previously, we restrict our attention to unramified representations since those will turn out to have unique $SL_2$ types (see proposition \ref{uniquenessoftypes}). We will also restrict our attention to symplectic and even orthogonal group so as not to complicate the proof too much. Our result is the following.
\begin{theorem}[Main results]
	\label{mainthm}
	Let $G_{2n}$ be $Sp_{2n}$ or $O_{2n}$, and let 
	$$\epsilon=\begin{cases}1&\text{if }G_{2n}=O_{2n}\\-1&\text{if }G_{2n}=Sp_{2n}\end{cases}.$$
	\begin{enumerate}[(1)]
		\item (Induction from $GL_n$ to $G_{2n}$) For an unramified representation $\sigma$ of a $GL_n$ Siegel Levi subgroup of any $SL_2$ type, the representation $\Ind_{GL_n}^{G_{2n}}\sigma$ is of type $\la 1,\dots,1\ra$ except in the case where $G_{2n}=O_{2n}$ and $\sigma$ is a one-dimensional character, in which case the type is $\la 2,\dots,2\ra$ if $n$ is even or $\la 1,1,2,\dots,2\ra$ if $n$ is odd.
		\item (Induction from $G_{2n}$ to $GL_{2n}$) For an unramified representation $\rho$ of $G_{2n}$ of any $SL_2$ type, the representation $\Ind_{G_{2n}}^{GL_{2n}}\rho$ is of type $\la 1,\dots,1\ra$, except in the case where $G_{2n}=Sp_{2n}$ and $\rho$ is one-dimensional, in which case the type is $\la 2,\dots,2\ra.$
		\item (Induction from $G_{2a}\times G_{2b}$ to $G_{2a+2b}$) Suppose $a\leq b$, and let $\rho_a$ and $\rho_b$ be unramified representations of $G_{2a}$ and $G_{2b}$ respectively. Suppose $\la \lambda_i\ra$ is the type of $\rho_b$. Then the type of $\Ind_{G_{2a}\times G_{2b}}^{G_{2a+2b}}\rho_a\boxtimes\rho_b$ is $\la \lambda_i-2a\ra$, except in the case where $\rho_a$ and $\rho_b$ are both one-dimensional characters and $G_{2n}=Sp_{2n}$, in which case the type is $\la 2b-2a+1,2,\dots,2\ra$.
		\item (Restriction from $G_{2n}$ to $GL_{n}$) Let $\rho$ be an unramified representation of $G_{2n}$ of type $\la \lambda_i\ra$, and let $GL_n$ be Siegel Levi subgroup. Then $\rho|_{GL_n}$ has type $\la \lambda_i-n+\epsilon\ra,$ except in the case where $G=O_{2n}$ and the type of $\rho$ is $\la \lambda_1,\lambda_2\ra$ with $1<\lambda_1\leq \lambda_2$, in which case the type is 
		$$\la \lambda_2-n+1,2,\dots,2\ra=\la n+1-\lambda_1,2,\dots,2\ra$$ 
		if $\lambda_1,\lambda_2$ are odd or $\la 2,\dots,2\ra$ if $\lambda_1=\lambda_2=n$ is even.
		\item (Restriction from $GL_{2n}$ to $G_{2n}$) Let $\sigma$ be an unramified representation of $GL_{2n}$ of type $\la \lambda_1,\dots,\lambda_k\ra$, and assume $\lambda_k$ is the largest number in the type. Then $\sigma|_{G_{2n}}$ has type $\la 2\lambda_k-2n-\epsilon,1,\dots,1\ra$ except when $k=2$ and $G_{2n}=Sp_{2n}$, in which case the type is $\la \lambda_2-\lambda_1+1,2,\dots,2\ra$ if $\lambda_1,\lambda_2$ are even or $\la \lambda_2-\lambda_1+1,1,1,2,\dots,2\ra$ if $\lambda_1,\lambda_2$ are odd.
		\item (Restriction from $G_{2n}$ to $G_{2m}$) Let $\rho$ be an unramified representation of $G_{2n}$ of type $\la \lambda_i\ra$, and suppose $m<n$. The type of $\rho|_{G_{2m}}$ is $\la \lambda_i+2m-2n\ra.$
		\item (Tensor product) Let $\rho_1$ and $\rho_2$ be unramified representations of $G_{2n}$ with types $\la \lambda_i\ra$ and $\la \tau_i\ra$ respectively. Then $\rho_1\otimes \rho_2$ has type $\la \lambda_i+\tau_j-2n+\epsilon\ra,$ except in the case where $G_{2n}=O_{2n}$, and the types of $\rho_1$ and $\rho_2$ are $\la \lambda_1,\lambda_2\ra$ and $\la \tau_1,\tau_2\ra$ with $1<\lambda_1\leq\lambda_2$, $1<\tau_1\leq\tau_2$, and $\lambda_2+\tau_2>2n-1$, in which case, assuming $\lambda_1\leq \tau_1$ (so that $\tau_1+\lambda_2-(2n-1)=\tau_1-\lambda_1+1>0$), the type is 
		$$\la \lambda_2+\tau_2-(2n-1),\lambda_2+\tau_1-(2n-1),2,\dots,2\ra$$ 
		if $\lambda_1,\lambda_2$ are odd or $\la 2,\dots,2\ra$ if $\lambda_1=\lambda_2=n$ is even.
	\end{enumerate}
\end{theorem}
The formulas for each case are remarkably similar to those from the $GL_n$ case, except for those curious exceptional cases (see remark \ref{exceptionalremark}). The differing $\epsilon$ for symplectic and orthogonal groups arises as a result of the type of the trivial representation, which is $\la 2n+1\ra$ for $Sp_{2n}$ and $\la 2n-1,1\ra$ for $O_{2n}$; hence one can verify that tensoring with the trivial representation results in the same type, and restricting the trivial representation gives the trivial type. One can also easily verify that the parity conditions on numbers in partitions are satisfied by these results.

Let us describe how to form an unramified representation that has a given $SL_2$ type. Let $\lambda$ be a valid $SL_2$ type for a $G$-representation of Arthur type, that is, a partition that satisfies the restriction on ``bad parity" entries described above. We may write
$$\lambda=\la \lambda_i\ra\oplus 2\la \tau_j\ra,$$
where $\oplus$ denotes concatenation, the $2$ means repeating those numbers twice, and we take all the $\lambda_i$ to be distinct (``good parity") numbers. For example, an unramified representation of $Sp_{22}$ of type $\la 1,2,2,3,5,5,5\ra$ can have its type written as $\la 1,3,5\ra\oplus 2\la 2,5\ra.$ Then let $G'$ be a group of the same type (symplectic or orthogonal) as $G$ with $\rank G'=\rank G-\sum \tau_j$. We form the Levi subgroup $M=G'\times \prod_j GL_{\tau_j}(\nC)$ corresponding to a parabolic $P$. Assuming for the moment that we may find an unramified representation $\rho$ of $G'$ of type $\la \lambda_i\ra$, we then choose unramified characters $\chi_j$ of $\GL_{\tau_j}$ for each $j$. Then the parabolically induced representation
\begin{equation}
\label{prestdrep}
\Ind_P^G \rho\otimes\bigotimes_j \chi_j
\end{equation}
is almost irreducible (it is either irreducible or decomposes into two direct summands), and it has a unique unramified subrepresentation of Arthur type $\lambda$. If $\sigma=\Ind_{\prod_j GL_{\tau_j}}^{GL_{\sum_j \tau_j}}\bigotimes_j \chi_j$ and $P'$ is the maximal parabolic with Levi subgroup $G'\times GL_{\sum_j \tau_j}$, then we may rewrite $(\ref{prestdrep})$ as the parabolically induced representation
\begin{equation}\label{stdrep}\Ind_{P'}^G\rho\otimes \sigma.\end{equation}

The representation $\rho$ can be taken to be an unramified unipotent representation, which for us means a representation obtained via iterated theta lifts from a one-dimensional unramified character with the theta lift being unramified at each step, which produces a representation of the expected type in the following way. Suppose $G$ is $O_{2n}$ or $Sp_{2n}$ and that $\lambda_1<\dots<\lambda_k$ are distinct odd positive integers whose sum is $2n-1$ or $2n+1$ for $G$ orthogonal or symplectic respectively.

If $k\not=1$ and $\lambda_1=1$, then a character $\chi$ of $O_{1+\lambda_2}$ has type $(1,\lambda_2)$. If $k=2$ we are done, but otherwise the theta lift of $\chi$ to $Sp_{\lambda_2+\lambda_3}$ has type $(1,\lambda_2,\lambda_3)$ if it is nonzero. If $k=3$ we are done, but otherwise the theta lift from there to $O_{1+\lambda_2+\lambda_3+\lambda_4}$ if nonzero has type $(1,\lambda_2,\lambda_3,\lambda_4)$. We proceed in this way until we reach $G$. 

If $k=1$ or $\lambda_1>1$, then a character $\chi$ of $Sp_{\lambda_1-1}$ has type $\lambda_1$. Iterated lifting to $O_{\lambda_1+\lambda_2}$ and so on up to $G$ will produce a representation of the desired type (in the case $k=1$ and $\lambda_1=1$, $G'$ is the trivial group, so $M$ has only the $GL_{\tau_j}$ factors and we need not worry about $\rho$). 

In the case where $G=SO_{2n+1}$, the iterated theta lifting proceeds similarly using metaplectic groups instead of symplectic groups in the dual pair correspondence. 

The use of a unipotent representation is necessary in this case because it's not possible to realize an Arthur type unramified representation with distinct odd numbers in its type as fully parabolically induced from a one-dimensional character of a Levi subgroup. This is a significant difference between the classical group case and the $GL_n$ case where unipotent representations are just unramified characters. 

We may analyze the restriction of an induced representation with Mackey theory. By finding an unramified representation which is weakly contained in the restriction of the induced representation, we will have identified the map (\ref{partitionmap}). This will not be sufficient for unipotent representations, which we will analyze separately using theta correspondence.

Similar to Venkatesh's proof for $GL_n$, our argument will be an inductive one, where for example the restriction problem may reduce to a lower rank tensor product problem. From the use of theta correspondence, we will also see that some symplectic cases reduce to even orthogonal cases and vice-versa. We focus only on maximal cases, such as induction from the Siegel Levi subgroup $GL_n$ for example, because the operations on types are transitive.

Although we do not solve the odd orthogonal case here, we shall state the corresponding results as a conjecture including passing between odd and even orthogonal groups. They are readily computable using the techniques of this paper.  We state them without any exceptional cases since we have not computed whether any should arise, and so some of the formulas below may fail in some exceptional cases.
\begin{conjecture}[Odd orthogonal groups]\label{oddorthogonalconjecture}
	\begin{enumerate}[(1)]
		\item (Induction from $GL_n$ to $O_{2n+1}$) For an unramified representation $\sigma$ of $GL_n$ of any $SL_2$ type, $\Ind_{GL_n}^{O_{2n+1}}\sigma$ is of type $\la 1,\dots,1\ra.$
		\item (Induction from $O_{2n+1}$ to $GL_{2n+1}$) For an unramified representation $\rho$ of $O_{2n+1}$ of any $SL_2$ type, $\Ind_{O_{2n+1}}^{GL_{2n+1}}\rho$ is of type $\la 1,\dots,1\ra.$
		\item (Induction from $O_a\times O_b$ to $O_{a+b}$) Assume $a\leq b$, and let $\rho_a$ and $\rho_b$ be unramified Arthur type representations of $O_a$ and $O_b$ respectively, where the type of $\rho_b$ is $\la \lambda_i\ra.$ Then $\Ind_{O_a\times O_b}^{O_{a+b}}\rho_a\boxtimes\rho_b$ has type $\la \lambda_i-a\ra.$
		\item (Restriction from $O_{2n+1}$ to $GL_{n}$) Let $\rho$ be an unramified Arthur type representation of $O_{2n+1}$ of type $\la\lambda_i\ra.$ Then $\rho|_{GL_n}$ has type $\la \lambda_i-n\ra$.
		\item (Restriction from $GL_{2n+1}$ to $O_{2n+1}$) Let $\sigma$ be an unramified representation of $O_{2n+1}$ of type $\la \lambda_1,\dots,\lambda_k\ra$, and assume $\lambda_k$ is the largest number in the partition. Then $\sigma|_{O_{2n+1}}$ has type $\la 2\lambda_k-2n-2,1,\dots,1\ra.$
		\item (Restriction from $O_n$ to $O_m$) Let $\rho$ be an unramified Arthur type representation of $O_n$ of type $\la \lambda_i\ra$. For any $m<n$, $\rho|_{O_m}$ has type $\la \lambda_i+m-n\ra.$
		\item (Tensor product for $O_{2n+1}$) Let $\rho$ and $\rho'$ be unramified Arthur type representations of $O_{2n+1}$ of types $\la \lambda_i\ra$ and $\la \tau_i\ra$ respectively. Then $\rho\otimes\rho'$ has type $\la \lambda_i+\tau_j-2n\ra.$
	\end{enumerate}
\end{conjecture}
Although it is included in the conjecture, part (2) is already known to be true (see remark \ref{gurevichremark}).
\subsection{Setting up the proof}
We are able to make some simplifications in order to make things easier. Since by the global argument of Clozel, the map of $SL_2$-types will be the same at every place, it suffices to solve the problem at just one place of $F$. Also, we are looking for the type in common among all representations in a direct integral decomposition, but that means we only need to find the type of just one unramified representation which occurs in the decomposition. We will work with $p$-adic fields and prove this weakened version of Theorem \ref{mainthm} in chapters 3, 4, and 5.
\begin{proposition}
\label{weakprop}
	For any $SL_2$ type $\la \lambda_i\ra$ for $G_{2n}$, there exists an unramified Arthur type representation of $G_{2n}(\nQ_p)$ of type $\la \lambda_i\ra$ whose restriction to $G_{2m}(\nQ_p)$ weakly contains an unramified Arthur type representation of type $\la \lambda_i+2m-2n\ra$. A similar assertion holds for each of the cases in theorem \ref{mainthm}.
\end{proposition}
We will show how theorem \ref{mainthm} follows from proposition \ref{weakprop} after proposition \ref{unramifiedspectrumofarthurtype}.

Here is a broad overview of the proof of proposition \ref{weakprop}. To simplify the notation, we give proofs only in the symplectic case with the understanding that the even orthogonal case uses almost exactly the same arguments. In cases where the orthogonal case does differ from the symplectic one, we give remarks at the end of those sections describing the differences and what happens as a result. These differences tend to arise either as a result of the exceptional cases of theorem \ref{mainthm} or from the fact that a skew-symmetric bilinear form on an odd-dimensional space is necessarily degenerate but the same is not true for a symmetric bilinear form.

The restriction and tensor product problems must be done in two parts, which we have divided into chapters 3 and 4. We begin in chapter 3 with the representation (\ref{stdrep}) of $G$ and consider its restriction using Mackey theory. This argument is fairly elementary as the geometry involved is very nice. The case of a tensor product of two representations, one of which is of the form (\ref{stdrep}), is analyzed using a tool we develop in section \ref{restrictiontoparabolicsection} using Howe's theory of low rank representations. 

The unipotent cases for restriction and tensor product are addressed in chapter 4 and use the explicit construction of the theta lift by Li \cite{Li1989}. We prove weak containment directly by approximating matrix coefficients. In each case we are exploiting a certain see-saw dual pair relationship. We are required to cope with some analytical difficulties and forced to make some dominated convergence arguments in this case. 

Chapter 5 deals with the induction problems. In the case for induction from a $GL_n$ to $Sp_{2n}$, we follow the technique of Venkatesh and induce through the Siegel parabolic, which allows us to use harmonic analysis on the abelian unipotent radical. For the remaining induction problems, we present a new type of deformation argument based in part on the boundary degeneration theory in $\cite{Sakellaridis2012}$ and the theory of deformations to the normal cone, though our calculations are very explicit and do not rely on any advanced machinery. 

Before we begin, let us establish some notation. If $P=MN$ is a parabolic subgroup of a group $G$, $\pi$ a representation of a subgroup $M'$ of $M$, and $\pi$ extends to a representation of $M'N$ with $N$ acting via a character $\chi$, we denote this representation of $M'N$ by $\pi\cdot \chi$. The trivial representation of a group $H$ will be denoted $\textbf{1}_H.$ To conserve space, we may write an induced representation from $H$ to $G$ as $\cI_H^G$ and the restriction of a representation $\pi$ from $G$ to $H$ as either $\pi|_H$ or $\cR^G_H \pi$. We use ``$\boxtimes$'' for the outer tensor product whenever necessary to distinguish it from the inner tensor product. Thus the parabolically induced representation (\ref{stdrep}) can be written as $\cI_{P'}^G (\rho\boxtimes\sigma\cdot\textbf{1}_{N'}).$

 If $\pi$ and $\pi'$ are two unitary representations of a group $G$, we will write $\pi\prec \pi'$ or $\pi'\succ \pi$ to indicate that $\pi$ is weakly contained in $\pi'$. All groups will be $p$-adic groups unless otherwise stated; we will write $Sp_{2m}$ for the symplectic group $Sp(2m,\nQ_p)$ and $O_{2m}$ for the split orthogonal group $O(2m,\nQ_p)$. 
 
For notational consistency, most of the time $\sigma$ will be used to denote representations of $GL_n$, $\rho$ will be used to denote unipotent representations of $Sp_{2n}$, $\pi$ will be used to denote an induced representation of $Sp_{2n}$, $\la \tau_i\ra$ will be used for the type of $\sigma$, and $\la \lambda_i\ra$ will be used for the type of $\rho$ or $\pi$. We also use $\gamma$ typically to denote a tempered representation.

\section{Preliminary results}
\subsection{Unramified representations of Arthur type}
Recall that we use the term \textit{Arthur type} to refer to irreducible unitary representations in local $A$-packets. We will prove some essential facts about unramified representations of Arthur type.
\begin{proposition}\label{uniquenessoftypes}
	Unramified representations of $Sp_{2n}(\nQ_p)$ and $O_n(\nQ_p)$ of Arthur type have unique $SL_2$ types.
\end{proposition}
\begin{proof}
	let $G$ be a symplectic or orthogonal group over $\nQ$. A local $A$-parameter for $G(\nQ_p)$ is a map
	$$\psi_p:W_{\nQ_p}\times SL_2(\nC)\times SL_2(\nC)\rightarrow \check{G},$$
	where $W_{\nQ_p}\times SL_2(\nC)$ is the Weil-Deligne group of $\nQ$. We refer to the first $SL_2(\nC)$ factor as the Deligne $SL_2(\nC)$ and the second one (the one that determines the $SL_2$ type) as the Arthur $SL_2(\nC)$. Proposition 6.4 of M\oe glin's paper \cite{Moeglin2009} states that the associated local $A$-packet $\Pi(\psi)$ contains an unramified representation if and only if $\psi|_{W_{\nQ_p}}$ is unramified and $\psi$ restricts trivially to the Deligne $SL_2(\nC)$, and when such an unramified representation does exist, it is the only unramified representation in $\Pi(\psi).$ Moreover Arthur has shown that each $A$-packet contains an associated $L$-packet \cite{Moeglin2011} with unramified L-parameter $\phi_{\psi_p}:W_{\nQ_p}\times SL_2(\nC)\rightarrow\check{G}$ defined by 
	$$\phi_{\psi_p}(w)=\psi_p\left(w,\begin{pmatrix}
	|w|^\frac{1}{2}&0\\0&|w|^\frac{-1}{2}
	\end{pmatrix}\right),$$
	where $|w|$ is the pullback to $W_{\nQ_p}\times SL_2(\nC)$ of the absolute value on $W_{\nQ_p}$, and so the unique unramified representation in the $A$-packet must in fact belong to this $L$-packet. The map $\psi_p\mapsto \phi_{\psi_p}$ is injective \cite{Arthur2013}, and so the unique unramified repesentation in an $A$-packet does not belong to any other $A$-packet. Therefore such an unramified representation has a unique $SL_2$ type.
\end{proof}
\begin{proposition}\label{unramifiedspectrumofarthurtype}
	Let $G$ be a general linear, orthogonal, or symplectic group and $H$ a general linear, orthogonal, or symplectic subgroup of $G$ over a local field $F_v$ of a number field $F$. If $\pi$ is an Arthur type representation of $G$, then every representation in the spectrum of $\pi|_H$ is of Arthur type.
\end{proposition}
\begin{proof}
	The Burger-Sarnak theorem (\ref{BSthm}) tells us the spectrum of $\pi|_H$ is contained in $L^2(H(F)\backslash H(\nA)).$ Arthur's book \cite{Arthur2013} tells us these representations are in Arthur packets.
\end{proof}
We are now able to show how theorem \ref{mainthm} follows from proposition \ref{weakprop}. We give a proof in the restriction $Sp_{2n}$ to $Sp_{2m}$ case. The other cases are similar. This argument is due to Clozel.
\begin{proof}[Proof of theorem \ref{mainthm} assuming proposition \ref{weakprop}]
	Let $F$ be a number field that has $p$ and $v$ among its places corresponding to local fields $\nQ_p$ and $F_v$ with $p\not=v$. Suppose $\pi$ is an unramified representation of $Sp_{2n}(F_v)$ in an Arthur packet for some place $v$. By proposition \ref{uniquenessoftypes}, $\pi$ does not occur in any other Arthur packet and has a unique $SL_2$-type $\la \lambda_i\ra.$ Suppose $\pi'$ is an unramified Arthur type representation of $Sp_{2n}(\nQ_p)$ for some place $p\not=v$ that has the same $SL_2$ type as $\pi$. 
	
	Let $m<n$, and suppose $\sigma$ and $\sigma'$ are two unramified representations in the spectrum of $\pi|_{Sp_{2m}(F_v)}.$ These representations are of Arthur type with unique $SL_2$ types by theorem \ref{BSthm}, proposition \ref{unramifiedspectrumofarthurtype}, and proposition \ref{uniquenessoftypes}. By proposition $\ref{weakprop}$, $\pi'|_{Sp_{2m}(\nQ_p)}$ weakly contains an unramified Arthur type representation $\rho$ of type $\la \lambda_i-2n+2m\ra.$ 
	
	By theorem $\ref{BSthm}$, since $\sigma\otimes \rho$ and $\sigma'\otimes \rho$ are weakly contained in $\pi\otimes\pi'|_{Sp_{2m}(F_v\times \nQ_p)}$, it is weakly contained in $L^2(Sp_{2m}(F)\backslash Sp_{2m}(\nA))$. Because they occur in $L^2(Sp_{2m}(F)\backslash Sp_{2m}(\nA))$, they are associated to some global $A$-parameter, and because each of the representations $\sigma,\sigma',$ and $\rho$ occur only in one local $A$-packet each, this means $\sigma$ and $\rho$ have the same $SL_2$ type and $\sigma'$ and $\rho$ have the same type. Then $\sigma$ and $\sigma'$ have the same type. Therefore the unramified spectrum of $\pi|_{Sp_{2m}(F_v)}$ all has the type $\la \lambda_i-2n+2m\ra$, and it follows that even ramified representations must have this type among their types. This concludes the proof.
\end{proof}
For the rest of this paper we are proving proposition \ref{weakprop}. We begin with some results which justify the construction (\ref{stdrep}).
\begin{proposition}\label{typeofthetalift}
	Let $\pi$ be an unramified Arthur type representation of $O_{2m}(\nQ_p)$ of type $\la \lambda_1,\dots,\lambda_k\ra$, and suppose $n\geq m$. Then the theta lift $\theta(\pi)$ of $\pi$ to $Sp_{2n}(\nQ_p)$ is nonzero, and in particular it is an unramified Arthur type representation of type $\la\lambda_1,\dots,\lambda_k,2n-2m+1\ra.$ The same is true when lifting from $Sp_{2m}$ to $O_{2n}$ but with the $+1$ in the formula replaced by $-1$ and under the assumption $n>m$.
\end{proposition}
\begin{proof}
	The unramified theta correspondence was described explicitly by Kudla and Rallis in proposition 7.1.1 of \cite{Kudla1994}. We will describe their result in the case of lifting from the split group $H=O_{2m}(\nQ_p)$ to $G=Sp_{2n}(\nQ_p)$ with $m\leq n$. Let $\pi$ be an irreducible unramified representation of $H=O_{2m}(\nQ_p)$ with associated $L$-parameter
	$$\phi_p:L_{\nQ_p}\rightarrow \check{H},$$
	where $L_{\nQ_p}$ is the Weil-Deligne group $W_{\nQ_p}\times SL_2(\nC)$, and $\check{H}=O_{2m}(\nC)$. Since $O_{2m}(\nQ_p)$ is the split orthogonal group, $\phi_p$ factors through $SO_{2m}(\nC)\subset \check{H}$. Suppose $n\geq m$, and consider the embedding 
	$$\iota:SO_{2m}(\nC)\times SO_{2n-2m+1}(\nC)\hookrightarrow SO_{2n+1}(\nC)=\check{G}.$$ 
	The result of Kudla and Rallis is that the theta lift $\theta(\pi)$ of $\pi$ to $G=Sp_{2n}(\nQ_p)$ is nonzero and unramified with $L$-parameter 
	$$\theta(\phi_p):L_{\nQ_p}\rightarrow \check{G}$$
	given by 
	$$\theta(\phi_p)(w)=\iota\left(\phi_p(w),\diag(|w|^{n-m},|w|^{n-m-1},\dots,|w|^{m-n+1},|w|^{m-n})\right).$$
	Now suppose $\pi$ is of Arthur type associated to a local $A$-parameter
	$$\psi_p:L_{\nQ_p}\times SL_2(\nC)\rightarrow \check{H}.$$
	Recall that the $L$-parameter in this case is given by
	$$\phi_{\psi_p}=\psi_p\left(w,\diag(|w|^\frac{1}{2},|w|^\frac{-1}{2})\right).$$
	Let $\tau:SL_2(\nC)\rightarrow SO_{2n-2m+1}(\nC)$ be the principal $SL_2$-type which determined on the diagonal torus by the map 
	$$\tau(\diag(a,a^{-1}))=\diag(a^{2(n-m)},a^{2(n-m-1)},\dots,a^{2(m-n+1)},a^{2(m-n)}).$$
	Then define the $A$-parameter $\Psi_p:L_{\nQ_p}\times SL_2(\nC)\rightarrow \check{G}$ by 
	$$\Psi_{p}(w,x)=\iota(\psi_p(w,x),\tau(x)).$$ 
	This unramified since $\psi_p$ is. The associated $L$-parameter is
	\begin{align*}
	\Phi_{\Psi_p}&=\Psi_p(w,\diag(|w|^\frac{1}{2},|w|^\frac{-1}{2}))\\
	&=\iota(\psi_p(w,\diag(|w|^\frac{1}{2},|w|^\frac{-1}{2})),\tau(\diag(|w|^\frac{1}{2},|w|^\frac{-1}{2}))\\
	&=\iota(\phi_{\psi_p}(w),\diag(|w|^{n-m},|w|^{n-m-1},\dots,|w|^{m-n+1},|w|^{m-n}))\\
	&=\theta(\phi_{\psi_p})(w).
	\end{align*}
	Therefore $\theta(\pi)$ belongs to the $A$-packet associated to $\Psi_p$. Its $SL_2$ type is the type of $\pi$ concatenated with the type $\tau$, which is $2n-2m+1$ as desired.
\end{proof}
\begin{corollary}\label{iteratedthetalifting}
	For each $SL_2$ type consisting of distinct odd numbers, there exists an unramified Arthur type representation of $Sp_{2n}$ or $O_{2n}$ with this type which is an iterated theta lift from the trivial representation of some smaller symplectic or even orthogonal group such that the theta lift is an unramified Arthur type representation at each step.
\end{corollary}
\begin{proof}
	 Proposition \ref{typeofthetalift} lets us theta lift repeatedly from a one-dimensional character (which we take to be trivial) to construct such a representation. In particular if $\lambda_1<\dots<\lambda_k$ are odd numbers, we may start with a character of type $\lambda_1$, then lift to a representation of type $\la \lambda_1,\lambda_2\ra$, and so on until we get a representation of type $\la \lambda_1,\dots,\lambda_k\ra.$
\end{proof}
From this corollary we obtain our unipotent representations. Then to obtain a representation of any type we need the following.
\begin{proposition}\label{parabolicarthurtype}
	Let $G_{2n}$ be $Sp_{2n}$ or $O_{2n}$. Suppose $\rho$ is an unramified Arthur type representation of $G_{2m}$ type $\la \lambda_i\ra$ with $m<n.$ Let $P=MN$ be a parabolic subgroup of $G_{2n}$ with Levi subgroup $G_{2m}\times GL_{\tau_1}\times\dots\times GL_{\tau_k}$, where $2m+2\sum_{k}\tau_k=2n.$ Suppose $\chi_i$ is an unramified character of $GL_{\tau_i}$ for each $i$. Then
	\begin{equation}\label{eqn211}\Ind_{P}^{G_{2n}}\rho\boxtimes\chi_1\boxtimes\dots\boxtimes \chi_k\cdot\textbf{1}_{N}\end{equation}
	has a unique unramified Arthur type subrepresentation of type $\la \lambda_i\ra\oplus 2\la \tau_i\ra.$
\end{proposition}
\begin{proof}
	The representation (\ref{eqn211}) has a unique irreducible unramified subrepresentation by the Iwasawa decomposition. It's of Arthur type by the work of M\oe glin \cite{Moeglin2011}.
\end{proof}
\subsection{The restriction of a unitary representation to a parabolic subgroup}\label{restrictiontoparabolicsection}
One of the difficulties we will encounter in the restriction and tensor product problems is when we restrict a representation to some non-reductive subgroup and do not have a trivial action on the unipotent radical. In order to read the type of a representation induced from a parabolic subgroup, a trivial action on the unipotent radical is necessary so that we have parabolic induction of a representation on the Levi subgroup. We will resolve much of that in this section by proving the following.
\begin{lemma}\label{restrictiontoparabolic}
	$G$ be $GL_n$, $O_{2n}$, or $Sp_{2n}$, and suppose $\pi$ is an irreducible unramified Arthur type representation of $G$. Let $P=MN$ be a parabolic subgroup of $G$. Then there exists an irreducible unramified Arthur type representation $\rho$ of $M$ such that $\pi|_P\succ\rho\cdot\textbf{1}_N.$
\end{lemma}
 We give a proof for $GL_n$ and one for $Sp_{2n}$, with the understanding that the $O_{2n}$ case is essentially the same as the $Sp_{2n}$ one. Observe that if we show that $\pi|_P$ weakly contains some $\pi'\cdot \textbf{1}_N$, where $\pi'$ is a representation of $M$, then continuity of restriction with respect to the Fell topology implies that $\pi|_M\succ \pi'$, which in turn means $\pi|_M\cdot\textbf{1}_N\succ \pi'\cdot\textbf{1}_N$. We also note that it suffices to prove the lemma for a restriction to maximal parabolic subgroup since restriction can be done in stages.
\begin{proof}[Proof ($GL_n$ case).]
	We will prove this by induction on the number of numbers in the type of $\pi$. Observe that if $\pi$ is a one-dimensional character, then since one-dimensional characters of $GL_n$ factor through the determininant, $\pi|_N$ is already trivial, and thus $\pi|_P=\pi|_M\cdot\textbf{1}_N.$
	
	Now suppose $\pi$ is an irreducible unramified Arthur type representation of $GL_n$ which is realized as $\pi=\Ind_{P'}\sigma_a\boxtimes\sigma_b\cdot\textbf{1}_{N'}$, where $\sigma_a$ is an unramified Arthur type representation of $GL_a$, $\sigma_b$ is a one-dimensional unramified character of $GL_b$, and $P'=M'N'$ is the maximal subgroup of type $(a,b)$ of $G$. The $SL_2$ type of $\pi$ is the type of $\sigma_a$ concatenated with $b$, and we may assume $b$ is the largest number in the type of $\pi$.
	
	We analyze the restriction to $P$ using Mackey theory. Note $P'$ acts on $G/P$ with a unique open orbit. If we assume $P$ is in this open orbit, then the stabilizer of the $P'$ action is $P'\cap P$, and Mackey theory says we will have
	$$\Res_P^G\Ind_{P'}^G\sigma_a\boxtimes\sigma_b\cdot\textbf{1}_{N'}\cong\Ind_{P\cap P'}^P\Res^{P'}_{P\cap P'}\sigma_a\boxtimes\sigma_b\cdot\textbf{1}_{N'}.$$
	
	We realize $P$ as the stabilizer of some $m$-dimensional subspace $W$ of an $n$-dimensional space $V$. Meanwhile $P'$ is the stabilizer of some $a$-dimensional subspace $W'$. Any elements of $P\cap P'$ stabilize $W\cap W'$ and so each possible dimension of $W\cap W'$ forms an orbit for the action of $P'$. The open orbit corresponds to the generic case where $W\cap W'$ have as small an intersection as possible. We take $\sigma_a$ to be a representation of $GL(W')$ while $\sigma_b$ is a representation of $GL(V/W')$.\\
	
	\textbf{Case 1:} Suppose $a+m\leq n.$ Then $W\cap W'$ should have trivial intersection. Choose a (possibly zero) complement $U$ to $W\oplus W'$ so that $V=W\oplus U\oplus W'.$ We then can realize $M$ and $M'$ as
	\begin{align*}M&=GL(W)\times GL(W'\oplus U)\\M'&=GL(W')\times GL(W\oplus U)\end{align*}
	so that $\sigma_a\in \widehat{GL(W')}$ and $\sigma_b\in \widehat{GL(W\oplus U)}.$ Any element stabilizing $W$ and $W'$ must stabilize $W\oplus W'$, so let $P_\oplus$ be the parabolic subgroup of $GL(V)$ of elements stabilizing $W\oplus W'.$ The Levi subgroup of $P_\oplus$ can be realized as $M_\oplus=GL(W\oplus W')\times GL(U)$, and $N_\oplus$ is the set of elements which act trivially on $W\oplus W'$ and $V/(W\oplus W')$. We have 
	$$N_\oplus\cong\Hom(U,W\oplus W')=\Hom(U,W)\times \Hom(U,W').$$
	Since $N_\oplus$ acts trivially on $W$ and $W'$, it stabilizes both of these spaces. Of course the $GL(U)$ factor of $M_\oplus$ stabilizes both spaces as well, and the elements of $GL(W\oplus W')$ stabilizing both $W$ and $W'$ is the Levi subgroup $GL(W)\times GL(W').$ Hence
	$$P\cap P'=(GL(W)\times GL(W')\times GL(U))N_\oplus.$$
	
	Of course the $\sigma_a$ is already a representation of $GL(W')$, and the $\sigma_b\in\widehat{GL(W\oplus U)}$ action restricts trivially to $GL(W')$ and to $N_\oplus$ since $\sigma_b$ is a one-dimensional character. Since $\sigma_b$ is one dimensional, we may write $\sigma_b|_{GL(W)\times GL(U)}=\sigma_b|_{GL(W)}\boxtimes \sigma_b|_{GL(U)}.$ Thus we have
	$$\pi|_P\cong\Ind_{P\cap P'}^P\sigma_a\boxtimes \sigma_b|_{GL(W)}\boxtimes \sigma_b|_{GL(U)}\cdot\textbf{1}_{N_\oplus}.$$
	We will induce in stages. According to our descriptions so far, we have 
	$$N=\Hom(W'\oplus U,W)=\Hom(W',W)\oplus \Hom(U,W).$$
	Also $P''=(GL(W')\times GL(U))\Hom(U,W')$ forms a parabolic subgroup of $GL(W'\oplus U)$. By inducing in stages and using proposition E.2.5 of \cite{Bekka2008}, we have
	\begin{align*}\pi|_P\cong \Ind_{(GL(W')\times P'')N}^{P}&\left(\sigma_b|_{GL(W)}\boxtimes\left(\sigma_a\boxtimes\sigma_b|_{GL(U)}\cdot \textbf{1}_{\Hom(U,W')}\right)\cdot \textbf{1}_N\right)\\&\otimes L^2\left((GL(W')\times P'')N/(GL(W')\times P'')\Hom(U,W)\right).\end{align*}
	Let us consider the representation 
	$$L^2\left((GL(W')\times P'')N/(GL(W')\times P'')\Hom(U,W)\right).$$ 
	We will analyze this in a manner similar to Venkatesh in section 3.3 of \cite{Venkatesh2005}. As a vector space, 
	$$L^2\left((GL(W')\times P'')N/(GL(W')\times P'')\Hom(U,W)\right)\cong L^2(\Hom(W',W)).$$
	The $(GL(W')\times P'')\Hom(U,W)$ action on this space is intertwined by the Fourier transform, so we may realize this representation on the space $L^2(\widehat{\Hom(W',W)}).$ The orbits of the $(GL(W')\times P'')\Hom(U,W)$-action on this space correspond to the different possible ranks of maps $W\rightarrow W'$. The unique open orbit corresponds to maps of maximal rank. If $\zeta\in \widehat{\Hom(W',W)}$ is an element of the open $(GL(W')\times P'')\Hom(U,W)$ orbit and $M_\zeta\subset (GL(W')\times P'')\Hom(U,W)$ its stabilizer, then Fourier analysis produces an isomorphism
	$$L^2\left((GL(W')\times P'')N/(GL(W')\times P'')\Hom(U,W)\right)\cong \Ind_{M_\zeta N}^{(GL(W')\times P'')N}\textbf{1}_{M_\zeta}\cdot \zeta.$$
	If we replace $\zeta$ by $\zeta_t(x)=\zeta(tx)$, then the stabilizer $M_{\zeta_t}$ is the same as $M_\zeta.$ The map $t\mapsto \zeta_t$ is continuous on the Fell topology of $\widehat{\Hom(W',W)}$. By continuity of induction with respect to the Fell topology, we may let $t\rightarrow 0$ to obtain 
	\begin{align*}L^2&\left((GL(W')\times P'')N/(GL(W')\times P'')\Hom(U,W)\right)\\&\hspace{3cm}\succ L^2((GL(W')\times P'')\Hom(U,W)/M_\zeta)\cdot\textbf{1}_{\Hom(W',W)}.\end{align*}
	Observe that $\Hom(U,W')$ and $GL(U)$ act trivially on $\Hom(W',W)$, and so we will have $GL(U)\Hom(U,W')\subset M_\zeta$. Then $M_\zeta$ will be of the form $(M_\zeta'\times GL(U))\Hom(U,W')$, where $M_\zeta'$ is the subgroup of elements of $GL(W)\times GL(W')$ stabilizing $\zeta.$ The representation $L^2((GL(W)\times GL(W'))/M_\zeta')$ was described in \cite{Venkatesh2005} section 3.3. It weakly contains a representation of the form $\tau_W\boxtimes \tau_{W'}$, where $\tau_W$ and $\tau_{W'}$ are irreducible Arthur type representations of $GL(W)$ and $GL(W')$, which we may take to be unramified. Hence 
	$$L^2\left((GL(W')\times P'')N/(GL(W')\times P'')\Hom(U,W)\right)\succ (\tau_W\boxtimes\tau_{W'}\boxtimes\textbf{1}_{GL(U)})\cdot\textbf{1}_N.$$
	Therefore we now have
$$\pi|_{P}\succ \cI_{(GL(W')\times P'')N}^P(\sigma_b|_{GL(W)}\otimes\tau_{W})\boxtimes\left((\sigma_a\otimes \tau_{W'})\boxtimes\sigma_b|_{GL(U)}\cdot\textbf{1}_{\Hom(U,W')}\right)\cdot\textbf{1}_N.$$
Now if we write 
$$\pi'=\Ind_{P''}^{GL(W'\oplus U)}(\sigma_a\otimes \tau_{W'})\boxtimes\sigma_b|_{GL(U)}\cdot\textbf{1}_{\Hom(U,W')},$$
then we have 
$$\pi|_{P}\succ (\sigma_b|_{GL(W)}\otimes\tau_{W})\boxtimes \pi'\cdot\textbf{1}_N.$$
The $M$-representation $(\sigma_b|_{GL(W)}\otimes\tau_{W})\boxtimes \pi'$ weakly contains an irreducible unramified Arthur type representation $\rho$, which proves the lemma in this case.
	
	\textbf{Case 2:} Suppose $a+m>n$. Then $W$ and $W'$ have nontrivial intersection, and in the open orbit case we will have $V=W+W'.$ Any elements stabilizing both spaces must stabilize the intersection. Let $P_\cap=M_\cap N_\cap$ be the parabolic subgroup of $GL(V)$ consisting of elements stabilizing $W\cap W'$. Choose complements $U$ and $U'$ in $W$ and $W'$ respectively to $W\cap W'$ so that 
	$$V=U\oplus (W\cap W')\oplus U'.$$
	We may then take as our Levi subgroups
	\begin{align*}
	M_\cap&=GL(W\cap W')\times GL(U\oplus U'),\\
	M&=GL(W)\times GL(U'),\\
	M'&=GL(W')\times GL(U)
	\end{align*}
	so that $\sigma_a\in \widehat{GL(W')}$ and $\sigma_b\in \widehat{GL(U)}.$ The elements of $N_\cap$ act trivially on $W\cap W'$ and $V/(W\cap W')$, and so they stabilize both $W$ and $W'$. The elements of $GL(U\oplus U')$ that stabilize both $W$ and $W'$ stabilize both $U$ and $U'$, and so we see that
	$$P\cap P'=(GL(U)\times GL(W\cap W')\times GL(U'))N_\cap.$$
	Of course $\sigma_a$ restricts trivially to $GL(U)$, and $\sigma_b$ is already a representation of $GL(U)$. Then we have
	$$\pi|_P=\Ind_{P\cap P'}^{P}\sigma_a|_{(GL(W\cap W')\times GL(U'))N_\cap}\otimes \sigma_b.$$
	Note that 
	$$N_\cap=\Hom(U\oplus U',W\cap W')=\Hom(U,W\cap W')\times \Hom(U',W\cap W').$$ 
	Of course $\sigma_a$ restricts trivially to the $\Hom(U,W\cap W')$ factor since that factor acts trivially on $W'$. We see that 
	$$(GL(W\cap W')\times GL(U'))\Hom(U',W\cap W')$$ 
	is a parabolic subgroup of $GL(W').$ Applying the induction hypothesis, we know that $\sigma_a|_{(GL(W\cap W')\times GL(U'))\Hom(U',W\cap W')}$ weakly contains a representation of the form $\sigma'\cdot \textbf{1}_{\Hom(U',W\cap W')}$, where $\sigma'$ is an irreducible Arthur type unramified representation of $GL(U')\times GL(W\cap W')$. Since $\sigma'$ is weakly contained in $\sigma_a$ restricted to a Levi subgroup, we know that we can take it to be of the form $\sigma_{W\cap W'}\boxtimes \sigma_{U'}.$ Hence 
	$$\pi|_{P}\succ\Ind_{P\cap P'}^P\sigma_{U'}\boxtimes \sigma_b\boxtimes (\sigma_{W\cap W'}\cdot \textbf{1}_{N_\cap}).$$
	Next observe that $P''=(GL(W\cap W')\times GL(U))\Hom(U,W\cap W')$ is a parabolic subgroup of $GL(W)$. By inducing in stages through $(GL(U')\times P'')N$ and repeating the same harmonic analysis as in the previous case, we find that for some unramified Arthur type representations $\tau_U$ and $\tau_{U'}$ of $GL(U)$ and $GL(U')$ that
	$$\pi|_{P}\succ \Ind_{(GL(U')\times P'')N}^P(\sigma_{U'}\otimes \tau_{U'})\boxtimes\left(\sigma_{W\cap W'}\boxtimes(\sigma_b\otimes\tau_U)\cdot\textbf{1}_{\Hom(U,W\cap W')}\right)\cdot\textbf{1}_N.$$
	Then let 
	$$\pi'=\Ind_{(GL(W\cap W')\times GL(U))\Hom(U,W\cap W')}^{GL(W)}\sigma_{W\cap W'}\boxtimes(\sigma_b\otimes\tau_U)\cdot \textbf{1}_{\Hom(U,W\cap W')}.$$
	By inducing in stages, we have
	$$\pi|_P\succ (\sigma_{U'}\otimes \tau_{U'})\boxtimes \pi'\cdot\textbf{1}_{N}.$$
	The $M$-representation $(\sigma_{U'}\otimes \tau_{U'})\boxtimes \pi'$ weakly contains an unramified Arthur type representation. This completes the proof in the $GL_n$ case.\phantom\qedhere
\end{proof}

The proof for the classical group case will be shorter and rely on a wholly different technique. We use Howe's theory of low rank representations \cite{Howe2010}. In principle, the theory of low rank representations exists for $GL_n$ representations as well and could perhaps provide a shorter proof of case 1 as well, but this theory has not been explicated for $GL_n$ in the literature.
\begin{proof}[Proof ($Sp_{2n}$ case):]
	Let $\pi$ be an irreducible unramified Arthur representation of $Sp_{2n}$. Let $P$ be the maximal parabolic subgroup of $Sp_{2n}$ with Levi subgroup $M=GL_m\times Sp_{2n-2m}$ and unipotent radical $N$.
	
	\textbf{Case 1:} Suppose the rank of $\pi$ in the sense of Howe (which is even by \cite{Howe2010}) is $2a\leq m.$ Then according to \cite{Li1989}, $\pi$ is a theta lift of a representation $\rho$ of $O_{2a}$, and $\pi|_{Z(N)}$ is supported on rank $2a$ characters, where $Z(N)$ is the center of $N$. The stabilizer in $M$ of a rank $2a$ character is a subgroup $(O_{2a}\times GL_{m-2a})N'\times Sp_{2n-2m}$, where $N'$ is the unipotent radical of a parabolic subgroup $P'$ of $GL_m$ with Levi subgroup $M'=GL_{2a}\times GL_{m-2a}.$ If $\chi$ is a rank $2a$ character of $Z(N)$ and $p_\chi$ denotes the (unique up to conjugation) representation of $N$ with central character $\chi$, then we have an oscillator representation $\omega$ of $(O_{2a}\times Sp_{2n-2m})N$ whose $N$-action is by $p_\chi.$ The restriction $\pi|_P$ is irreducible and by the arguments of Li \cite{Li1989} has the description
	\begin{equation}\label{eqn221}\pi|_P\cong \Ind_{((O_{2a}\times GL_{m-2a})N'\times Sp_{2n-2m})N}^P\rho\otimes \omega\boxtimes \textbf{1}_{GL_{m-2a}}\cdot\textbf{1}_{N'}.\end{equation}
	(In the inducing data, $O_{2a}$ acts by $\rho\otimes \omega$, $Sp_{2n-2m}$ and $N$ act by $\omega.$) If we replace $\chi$ by $\chi_t(x)=\chi(t^2x)$, the isomorphism class of the corresponding oscillator representation $\omega_t$ does not change as $t$ varies. Hence the isomorphism class of the representation (\ref{eqn221}) does not change when replacing $\chi$ by $\chi_t$. If we let $t\rightarrow 0$, then $\chi_t$ approaches the trivial representation. The closure in $\hat{N}$ of the set of $p_t$, $t\not=0$ is $\{p_t:t\not=0\}\cup \widehat{N/Z(N)}$, so we may let $p_t$ approach the trivial representation continuously in the Fell topology on $N$ as $t\rightarrow 0$. Since induction is continuous with respect to the Fell topology, this tells us
	$$\pi|_P\succ \left(\Ind_{(O_{2a}\times GL_{m-2a})N'\times Sp_{2n-2m}}^M\rho\otimes\omega|_{O_{2a}\times Sp_{2n-2m}}\cdot \textbf{1}_{GL_{m-2a}N'}\right)\cdot\textbf{1}_N.$$
	The $M$-representation $\Ind_{(O_{2a}\times GL_{m-2a})N'\times Sp_{2n-2m}}^M\rho\otimes\omega|_{O_{2a}\times Sp_{2n-2m}}\cdot \textbf{1}_{GL_{m-2a}N'}$ has nonzero unramified spectrum supported on Arthur type representations by theorem \ref{BSthm} and proposition \ref{parabolicarthurtype}.
	
	\textbf{Case 2:} Suppose the rank of $\pi$ in the sense of Howe is greater than $m$. Then $\pi|_{Z(N)}$ is supported on maximal rank characters. The $M$-stabilizer of such a character is $O_m\times Sp_{2n-2m}.$ Any irreducible representation of $P$ weakly contained in $\pi|_P$ must be of the form
	$$\pi'\otimes\Ind_{(O_m\times Sp_{2n-2m})}^{M}\rho\otimes \omega,$$
	where $\pi'$ is some representation of $M$, $\rho$ is some representation of $O_m\times Sp_{2n-2m}$, and $\omega$ is an oscillator representation for $(O_m\times Sp_{2n-2m})N$ whose $Z(N)$ action is by a full rank character. Once again we may deform such a character to get a trivial action on $M$ as desired. This concludes the proof.
\end{proof}
\subsection{Decay of matrix coefficients}\label{sectiondecayofmatrixcoefficients}
Our main tool for the unipotent case will be theta correspondence. An explicit construction for the theta lift was given in \cite{Li1989} in the stable range. This was later extended to the so called ``semistable" range \cite{Hongyu2000}. For the dual pair $(O_{2m},Sp_{2n})$ with $m<n$, the construction of the lift of an $O_{2m}$-representation $\sigma$ to $Sp_{2n}$ relies on the absolute convergence of the integral
\begin{equation}\label{Liintegral}\int_{O_{2m}}\la\omega(gh)\phi,\phi'\ra\la\sigma(h)v,v'\ra dh,\end{equation}
where $(\omega, V_\omega)$ is a Weil representation for the dual pair $(O_{2m},Sp_{2n})$, $g\in Sp_{2n}$, $\phi,\phi'\in V_\omega$, and $\la\sigma(h)v,v'\ra$ is a matrix coefficient for $\sigma.$ That requires certain $L^p$ decay of the matrix coefficients of our representations in order that the matrix coefficient integrals in his construction converge. 

The $SL_2$-type of a representation is expected to determine a bound for the decay of its $K$-finite matrix coefficients. In particular, this bound should depend on the rank of the group, the largest number in the partition, and that number's multiplicity in the partition. We will establish the following.
\begin{theorem}\label{convergencetheorem}
Let $(H,G)$ be one of the dual pairs $(GL_m,GL_n)$, $(O_{m},Sp_{n})$, or $(Sp_{m},O_{n})$ ($m,n$ even in the latter two cases) over $\nQ$ with $m\leq n$, and suppose $\sigma$ is an irreducible unitary representation of $H(\nQ_p)$ occurring in some local $A$-packet. Suppose $\lambda$ is the largest number in the partition associated to an $SL_2$-type of $\sigma$ and that all numbers in the type of $\sigma$ greater than $2$ (if any) are distinct.

Let $\omega$ be a Weil representation for the dual pair $(H,G)$ and $V_\omega$ and $V_\sigma$ the spaces on which $\omega$ and $\sigma$ act respectively. Then for all $\phi,\phi'\in V_\omega$ and $v,v'\in V_\sigma$,
\begin{equation*}\int_{H(\nQ_p)} \la \omega(h)\phi,\phi'\ra\la\sigma(h)v,v'\ra dh\end{equation*}
converges absolutely if $m+\lambda\leq n$ when $H=GL_m$ or $H=O_{m},$ or if $m+\lambda< n$ when $H=Sp_{m}.$
\end{theorem}
Note that we do not assume uniqueness of $SL_2$ types, and so this theorem applies to representations of Arthur type which are not unramified as well; this is a crucial point, as we will apply lemma 6.2(b) of \cite{Howe2010} to representations $\pi$ whose spectrum have a common $SL_2$-type in order to obtain a matrix coefficient bound for $\pi$ itself. 

Theorem \ref{convergencetheorem} was established by Li in the stable range $(n\geq 2m)$ in his original paper \cite{Li1989} and by others in some cases since: for example in the real unipotent case in \cite{Ma2017}. Ideally, we would like to establish this theorem without the restriction on the type of $\sigma$, but the proof is computation heavy, and it may even be the case that the theorem fails without some such restrictions.

The integral (\ref{Liintegral}) defines a (possibly degenerate) $G$-invariant inner product $(\cdot,\cdot)$ on $V_\omega\otimes V_\sigma$. It follows from theorem A.5 of \cite{Harris2011} that $(\cdot,\cdot)$ is positive semi-definite. Hence if $R$ denotes the kernel of $(\cdot,\cdot)$, then $(\cdot,\cdot)$ descends to give a positive-definite inner product on $(V_\omega\otimes V_\sigma)/R$. In particular, if $(\cdot,\cdot)$ is nondegenerate, then the theta lift of $\sigma$ is nonzero and $\theta(\sigma)$ is realized on $(V_\omega\otimes V_\sigma)/R$. We will refer to this as Li's construction of the theta lift.

\begin{proof}We give a proof in the $(H,G)=(O_{2m},Sp_{2n})$ case with the understanding that the other cases are similar. The decay of the matrix coefficients of $\omega$ is well understood. It follows from the proof of theorem 3.2 in \cite{Li1989} and (7.7) in \cite{Howe2010} that we have the estimate
\begin{equation}\label{eqn231}|\la\omega(h)\phi,\phi'\ra|\leq C(\phi,\phi')\Xi_{H}(h)^{\frac{n}{m-1}}\end{equation}
for $h\in H$ on $K$-finite vectors $\phi,\phi'\in V_\omega$, where $C(\phi,\phi')$ is some constant and $\Xi_H$ is Harish-Chandra's $\Xi$ function on $H$. Thus $\omega|_H$ has matrix coefficients in $L^p$ for $p>\frac{2m-2}{n}.$ If $2m-2<n$, then the boundedness of the matrix coefficients of $\sigma$ implies the absolute convergence of $(\ref{Liintegral})$ already. 

If $2m-2\geq n$, then $\sigma$ having $K$-finite matrix coefficients in $L^q$ (i.e. $\sigma$ being a strongly $L^q$ representation) for some $q<\frac{2m-2}{2m-2-n}$ by H\"{o}lder's inequality implies the convergence of $(\ref{Liintegral})$ (cf. corollary 3.4 of \cite{Li1989}). Given the assumption $2m+\lambda\leq 2n$ and $2m-2\geq n$ (so $m\geq 2$), theorem \ref{convergencetheorem} is implied by the following proposition, for which we also state the symplectic case for later use.\phantom\qedhere\end{proof}
\begin{proposition}\label{matrixcoefficientsdecayrate}
	Suppose $\sigma$ is a unitary irreducible Arthur type representation of $G=O_{2m}$ or $G=Sp_{2m}$ with nontrivial $SL_2$ type, $m\geq 2$, and that $\lambda$ is the largest entry of the partition associated to an $SL_2$ type of $\sigma$. Suppose that all numbers in the type of $\sigma$ greater than $2$ (if any) are distinct (and hence odd). Then there exists $\delta>0$ such that the $K$-finite matrix coefficients of $\sigma$ are bounded by a multiple of $\Xi_G^{\frac{2}{q}}$, where
	$$q\geq\begin{cases}\frac{2m-2}{m-3+\delta}&G=O_{2m}\text{ and }\lambda=2\\
	\frac{2(2m-2)}{2m-5-\lambda+\delta}&G=O_{2m}\text{ and }\lambda\text{ is odd}\\
	\frac{2m}{m-2+\delta}&G=Sp_{2m}\text{ and }\lambda=2\\
	\frac{4m}{2m-1-\lambda+\delta}&G=Sp_{2m}\text{ and }\lambda\text{ is odd}\end{cases}.$$
	A similar result holds for $G=GL_n.$
\end{proposition}
For example, when $\lambda=1$ and $G=O_{2m}$ so that $\sigma$ has type $(1,\dots,1)$, we expect this to be true for all $\delta<4$, which implies $\sigma$ is tempered. This result is stronger than we are actually able to show because we do not know the Ramanujan conjecture. Instead, estimates for the Ramanujan conjecture for classical groups (see \cite{Cogdell2004}) will give a smaller $\delta$ value but one which is still positive. In fact, the trivial estimate for the Ramanujan conjecture is sufficient for this application (see remark \ref{ramanujanestimate} for an example). For our proof, we will make a computation assuming the Ramanujan conjecture for ease of computation, but this will not affect the validity of the proof in general because the estimates known for the conjecture are sufficient to establish a positive value for $\delta.$
\begin{proof}
We will just do the $G=Sp_{2m}$ case. We use the theory of leading exponents to compute the $L^p$ decay of $\sigma$. The connection between exponents and decay of matrix coefficients is a result of Casselman (see \cite{Knapp1986} for the archimedean case and \cite{Casselman1995} for the nonarchimedean case). 

In \cite{Moeglin2009}, M\oe glin describes the exponents of representations in $p$-adic $A$-packets. Her description of the exponents differs from Casselman's theory, but they are related. For M\oe glin, a representation $\sigma$ occurs as the unique irreducible quotient of a parabolically induced representation $\Ind_P^G \chi\otimes\gamma$, where $\gamma$ is tempered and $\chi$ is a positive character. The character $\chi$, which is uniquely determined, will be denoted $Exp(\sigma)$. From lemma 5.2.8 in \cite{Wallach1988} (which is also valid in the nonarchimedean case), the matrix coefficients of $\sigma$ are bounded in terms of $\chi$ which gives a particular $L^p$ bound, and from Casselman's theory, the leading exponents of $\sigma$ must give the same $L^p$ bound as well. We will describe how to translate M\oe glin's $Exp(\sigma)$ data into leading exponents from which we will compute an $L^p$ bound.

M\oe glin showed that the exponents $Exp(\sigma)$ of a representation $\sigma$ as a character of a maximal torus are entry-wise bounded by a certain set of values, which she denotes $Exp(\psi)$ where $\psi$ is the $A$-parameter associated to $\sigma$. It follows from her results that the $A$-parameters which maximize $Exp(\psi)$ are the unramified $A$-parameters, and in that case, the $Exp(\psi)=Exp(\sigma)$ if $\sigma$ is the unique unramified representation in the $A$-packet associated to $\psi$. Moreover, the leading exponents of this representation depend only on the $SL_2$-type of $\sigma$. Hence we may assume that $\sigma$ is precisely this unramified representation with the largest exponents.

 Let $\la \lambda_1,\dots,\lambda_k\ra$ be the $SL_2$ type of $\sigma$ with $\lambda_1\leq\dots\leq\lambda_k=\lambda$. Let $Exp(\psi)=Exp(\sigma)$ be the sequence formed by concatenating the sequences 
 $$\left[\frac{\lambda_i-1}{2}, \frac{\lambda_i-3}{2},\dots,-\frac{\lambda_i-1}{2}\right],$$ 
 re-ordering the entries of the new sequence in decreasing order, and then removing all but the first $n$ entries. For example, if the $SL_2$-type is $\la 2,2,7\ra$ (so that $\sigma$ is a representation of $Sp_{10}$), then $Exp(\sigma)=\left(3,2,1,\frac{1}{2},\frac{1}{2}\right).$ Let us in general write $Exp(\sigma)=(a_1,\dots,a_m)$ with $a_i\geq a_{i+1}$. 

Recall that for $Sp_{2m}$, the half-sum of positive roots is $\rho=me_1+(m-1)e_2+\dots+e_m$, which we rewrite as the sequence $(m,m-1,\dots,1).$ The dominant weights here are $\omega_i=e_1+\dots+e_i$ for $i=1,\dots,m$, and so we have
\begin{align*}
\la\rho,\omega_i\ra&=m+(m-1)+\dots+(m-(i-1))=im-\frac{i(i-1)}{2},\\
\la Exp(\sigma),\omega_i\ra&=a_1+\dots+a_i.
\end{align*}
It follows from M\oe glin's result and the theory of exponents (\`a la Casselman) that $\sigma$ has $K$-finite matrix coefficients bounded by a multiple of $\Xi_G^\frac{2}{p}$ if and only if
$$p>\frac{2}{1-\frac{\la Exp(\sigma),\omega_i\ra}{\la \rho,\omega_i\ra}}$$
for all $i.$ Thus the problem of computing the $L^p$ decay for the matrix coefficients of $\sigma$ has been reduced to a purely combinatorial one.

Let us consider the most simple case first. Suppose that $\lambda=1$ so that $Exp(\sigma)=(0,\dots,0).$ Then it follows from the above that $\sigma$ has $K$-finite matrix coefficients in $L^p$ for $p>2$ and is bounded by $\Xi_G$ (recall that in actuality, there is a slight error due to the Ramanujan conjecture estimates, so we are not actually saying $\sigma$ is tempered).

Now suppose $\lambda=2.$ Then $Exp(\sigma)=\left(\frac{1}{2},\dots,\frac{1}{2},0,\dots,0\right)$, where the number of $\frac{1}{2}$ entries corresponds to the number of $2$'s in the type of $\sigma.$ Suppose there are $2k$ such $2$'s (recall that there must be an even number of them, and note $2k\leq m$). Then for $i\leq 2k$ we have 
$$\frac{\la Exp(\sigma),\omega_i\ra}{\la \rho,\omega_i\ra}=\frac{i}{2im-i(i-1)}=\frac{1}{2m-i+1}.$$
This expression is maximized in the range $1\leq i\leq 2k$ when $i=2k$, in which case we obtain
$$\frac{\la Exp(\sigma),\omega_{2k}\ra}{\la \rho,\omega_{2k}\ra}=\frac{1}{2m-2k+1}.$$
It is clear that $\frac{\la Exp(\sigma),\omega_{2k}\ra}{\la \rho,\omega_{2k}\ra}$ is maximal among the values $\frac{\la Exp(\sigma),\omega_i\ra}{\la \rho,\omega_i\ra}$ since for all values $i>2k$, the $i$th entry of $Exp(\sigma)$ is $0$, and so in that case we have
$$\frac{\la Exp(\sigma),\omega_i\ra}{\la \rho,\omega_i\ra}=\frac{\la Exp(\sigma),\omega_{2k}\ra}{\la \rho,\omega_{2k}\ra+x}$$
for some $x>0.$ Then $\sigma$ has $K$-finite matrix coefficients in $L^p$ for 
$$p>\frac{2m-2k+1}{m-k}.$$
The maximum possible value for this expression happens if $2k=m$, in which case we obtain $p>\frac{2(m+1)}{m}.$ It is easy to check that this value is strictly less than the value of $\frac{2m}{m-2}$ given in the statement of proposition \ref{matrixcoefficientsdecayrate}.

Now suppose $\lambda>2$. The computation for $Exp(\sigma)$ relies only on the $SL_2$ type of $\sigma$, and so we may choose any $\sigma$ with that type so long as it is the element of the associated $L$-packet with the maximal exponents. In particular it follows from proposition \ref{typeofthetalift} that we may take $\sigma$ to be an iterated theta lift from an unramified representation whose type consists only of $1$'s and $2$'s. We will give an inductive argument to show this has the desired growth.

Suppose $\sigma=\theta(\rho)$ is a theta lift of a representation $\rho$ of $O_{2n}$ where $2m-2n+1=\lambda$ and that proposition \ref{matrixcoefficientsdecayrate} (and hence theorem \ref{convergencetheorem}) holds for $\rho$. Then $\sigma$ may be given by Li's construction of the theta lift. In particular, its matrix coefficients are of the form 
$$g\mapsto\int_{O_{2n}}\la\omega(gh)\phi,\phi'\ra\la\rho(h)v,v'\ra dh.$$
We may generalize the bound (\ref{eqn231}) by considering both factors of the dual pair $(O_{2n},Sp_{2m})$. We may obtain the bound (cf. \cite{Ma2017}, lemma 3.8)
$$|\la\omega(gh)\phi,\phi'\ra|\leq C(\phi,\phi')\Xi_G(g)^{\frac{n}{m}}\Xi_{H}(h)^{\frac{m}{n-1}},$$
where $G=Sp_{2m}$ and $H=O_{2n}.$ Then
$$\left|\int_{H}\la\omega(gh)\phi,\phi'\ra\la\rho(h)v,v'\ra dh\right|\leq C(\phi,\phi')\Xi_G(g)^\frac{m}{n}\int_{H}\Xi_H(h)^\frac{m}{n-1}|\la\rho(h)v,v'\ra|dh.$$
We have already established the convergence of the integral on the right by the induction hypothesis. Then $\sigma$ has matrix coefficients bounded by $\Xi_G^\frac{2}{p}$, where 
$$p=\frac{2m}{n}=\frac{4m}{2m-\lambda+1}.$$
Note that this implies
$$\frac{4m}{2m-\lambda+1}>\frac{2}{1-\frac{\la Exp(\sigma),\omega_i\ra}{\la \rho,\omega_i\ra}},$$
and so because $Exp(\sigma)$ bounds $Exp(\pi)$ for any $\pi$ with the same $SL_2$-type, this bound applies (with some small error due to the Ramanujan conjecture estimates) to all representations of this type. This concludes the proof.
\end{proof}
\begin{remark}\label{ramanujanestimate}
	For the sake of illustration, let us see how the trivial estimate for the Ramanujan conjecture applies in the case $\lambda=2$. We must add an error factor of $\epsilon$ to each entry of $Exp(\sigma)$, from which we obtain
	$$\frac{\la Exp(\sigma)+(\epsilon,\dots,\epsilon),\omega_i\ra}{\la \rho,\omega_i\ra}=\frac{1+2\epsilon}{2m-i+1}.$$
	Taking the maximal value of $i=m$, we obtain a $p$ value of
	$$p>\frac{2(m+1)}{m-2\epsilon}$$
	This is less than $\frac{2m}{m-2}$ if and only if $\epsilon<\frac{1}{2}+\frac{1}{m}.$ The trivial estimate for the Ramanujan conjecture is $\epsilon<\frac{1}{2}$, and so adding this error term indeed does not affect the proof.
	
	Moreover, if in theorem \ref{convergencetheorem} we allow for the case $\lambda=2$, $m+2=n$ when $G=Sp_{2n}$, then the corresponding case of proposition \ref{matrixcoefficientsdecayrate} would require  $q\geq \frac{2m}{m-1+\delta}$ for some $\delta>0$. We have $\frac{2(m+1)}{m-2\epsilon}<\frac{2m}{m-1}$ if and only if $\epsilon<\frac{1}{2m}$, and so the truth of theorem 2.3.2 in the case $m+2=n$, $G=Sp_{2n}$ would follow if we knew the Ramanujan conjecture. Fortunately, we do not require this case in our application.
\end{remark}
It should be noted from the proof of proposition \ref{matrixcoefficientsdecayrate} that when we do have a representation which is a theta lift given by Li's construction that we obtain a bound for the matrix coefficients which does not rely on any estimates to the Ramanujan conjecture. This will be useful later, so let us state it as a corollary.
\begin{corollary}\label{matrixcoefficientsdecayratethetalifted}
Let $\rho$ be an Arthur type representation of $H=O_{2m}$ (or $H=Sp_{2m}$), and suppose the theta lift $\theta(\rho)$ of $\rho$ to $G=Sp_{2n}$ (or $G=O_{2n}$) is nonzero and may be given by Li's construction. Then $\theta(\rho)$ has $K$-finite matrix coefficients in $L^p$ for all 
$$p>\frac{2n-2\epsilon}{m},$$
where $\epsilon=0$ if $H$ is orthogonal or $\epsilon=1$ if $H$ is symplectic.
\end{corollary}

An essential fact for us relates the convergence of the integral $\ref{Liintegral}$ to the nonvanishing of the theta lift for unramified representations, which we now state.
\begin{corollary}
Let $(H,G)$ be one of the dual pairs $(GL_m,GL_n)$, $(O_{m},Sp_{n})$, or $(Sp_{m},O_{n})$ ($m,n$ even in the latter two cases) over $\nQ$ with $m\leq n$, and suppose $\sigma$ is a unitary irreducible representation of $H(\nQ_p)$ occurring in some local Arthur packet. Suppose $\lambda$ is the largest number in the partition associated to an $SL_2$-type of $\sigma$, and suppose that all numbers in the type of $\sigma$ greater than $2$ (if any) are distinct (and hence odd). Then the theta lift $\theta(\sigma)$ of $\sigma$ to $G$ is nonzero and may be constructed via Li's construction if $m+\lambda\leq n.$
\end{corollary}
\begin{proof}
Since $\sigma$ is unramified, proposition \ref{typeofthetalift} says its abstract theta lift is nonzero. By proposition 11.6 of \cite{Gan2014a}, when the integral (\ref{Liintegral}) converges, the abstract theta lift of $\sigma$ is nonzero if and only if the integral (\ref{Liintegral}) is nonvanishing. Therefore we may use Li's construction for the theta lift of $\sigma.$
\end{proof}
In particular, the following cases will be particularly useful to us.
\begin{corollary}\label{convergencecorollary}
	\begin{enumerate}[(i)]
		\item Suppose $\rho=\theta(\sigma)$ is an unramified Arthur type representation of $O_{2n}$ or $Sp_{2n}$ of type $(\lambda_1,\dots,\lambda_k)$ with $\lambda_1<\dots<\lambda_k$ which is a theta lift from an unramified Arthur type representation $\sigma$ of type $(\lambda_1,\dots,\lambda_{k-1})$. Then the integral $(\ref{Liintegral})$ converges absolutely, $\rho$ may be constructed using Li's construction, and its matrix coefficients have the bound established in corollary \ref{matrixcoefficientsdecayratethetalifted}.
		\item Suppose $a\in \{1,2\}$ and $\rho=\theta(\sigma)$ is an unramified Arthur type representation of $O_{2n}$ or  $Sp_{2n}$ of type $\la a,\dots,a,\lambda_1,\dots,\lambda_k\ra$ with $\lambda_1<\dots<\lambda_k$ and $a\in\{1,2\}$ which is a theta lift of a representation $\sigma$ of type $\la a,\dots,a,\lambda_1,\dots,\lambda_{k-1}\ra$ (or just $\la a,\dots,a\ra$ with $\lambda_1>a$ if $k=1$). Then $\rho$ may be constructed using Li's construction, and its matrix coefficients have the bound established in corollary \ref{matrixcoefficientsdecayratethetalifted}.
	\end{enumerate}
\end{corollary}
\section{Restriction and tensor product part 1: induced representations}
Our goal in this chapter is to reduce the restriction problem from that of a general induced representation of Arthur type to the case of a unipotent representation and to show how the tensor product problem reduces to the restriction one. That is, we will prove the restriction cases of proposition \ref{weakprop} for restrictions of representations of $Sp_{2n}$ using the induction hypothesis that all cases of proposition \ref{weakprop} hold for representations of $Sp_{2m}$, $m<n$, but this proof will only be for representations whose $SL_2$ types are not made entirely of distinct odd numbers. 

We will prove the tensor product case of proposition \ref{weakprop} for representations of $Sp_{2n}$ using the induction hypothesis that the restriction cases of proposition \ref{weakprop} hold for representations of $Sp_{2n}$ and that the tensor product result of proposition \ref{weakprop} holds for representations of $Sp_{2m}$ for $m<n.$ This proof will only work for the tensor product of two representation such that at least one of them does not have only distinct odd numbers it its type.

In chapter 4, we prove separately the unipotent cases, where a representation has only distinct odd numbers in its type.

The main tool used in this chapter is Mackey theory. Our proofs in the restriction cases will take the following form.
\begin{enumerate}[(1)]
	\item \label{stabilizerfind} In each section, we will have a group $H\subset G$ acting on $G/P$, where $P$ is a maximal parabolic subgroup of $G$. We realize $G/P$ as the Grassmannian of isotropic subspaces of a particular dimension. The group $H$ will act with open orbit. We will compute the stabilizer $H'$ in $H$ of a point in the open orbit in $G/P$ and assume the identity coset $P$ belongs to the open orbit (so that $H'=H\cap P$). For a representation $\omega$ of $P$, Mackey theory then says
	$$\Res_H^G\Ind_P^G\omega\succ \Ind_{H'}^H\Res^P_{H'}\omega.$$
	\item \label{stabilizerrestrict} We will then describe explicitly the restriction of $\omega$ to $H'$.
	\item \label{stabilizerinduce} We will induce in stages through a parabolic subgroup $P'$ of $H$ in order to obtain a representation of the form (\ref{stdrep}), from which we can read the type of the representation.
\end{enumerate}
The tensor product case will be much easier as it is almost a direct consequence of lemma \ref{restrictiontoparabolic}.
\subsection{Restriction from $Sp_{2n}$ to $Sp_{2m}$}\label{restrictionsphericalinduced}
In this section, we prove the following sub-case of proposition \ref{weakprop}.
\begin{lemma}\label{restrictionsphericallemma}
	Let $M=GL_a\times Sp_{2b}$ be a Levi subgroup of a parabolic subgroup $P$ of $G=Sp_{2n}(\nQ_p)$, and suppose $\pi=\Ind_P^G\sigma\boxtimes\rho\cdot\textbf{1}$ is an unramified Arthur type representation of the form (\ref{stdrep}), where $\sigma$ is an unramified Arthur type representation of $GL_a$ and $\rho$ is a unipotent unramified Arthur type representation of $Sp_{2b}.$ If $\la \lambda_i\ra$ is the type of $\pi$, then $\pi|_{Sp_{2m}(\nQ_p)}$ weakly contains an irreducible unramified Arthur type representation of type $\la\lambda_i+2n-2m\ra.$
\end{lemma}
We will show by using Mackey theory that this lemma reduces to the restriction problem for the unipotent $\rho$ and the case of induction from an $Sp_{2c}\times Sp_{2(m-c)}$ subgroup to $Sp_{2m}$. In other words, our induction hypotheses for this section is to use assume the restriction result for representations of $Sp_{2b}$ and the induction result for $Sp_{2c}\times Sp_{2(m-c)}$ to $Sp_{2m}$ from which we prove the result for $Sp_{2n}.$ In the orthogonal case, we also use the induction hypothesis for tensor product (see remark \ref{restrictionorthogonalremark})

\begin{proof} We realize $G/P$ as the isotropic Grassmannian of $a$-dimensional isotropic subspaces in a $2n$ dimensional vector space $V$ with skew-symmetric bilinear form $B$. Let $H=Sp(S)\times Sp(T)$ be the stabilizer of $V=S\oplus T$ where $S$ and $T$ are symplectic subspaces of dimension $2m$ and $2(n-m)$ respectively. We assume $m\leq n-m.$ We will consider three cases:
\begin{enumerate}[(1)]
	\item $a\leq2m\leq2(n-m)$ and $a$ is even,
	\item $a<2m\leq 2(n-m)$ and $a$ is odd, and
	\item $2m<a\leq 2(n-m).$
\end{enumerate}
In each case, we will show the restriction of $\pi$ to 	each factor of $H$ weakly contains a representation of the form $(\ref{stdrep})$, and we will demonstrate that it is of the predicted type $\la \lambda_i-2m+2n\ra$ assuming the induction hypothesis holds.
\\

\textbf{Case 1:} Suppose $a\leq 2m\leq 2(n-m)$ and $a$ is even. Let $W$ be an $a$-dimensional isotropic subspace of $V$. In the most generic case, the projections $W_S$ and $W_T$ of $W$ to $S$ and $T$ respectively are $a$-dimensional subspaces. The space $W\subset W_S\oplus W_T$ defines an isomorphism $\phi:W_T\rightarrow W_S$ with $\graph(\phi)=W$, and since $W$ is isotropic and $B=B|_S+B|_T$, it must be that $B(\phi(v),\phi(w))=-B(v,w)$ for all $v,w\in W_T.$ Hence $B|_{W_S}$ and $B|_{W_T}$ have kernels of equal dimension. 

The action of $Sp_{2m}=Sp(S)$ on the collection of all $a$-dimensional subspaces of $S$ with kernel of dimension $k$ is transitive, and so each possible dimension for the kernel of $B|_{W_S}$ corresponds to an $H$ orbit in $G/P.$ The largest orbit corresponds to the most non-degenerate case, i.e. the case where $W_S$ is a non-degenerate symplectic space of dimension $a$ (since $a$ is even). Thus a point in the open $H$-orbit of $G/P$ can be described by giving two $a$-dimensional symplectic subspaces of $S$ and $T$ together with an isomorphism $\phi$ between them that satisfies $B(\phi(v),\phi(w))=-B(v,w).$

\hyperref[stabilizerfind]{\textbf{The stabilizer:}} Let us describe the stabilizer in $H$ of a generic point in the open $H$-orbit of $G/P$. Suppose $P$ is the stabilizer of such a subspace $W$ so that the identity coset in $G/P$ belongs to the open $H$-orbit on $G/P$. We have that $S=W_S\oplus W_S^\perp$ and $T=W_T\oplus W_T^\perp,$ where $\perp$ denotes the symplectic complement in the spaces $S$ and $T$ respectively. Thus
$$V=S\oplus T=W_S\oplus W_S^\perp\oplus W_T\oplus W_T^\perp,$$
and 
$$W=\{w+\phi(w):w\in W_T\}\subset W_S\oplus W_T.$$

If $(g_1,g_2)\in Sp(S)\times Sp(T)$ stabilizes $W$, then for $(v,w)\in W_S\oplus W_T=W$ we have $(g_1v,g_2w)\in W$, which implies $g_1v\in W_S$ and $g_2w\in W_T$ for all $v\in W_S$ and $w\in W_T.$ Hence any elements in $H$ which stabilize $W$ must also stabilize $W_S$ and $W_T.$ Let $H_S$ and $H_T$ be the subgroups of elements of $Sp(S)$ and $Sp(T)$ which stabilize $W_S$ and $W_T$ respectively. Any elements of $Sp(S)$ or $Sp(T)$ stabilizing $W_S$ or $W_T$ must also stabilize $W_S^\perp$ or $W_T^\perp$ respectively. Thus
$$H\supset H_S\times H_T=Sp(W_S)\times Sp(W_S^\perp)\times Sp(W_T)\times Sp(W_T^\perp).$$

Now let $H_W$ be the subgroup of elements of $H_S\times H_T$ which stabilize $W=\graph\phi$, i.e those $(g_1,g_2)$ such that $g_1\phi(w)=\phi(g_2w).$ The isomorphism $\phi$ induces a symplectic form $(\cdot,\cdot)$ on $W=\graph\phi$ via 
$$(w_1+\phi(w_1),w_2+\phi(w_2)):=B(w_1,w_2)=-B(\phi(w_1),\phi(w_2))$$ 
for $w_1,w_2\in T$, and the set of elements of $Sp(W_S)\times Sp(W_T)$ preserving this form is a diagonally embedded subgroup $Sp(W)\subset Sp(W_S)\times Sp(W_T).$ Hence we have 
$$H_W=Sp(W)\times Sp(W_S^\perp)\times Sp(W_T^\perp).$$ 
The subgroup of elements of $H$ stabilizing the point $P$ in the open orbit of $G/P$ is then seen to be $H\cap P=H_W.$ By Mackey theory,
$$\pi|_{H}\succ\Ind_{H_W}^{H}\Res^P_{H_W}\sigma\boxtimes\rho \cdot \textbf{1}_N.$$

\hyperref[stabilizerrestrict]{\textbf{Restriction to the stabilizer:}} Let us describe $\Res^P_{H_W}\sigma\boxtimes\rho\cdot \textbf{1}_N$. We may choose the Levi subgroup $M$ of $P$ to be 
$$M=GL(W)\times Sp(W_S^\perp\oplus W_T^\perp).$$
With this choice, since $Sp(W)\subset GL(W)$ and $Sp(W_S^\perp)\times Sp(W_T^\perp)\subset Sp(W_S^\perp\oplus W_T^\perp)$, $\sigma$ restricts trivially to $Sp(W_S^\perp)\times Sp(W_T^\perp)$ while $\rho$ restricts trivially to $Sp(W)$. Inducing in stages, we have 
\begin{align*}
\pi|_{H}&\succ\Ind_{H_S\times H_T}^{H}\Ind_{H_W}^{H_S\times H_T}\sigma|_{Sp(W)}\boxtimes \rho|_{Sp(W_S^\perp)\times Sp(W_T^\perp)}\\
&=\Ind_{H_S\times H_T}^{H}\left(\Ind_{Sp(W)}^{Sp(W_S)\times Sp(W_T)}\sigma|_{Sp(W)}\right)\boxtimes \rho|_{Sp(W_S^\perp)\times Sp(W_T^\perp)}.
\end{align*}

\hyperref[stabilizerinduce]{\textbf{Induction in stages:}} If $\sigma'$ is any irreducible unitary representation weakly contained in $\sigma|_{Sp(W)}$, then by corollary 1 of \cite{Venkatesh2005}, $\Ind_{Sp(W)}^{Sp(W_S)\times Sp(W_T)}\sigma'$ is supported on tempered representations. By continuity of induction with respect to the Fell topology, that means $\Ind_{Sp(W_S)}^{Sp(W_S)\times Sp(W_T)}\sigma|_{W_S}$ weakly contains a representation $\gamma_S\boxtimes \gamma_T$, where $\gamma_S$ and $\gamma_T$ are tempered unramified representations of $Sp(W_S)$ and $Sp(W_T)$ respectively. Also, $\rho|_{Sp(W_S^\perp)\times Sp(W_T^\perp)}$ weakly contains a representation $\rho_S\boxtimes \rho_T$, where $\rho_S$ and $\rho_T$ are irreducible unramified Arthur type representations of $Sp(W_S^\perp)$ and $Sp(W_T^\perp)$. Splitting the induction into the $S$ and $T$ parts, we therefore have
\begin{equation}\label{eqn311}\pi|_{H}\succ\left(\Ind_{H_S}^{Sp(S)}\gamma_S\boxtimes \rho_S\right)\boxtimes\left(\Ind_{H_T}^{Sp(T)}\gamma_T\boxtimes\rho_T\right).\end{equation}
Applying the induction hypotheses for restriction to get the type of $\rho_S$ and induction from $Sp(W_S)\times Sp(W_S^\perp)$ to $Sp(S)$ (similarly for $T$), we can compute the type of the representation (\ref{eqn311}). Restricting $\rho$ to $Sp(W_S^\perp)=Sp_{2m-2a}$ adds $2m-a-2b$ to its type, and inducing further to $Sp(S)$ subtracts $a$ from its type again resulting in a type $\la\lambda_i+2m-(2a+2b)\ra=\la \lambda_i+2m-2n\ra$, which proves lemma \ref{restrictionsphericallemma} in this case.
\\

\textbf{Case 2:} Suppose $a< 2m\leq 2(n-m)$ and $a$ is odd. Let $W$ be an $a$-dimensional isotropic subspace of $V=S\oplus T$. In the most generic case, the projections $W_S$ and $W_T$ of $W$ to $S$ and $T$ respectively are $a$-dimensional subspaces, and we have an isomorphism $\phi:W_T\rightarrow W_S$ with $\graph(\phi)=W$ with $B(\phi(v),\phi(w))=-B(v,w)$ as before. Since $a$ is odd, the open $H$ orbit corresponds to those subspaces $W_S$ and $W_T$ with kernels of minimal dimension (in this case $1$). Thus we assume $B|_{W_S}$ and $B|_{W_T}$ have kernels $K_S$ and $K_T$ of dimension $1$, and hence $\phi(K_S)=K_T.$

Suppose $P$ is the stabilizer of such a space $W$ so that the identity coset in $G/P$ belongs to the open $H$ orbit. Choose $1$-dimensional isotropic spaces $K_S'$ and $K_T'$ which complement $K_S$ and $K_T$ so that $W_S\oplus K_S'$ and $W_T\oplus K_T'$ are nondegenerate $a+1$-dimensional spaces. Let 
\begin{align*}
S'&=(W_S\oplus K_S')^\perp,\\ 
T'&=(W_T\oplus K_T')^\perp,\\ 
U_S&=(K_S\oplus K_S'\oplus S')^\perp\subset W_S,\\
U_T&=(K_T\oplus K_T'\oplus T')^\perp\subset W_T
\end{align*}
where once again $\perp$ denotes the symplectic complement either in $S$ or in $T$ for those spaces as appropriate. We have 
$$S=U_S \oplus K_S\oplus K_S'\oplus S'\text{ and }T=U_T\oplus K_T\oplus K_T'\oplus T'.$$
By choosing $K_S'$ and $K_T'$ appropriately, we have $\phi(U_T)=U_S.$ 

\hyperref[stabilizerfind]{\textbf{The stabilizer:}} As in the previous case, any elements of $H$ that stabilize $W$ must stabilize $W_S$ and $W_T$. Since elements of $H$ preserve the bilinear form, they must stabilize $K_S$ and $K_T$ as well. Let $P_S=M_SN_S$ and $P_T=M_TN_T$ be the parabolic subgroups of $Sp(S)$ and $Sp(T)$ stabilizing $K_S$ and $K_T$ respectively. Then $P_S\times P_T$ is a parabolic subgroup $H$, and we have Levi subgroups
\begin{align*}M_S&= GL(K_S)\times Sp(U_S\oplus S'),\\
M_T&= GL(K_T)\times Sp(U_T\oplus T').\end{align*}
and the unipotent radical $N_S\times N_T$ of $P_S\times P_T$ acts by the identity endomorphism on $K_S, K_T, (U_S\oplus K_S\oplus S')/K_S, (U_T\oplus K_T\oplus T')/K_T, S/(U_S\oplus K_S\oplus S')$, and $T/(U_T\oplus K_T\oplus T').$ It follows that $N_S\times N_T$ stabilizes $W_S$ and $W_T$.

Let $M_H$ be the subgroup of elements of $M_S\times M_T$ which preserve $\phi.$ The isomorphism $\phi|_{U_T}:U_T\rightarrow U_S$ induces a symplectic form on the $a-1$-dimensional space $U:=\graph \phi|_{U_T}\subset W$, and the set of elements of $Sp(U_S)\times Sp(U_T)$ preserving $\phi$ form a diagonally embedded subgroup $Sp(U).$ If $K:=\graph \phi|_{K_T}\subset W$, then the set of elements of $GL(K_S)\times GL(K_T)$ preserving $\phi$ is a diagonally embedded $GL(K)$ subgroup. Since $Sp(S')\times Sp(T')$ acts by the identity on $W_S$ and $W_T$, we have
$$M_H=Sp(U)\times GL(K)\times Sp(S')\times Sp(T')$$

Let $N_H$ be the subgroup of elements of $N_S\times N_T$ which preserve the isomorphism $\phi.$ The center $Z(N_S\times N_T)=Z(N_S)\times Z(N_T)$ of $N_S\times N_T$ acts by the identity on $W_S$ and $W_T$, and so $Z(N_S\times N_T)\subset N_H.$ Since $Z(N_S\times N_T)$ acts trivially on $W_S\times W_T$ and $N_S\times N_T$ preserves $W_S\times W_T$, the action of $N_S\times N_T$ on $W_S\times W_T$ factors through $(N_S\times N_T)/Z(N_S\times N_T).$ We have 
\begin{align*}(N_S\times N_T)/&Z(N_S\times N_T)\cong \Hom(U_S\oplus S',K_S)\times \Hom(U_T\oplus T',K_T)\\
&=\Hom(U_S,K_S)\times \Hom(S',K_S)\times \Hom(U_T,K_T)\times \Hom(T',K_T).\end{align*}
The $\Hom(S',K_S)$ and $\Hom(T',K_T)$ factors act by the identity on $W_S$ and $W_T$ and so preserve $\phi$, but the elements of $\Hom(U_S,K_S)\times \Hom(U_T,K_T)$ that preserve $\phi$ are a diagonally embedded subgroup $N_\Delta.$ We have
\begin{align*}Z(N_H)&=Z(N_S)\times  Z(N_T),\\N_H/Z(N_H)&=\Hom(S',K_S)\times \Hom(T',K_T)\times N_\Delta.\end{align*}
The stabilizer of the identity coset in $G/P$ is then seen to be $M_HN_H$, and we have
$$\pi|_{M'}\succ\Ind_{M_HN_H}^{M'}\Res^P_{M_HN_H}\sigma\boxtimes\rho\otimes\textbf{1}_N.$$

\hyperref[stabilizerrestrict]{\textbf{Restriction to the stabilizer:}} We analyze this restriction. We may choose
$$M=GL(W)\times Sp(K_S'\oplus K_T'\oplus S'\oplus T')=GL(W)\times Sp(U_S^\perp\oplus U_T^\perp).$$
Note that $Z(N_H)$ acts by the identity on $W$, and so 
$$Z(N_H)\subset Sp(K_S'\oplus K_T'\oplus S'\oplus T')N.$$ 
Since $M=P/N$, we may write
$$(Sp(U)\times GL(K))N_\Delta\subset GL(W).$$
We realize $(Sp(U)\times GL(K))N_\Delta$ as a subgroup of the parabolic subgroup of $GL(W)$ of elements stabilizing $K$ with Levi subgroup $GL(U)\times GL(K).$ Also, we may write
$$Sp(S')\Hom(S',K_S)\times Sp(T')\Hom(T',K_T)\subset Sp(K_S'\oplus K_T'\oplus S'\oplus T')N/Z(N_H).$$
Hence we see $\sigma$ restricts trivially to $Sp(S')$, $Sp(T')$, $Z(N_H)$, $\Hom(S',K_S)$, and $\Hom(T',K_T)$, and $\rho$ restricts trivially to  $Sp(U)$, $GL(K)$, and $N_\Delta$. Then we have
\begin{equation}\label{eqn312}\pi|_{H}\succ \Ind_{M_HN_H}^{H}\sigma|_{(Sp(U)\times GL(K))N_\Delta}\boxtimes \rho|_{(Sp(S')\times Sp(T'))(N_H\cap Sp(U_S^\perp\oplus U_T^\perp))}.\end{equation}
By lemma \ref{restrictiontoparabolic} and by continuity of restriction with respect to the Fell topology, $\sigma|_{(Sp(U)\times GL(K))N_\Delta}$ weakly contains a representation of the form $\tau\cdot\textbf{1}_{N_\Delta}$, where $\tau$ is an irreducible unramified Arthur type representation of $Sp(U)\times GL(K)$.

Similarly, deforming the $\rho$ action of $N_H\cap Sp(U_S^\perp\oplus U_T^\perp))$ to the trivial representation as in the proof of lemma \ref{restrictiontoparabolic} tells us $\rho|_{(Sp(S')\times Sp(T'))(N_H\cap Sp(U_S^\perp\oplus U_T^\perp))}$ weakly contains a representation of the form $\rho_S\boxtimes\rho_T\cdot\textbf{1}_{N_H\cap Sp(U_S^\perp\oplus U_T^\perp)}$ such that $\rho_S$ and $\rho_T$ are irreducible unramified Arthur type representations of $Sp(S')$ and $Sp(T')$. Therefore
\begin{equation}\label{eqn313}\pi|_{H}\succ\Ind_{M_HN_H}^H\tau\boxtimes\rho_S\boxtimes \rho_T\cdot \textbf{1}_{N_H}.\end{equation}

\hyperref[stabilizerinduce]{\textbf{Induction in stages:}} Let us induce (\ref{eqn313}) in stages from $M_HN_H$ to $M_H(N_S\times N_T).$ We may rewrite (\ref{eqn313}) as 
$$\cI_{M_H(N_S\times N_T)}^H\tau\boxtimes\rho_S\boxtimes\rho_T\cdot \textbf{1}_{N_S\times N_T}\otimes L^2(M_H(N_S\times N_T)/M_HN_H).$$
Note that 
$$M_H(N_S\times N_T)/M_HN_H\cong (\Hom(U_S,K_S)\times \Hom(U_T,K_T))/N_\Delta.$$ 
Note that $L^2((\Hom(U_S,K_S)\times \Hom(U_T,K_T))/N_\Delta)\cong L^2(N_\Delta)$ as a $\Hom(U_S,K_S)\times \Hom(U_T,K_T)$ representation. We see  that 
\begin{align*}L^2((\Hom(U_S,K_S)\times \Hom(U_T,K_T))/N_\Delta)|_{\Hom(U_S,K_S)}&\cong L^2(\Hom(U_S,K_S))\\&\succ \textbf{1}_{\Hom(U_S,K_S)}\end{align*}
since $\Hom(U_S,K_S)$ is abelian, and this is similar on the $T$ side. Thus 
$$L^2((\Hom(U_S,K_S)\times \Hom(U_T,K_T))/N_\Delta)\cong L^2(N_\Delta)$$ 
weakly contains the trivial representation of $\Hom(U_S,K_S)\times \Hom(U_T,K_T).$ It follows that we now have
$$\pi|_H\succ \Ind_{M_H(N_S\times N_T)}^{H}\tau\boxtimes\rho_S\boxtimes\rho_T\cdot\textbf{1}_{N_S\times N_T}.$$

Once again by corollary 1 of \cite{Venkatesh2005}, $\Ind_{Sp(U)\times GL(K)}^{Sp(U_S)\times GL(K_S)\times Sp(U_T)\times GL(K_T)}\tau$ weakly contains a representation of the form $\gamma_{U_S}\boxtimes \gamma_{K_S}\boxtimes \gamma_{U_T}\times \gamma_{K_T}$, where each of the $\gamma$ representations are tempered unramified representations of their respective groups. Then
\begin{align*}\pi|_H&\succ\left(\Ind_{(GL(K_S)\times Sp(U_S)\times Sp(S'))N_S}^{Sp(S)}\gamma_{K_S}\boxtimes\gamma_{U_S}\boxtimes\rho_S\cdot\textbf{1}_{N_S}\right)\\&\hspace{3.5cm}\boxtimes\left(\Ind_{(GL(K_T)\times Sp(U_T)\times Sp(T'))}^{Sp(T)}\gamma_{K_T}\boxtimes\gamma_{U_T}\boxtimes\rho_T\cdot\textbf{1}_{N_T}\right).\end{align*}
We focus just on the $S$ part for now as the $T$ part will be similar. By inducing in stages, we see that
$$\pi|_{Sp(S)}\succ\Ind_{P_S}^{Sp(S)}\gamma_{K_S}\boxtimes\left(\Ind_{Sp(U_S)\times Sp(S')}^{Sp(U_S\oplus S')}\gamma_{U_S}\boxtimes\rho_S\right)\cdot\textbf{1}_{N_S}$$
We now have enough information to compute the type of this representation. Since $\rho_S\prec \rho|_{Sp(S')}$, using the induction hypothesis for restriction on $\rho$, if $\la \lambda_i\ra$ is the type of $\rho$, then $\rho_S$ has type $\la \lambda_i+2m-(a+1)-2b\ra$. Using the induction hypothesis for induction from $Sp(U_S)\times Sp(S')=Sp_{a-1}\times Sp_{2m-a+1}$ to $Sp(U_S\oplus S')$ subtracts a further $a-1$ from the type resulting in a type of $\la \lambda_i+2m-2n\ra$, which proves lemma \ref{restrictionsphericallemma} in this case.

\begin{remark}\label{restrictionorthogonalremark}
	The orthogonal group version of case 2 is much more similar to case 1 here. No unipotent radical appears in the $H$-stabilizer of a point in the open orbit of $G/P$, but instead when $a$ is odd we end up with induction from a product of odd orthogonal groups. We briefly sketch the argument. After handling the diagonal induction, the restriction to $O(S)$ would weakly contain something like 
	$$\Ind_{O_{2m-a}\times O_{a}}^{O_{2m}}\rho|_{O_{2m-a}}\boxtimes \gamma,$$ 
	where $\gamma$ is tempered. We may avoid solving this odd orthogonal induction problem by replacing $\gamma$ with $L^2(O_{a})$, which does not change the type for the induced representation above since $\gamma\prec L^2(O_a)$ and $L^2(O_a)$ has unramified spectrum concentrated on Arthur type representations of type $\la 1,\dots,1\ra$. Then since $L^2(O_a)$ is induced from the trivial group, we may use induction in stages and use the projection formula (\ref{projectionformula}) to obtain a representation of the form
	$$\Ind_{O_{2m-a+1}}^{O_{2m}}\rho|_{O_{2m-a+1}}\otimes L^2(O_{2m-a+1}/O_{2m-a}).$$
	The irreducible constituents of $L^2(O_{2m-a+1}/O_{2m-a})$ are described in \cite{Gan2014} as theta lifts from tempered representations of $Sp_{2}$, so by proposition \ref{typeofthetalift} it is of type $\la 2m-a-2,1,1,1\ra.$ Then using the induction hypothesis for tensor product and induction from $O_{2m-a+1}$ to $O_{2m}$, we may compute the type of the $O_{2m}$-representation to obtain the desired result.
\end{remark}	
\textbf{Case 3:} Suppose $2m<a\leq 2(n-m)$. Suppose $W$ is an $a$-dimensional isotropic subspace of $V$. In this most generic case, the projection of $W$ to $S$ is all of $S$ and the projection $W_T$ of $W$ to $T$ is an $a$-dimensional space. The kernel $K_T=W\cap W_T$ of the projection of $W$ to $T$ is an $a-2m$-dimensional isotropic subspace of $W_T$. The space $W$ defines the graph of a surjection $\phi:W_T\rightarrow S$ with $B(\phi(v),\phi(w))=-B(v,w)$ and $K_T=\ker \phi$, which induces an isomorphism $W_T/K_T\xrightarrow{\sim}S.$ It's clear that the collection of such subspaces forms an open $M'$ orbit and even an open $Sp(T)$ orbit in $G/P.$ Suppose $P$ is the stabilizer of such a subspace.

We first note that the restriction of $\Ind_P^G\sigma\otimes\rho\otimes\textbf{1}_N$ to $Sp(S)$ is tempered. Essentially, if we carry out the same procedure as in the first two cases, the stabilizer of $W$ in $M'$ will contain a copy of $Sp(S)$ which is diagonally embedded in $Sp(S)\times Sp(T)$, and upon inducing this, we will get a tempered representation of $Sp(S)$ by corollary 1 of \cite{Venkatesh2005}. For this case then, we will focus just on the stabilizer in $Sp(T)$ of $W$. Indeed, the $Sp(T)$ factor already acts with open orbit on $G/P$.

Choose an $a-2m$-dimensional isotropic space $K_T'$ which complements $K_T$ so that $W_T\oplus K_T'$ is a $2a-2m$-dimensional non-degenerate space. Let $T'=(W_T\oplus K_T')^\perp$ and $U_T=(K_T\oplus K_T'\oplus T')^\perp$ so that $W_T=U_T\oplus K_T.$ We may realize the Levi subgroup of $P$ as $M=GL(W)\times Sp(T').$

\hyperref[stabilizerfind]{\textbf{The stabilizer:}} The set of elements of $Sp(T)$ that stabilize $W$ must stabilize $W_T$ and hence the kernel $K_T$. Let $P_K$ be the parabolic subgroup of elements of $Sp(T)$ stabilizing $K_T.$ It has Levi subgroup $M_K=GL(K_T)\times Sp(U_T\oplus T')$ and unipotent radical $N_K.$ Since $P_K$ acts trivially on $S$, the only elements of $P_K$ which preserve $\phi$ are those which stabilize $K_T$ and act by the identity on $W_T/K_T\cong U_T.$

Of course the $GL(K_T)$ factor of $M_K$ is in this stabilizer. The only elements of $Sp(U_T\oplus T')$ acting trivially on $U_T$ are precisely the elements of $Sp(T').$ Since $N_K$ acts trivially on $K_T$ and $U_T\oplus T'/K_T$, all of $N_K$ will preserve $\phi$. The stabilizer of $\phi$ in $Sp(T)$ is then seen to be $(GL(K_T)\times Sp(T'))N_K.$ By Mackey theory, we have
$$\pi|_{Sp(T)}\succ\Ind_{(GL(K_T)\times Sp(T'))N_K}^{Sp(T)}\Res^{P}_{(GL(K_T)\times Sp(T'))N_K}\sigma\boxtimes\rho\cdot\textbf{1}_N.$$

\hyperref[stabilizerrestrict]{\textbf{Restriction to the stabilizer:}} Since $\rho$ is already a representation of $Sp(T')$, it has trivial $GL(K_T)$ and $N_K$ actions, and then of course $\sigma$ acts trivially by elements of $Sp(T')$. We have
$$\pi|_{Sp(T)}\succ\Ind_{(GL(K_T)\times Sp(T'))N_K}^{Sp(T)}\sigma|_{GL(K_T)N_K}\boxtimes\rho.$$

\hyperref[stabilizerinduce]{\textbf{Induction in stages:}} We must try to obtain a trivial action on $N_K$. Note that $Z(N_K)$ acts trivially on $W_T$, and so $\sigma$ acts trivially by elements of $Z(N_K)$. Since $N_K$ stabilizes $W_T$ and its center acts trivially, the $N_K$ action factors through $N_K/Z(N_K).$ We have 
$$N_K/Z(N_K)\cong\Hom(U_T\oplus T',K_T)=\Hom(U_T,K_T)\times \Hom(T',K_T),$$
and since $\Hom(T',K_T)$ acts trivially on $W_T$, the $\sigma$ action by elements of $N$ factors through $Z(N_K)$ to $\Hom(U_T,K_T).$

Let $P_W=M_WN_W$ be the parabolic subgroup of $GL(W)$ stabilizing $K_T\subset W$. For any complement $U_W$ to $K_T$ in $W$, we may realize $M_W$ as $GL(K_T)\times GL(U_W)$ and $N_W=\Hom(U_W,K_T)\cong \Hom(U_T,N_W).$ With this identification, we recognize $\sigma|_{GL(K_T)N_K}$ as the pullback of $\sigma|_{GL(K_T)N_W}$ via the trivial actions on $\Hom(T',K_T)$ and $Z(N_K).$

By Lemma \ref{restrictiontoparabolic}, $\sigma|_{GL(K_T)N_W}$ weakly contains a representation of the form $\sigma'\cdot\textbf{1}_{N_W}$, where $\sigma|_{GL(K_T)}\succ \sigma'.$ Then by inducing in stages, we have
\begin{align*}\pi|_{Sp(T)}&\succ\Ind_{(GL(K_T)\times Sp(T'))N_K}^{Sp(T)}\sigma'\boxtimes\rho\cdot\textbf{1}_{N_K}\\&=\Ind_{P_K}^{Sp(T)}\sigma'\boxtimes\left(\Ind_{Sp(U_T)\times Sp(T)}^{Sp(U_T\oplus T)}L^2(Sp(U_T))\boxtimes \rho\right)\cdot\textbf{1}_{N_K}\end{align*}

We apply the induction hypothesis for induction from $Sp(U_T)\times Sp(T)$ to $Sp(U_T\oplus T)$ and the restriction result for the $GL$ case to get the types of $\Ind_{Sp(U_T)\times Sp(T)}^{Sp(U_T\oplus T)}L^2(Sp(U_T))\boxtimes \rho$ and $\sigma'$ respectively. In particular, if $\la \lambda_i\ra$ is the type of $\rho$, then inducing from $Sp(U_T)\times Sp(T)=Sp_{2m}\times Sp_{2b}$ to $Sp_{2m+2b}$ results in a type of $\la \lambda_i-2m\ra$, and if $\la\tau_i\ra$ is the type of $\sigma$, then restricting from $GL(W)=GL_a$ to $GL(K_T)=GL_{a-2m}$ results in a type of $\la\tau_i-2m\ra$. In the end, we obtain the type $\la \lambda_i-2m\ra\oplus 2\la \tau_i-2m\ra$ for $\pi|_{Sp(T)}$, which proves lemma \ref{restrictionsphericallemma} in this case (replacing $2m$ in the claim with $2n-2m$ since we are restricting to $Sp(T)=Sp_{2n-2m}$ instead of $Sp(S)=Sp_{2m}$).

This concludes the proof of lemma \ref{restrictionsphericallemma}.\end{proof}
\subsection{Restriction from $Sp_{2n}$ to $GL_n$}\label{restrictionGLinduced}
In this section, we will prove the following subcase of proposition \ref{weakprop}.
\begin{lemma}\label{restrictionlevilemma}
	Let $M=GL_a\times Sp_{2b}$ be a Levi subgroup of a parabolic subgroup $P$ of $G=Sp_{2n}(\nQ_p)$, and suppose $\pi=\Ind_P^G\sigma\boxtimes\rho\cdot\textbf{1}$ is an unramified Arthur type representation of the form (\ref{stdrep}), where $\sigma$ is an unramified Arthur type representation of $GL_a$ and $\rho$ is a unipotent unramified Arthur type representation of $Sp_{2b}.$ If $\la \lambda_i\ra$ is the type of $\pi$, then $\pi|_{GL_n}$ weakly contains an unramified Arthur type representation of type $\la \lambda_i-n-1\ra.$
\end{lemma}

\begin{proof}The $GL_n$ Siegel Levi subgroup can be realized as the stabilizer of an isotropic splitting $V=X\oplus X^\ast$ of a $2n$-dimensional symplectic space $V$, where $\dim X= \dim X^\ast = n$, and the subgroup $P$ will be the stabilizer of an $a$-dimensional isotropic subspace $Y$ (here $P=MN$ with $M=GL_a\times Sp_{2b}).$ In the most generic case, the projections $W_Y$ and $W_Y^\ast$ of $Y$ to $X$ and $X^\ast$ respectively should be $a$-dimensional subspaces, and $Y$ should be the graph of an isomorphism $\phi:W_Y^\ast\rightarrow W_Y$. 

\hyperref[stabilizerfind]{\textbf{The stabilizer:}} Clearly any elements of $GL(X)$ which stabilize $Y$ must stabilize both $W_Y$ and $W_Y^\ast$. Let $P'=M'N'$ be the parabolic subgroup of $GL(X)$ consisting of elements which stabilize $W_Y$.  Choose complements $U$ and $U^\ast$ to $W_Y$ and $W_Y^\ast$ respectively so that $X=W_Y\oplus U$ and $X^\ast= W_Y^\ast \oplus U^\ast,$ and thus $V=W_Y\oplus U\oplus U^\ast\oplus W_Y^\ast.$ We have 
\begin{align*}
M&=GL(Y)\times Sp(U\oplus U^\ast)\\
M'&=GL(W_Y)\times GL(U),
\end{align*}
and $N'$ is the set of elements of $GL(X)$ which act by the identity on $W_Y$ and $X/W_Y.$ Here the $GL(U)$ factor of $M'$ is realized as a Siegel Levi subgroup of $Sp(U\oplus U^\ast).$

Observe that $W_Y\oplus W_Y^\ast$ is a $2a$-dimensional symplectic space. Since the elements of $GL(X)$ which stabilize $Y$ must stabilize this symplectic space they must also stabilize the elements of $(W_Y\oplus W_Y^\ast)^\perp=U\oplus U^\ast.$ No nontrivial elements of $N'$ stabilize $U\oplus U^\ast$, so no elements of $N'$ will be in the $GL(X)$-stabilizer of $Y$.

Consider now the isomorphism $\phi$. Since $Y=\graph \phi$ is isotropic, for all $x,y\in W_Y^\ast$ we have 
$$0=\la x+\phi(x),y+\phi(y)\ra=\la x,\phi(y)\ra-\la y,\phi(x)\ra,$$
where $\la\cdot,\cdot\ra$ denotes the symplectic form on $W$. Thus $(x,y)=\la x,\phi(y)\ra$ defines a nondegenerate symmetric bilinear form on $W_Y^\ast$, and similarly on $W_Y$ and $Y$. Thus the set of elements of $GL(W_Y)$ stabilizing $\phi$ is the set of elements stabilizing $(\cdot,\cdot)$, which is an orthogonal group denoted $O(Y)= GL(Y)\cap GL(W_Y)$.

The $GL(X)$-stabilizer of the point $P$ in the open orbit of $G/P$ is then $GL(U)\times O(Y).$ By Mackey theory, we have 
$$\Res^{Sp_{2n}}_{GL_n}\Ind_P^{Sp_{2n}}\sigma\boxtimes\rho\succ\Ind_{GL(U)\times O(Y)}^{GL_n}\Res^{P}_{GL(U)\times O(Y)}\sigma\boxtimes\rho.$$

\hyperref[stabilizerrestrict]{\textbf{Restriction to the stabilizer:}} In this case, the restriction to the stabilizer is easy to describe. Because $GL(U)$ induces the identity endomorphism on $Y$, the representation of $\sigma$ to $GL(U)$ is trivial. Meanwhile, $\rho$ is a representation of $Sp_{2b}=Sp(U\oplus U^\ast)$, so $\rho$ restricts trivially on $O(Y).$ Then, we simply have to consider 
$$\Ind_{O_a\times GL_b}^{GL_n}\sigma|_{O_a}\times \rho|_{GL_b}.$$

\hyperref[stabilizerinduce]{\textbf{Induction in stages:}} Here $GL_b$ occurs as the Siegel Levi subgroup of $Sp_{2b}$. Using the induction hypothesis for restriction to a $GL_b$, restricting $\rho$ will have its type reduced by $b+1$. We induce through the Levi subgroup $GL_a\times GL_b$ (this results in a tempered representation of $GL_a$ using the induction from $O_a$ to $GL_a$ case) and using the result of Venkatesh to induce from the Levi subgroup $GL_a\times GL_b$ of $GL_n$, the type will be further reduced by $a$. Thus we will get a reduction of the type of $\rho$ by $a+b+1=n+1$, which proves lemma \ref{restrictionlevilemma}.\end{proof}

\begin{remark}\label{restrictionGLinducedremark}
In the orthogonal case, instead of $O_a$ we get either $H=Sp_a$ if $a$ is even or a stabilizer of a degenerate skew-symmetric form with one-dimensional kernel if $a$ is odd, which will be a subgroup $H=(GL_1\times Sp_{a-1})N_H$ of a parabolic subgroup $P_H(GL_1\times GL_{a-1})N_H$ of $GL_a$. This weakly contains a representation which is trivial on $N_H$ by lemma \ref{restrictiontoparabolic}, we may take the action on the unipotent radical $N_H$ to be trivial, so we can induce in stages up to $GL_a$ from there through the parabolic. 

If $\sigma$ is one-dimensional, the induction $\Ind_H^{GL_a}\textbf{1}_H=L^2(GL_a/H)$ is not tempered. Instead, 
\begin{itemize}
	\item If $a$ is even so that $H=Sp_{a}$, then $L^2(GL_a/H)$ has type $\la 2,\dots,2\ra$
	\item If $a$ is odd so that $H=(GL_1\times Sp_{a-1})N_H$, then 
	$$L^2(GL_a/H)=\Ind_{P_H}^{GL_a}\textbf{1}_{GL_1}\boxtimes L^2(GL_{a-1}/Sp_{a-1})\cdot\textbf{1}_{N_H}$$
	has type $\la 1,2,\dots,2\ra.$
\end{itemize}
However, this does not matter as long as $b>0$, since inducing from the Levi subgroup $GL_a\times GL_b$ to $GL_n$ will subtract a further $b$ from this type, reducing it to $\la 1,\dots,1\ra.$ 

If $b=0$ however so that $GL_a=GL_n$ is a Levi subgroup of $O_{2n}$, then $\pi$ has type $\la n,n\ra$, and 
$$\pi|_{GL_n}=\text{ has type }\begin{cases}
\la 2,\dots,2\ra&\text{ if }n\text{ is even}\\
\la 1,2,\dots,2\ra&\text{ if }n\text{ is odd}
\end{cases}.$$
This accounts for the equality case for the exception in theorem \ref{mainthm} (3).
\end{remark}
\subsection{Restriction from $GL_{2n}$ to $Sp_{2n}$}
In this section, we will prove the following sub-case of proposition \ref{weakprop}
\begin{lemma}\label{restrictionGLtoSPlemma}
Let $M=GL_{\lambda_1}\times\dots\times GL_{\lambda_k}$ be a Levi subgroup for a parabolic subgroup $P=MN$ of $GL_{2n}(\nQ_p)$, where $2n=\sum_i \lambda_i$, and let $\pi=\Ind_P^{GL_n}\chi_1\boxtimes\dots\chi_k\cdot \textbf{1}_N$ be an unramified principal series Arthur type representation of $GL_{2n}$ of type $\la \lambda_1,\dots,\lambda_k\ra$ where the $\chi_i$ are unramified characters. Suppose $\lambda_k$ is the largest number in the type. Then $\pi|_{Sp_{2n}}$ weakly contains an Arthur type representation of type $\la 2\lambda_k-2n+1,1,\dots,1\ra$, except if $k=2$, in which case the type is $\la \lambda_2-\lambda_1+1,2,\dots,2\ra$ if $\lambda_1,\lambda_2$ are even or $\la \lambda_2-\lambda_1+1,1,1,2,\dots,2\ra$ if $\lambda_1,\lambda_2$ are odd.
\end{lemma}
\begin{proof}
We argue by induction: assuming the result holds true for $GL_{2m}$ for $m<n$. Let $m=\lambda_1+\dots+\lambda_{k-1}$. By inducing in stages, we may find a parabolic subgroup $P'=(GL_{\lambda_1}\times \dots\times GL_{\lambda_{k-1}})N'$ of $GL_{m}$ and a parabolic subgroup $P''=(GL_{m}\times GL_{\lambda_k})N''$ of $GL_n$ and take 
$$\sigma=\Ind_{P'}^{GL_{m}}\chi_1\boxtimes\dots\boxtimes\chi_{k-1}\cdot \textbf{1}_{N'}$$
so that 
$$\pi=\Ind_{P''}^{GL_n}\sigma\boxtimes \chi_k\cdot \textbf{1}_{N''}.$$
The set $G/P''$ can be identified with the set of all $m$-dimensional subspaces of a $2n$-dimensional vector space $V$. After fixing a nondegenerate skew-symmetric form $(\cdot,\cdot)$ on $V$, the subgroup $Sp_{2n}\subset GL_{2n}$ is the stabilizer of this form, and $Sp_{2n}$ acts with open orbit on $G/P''$. The open orbit corresponds to those $m$-dimensional spaces which are as nondegenerate as possible. We may assume $P''$ is the stabilizer of such a space $W$. There are two cases to consider: whether $m$ is even or $m$ is odd.

\textbf{Case 1:} Suppose $m$ is even. Then $\lambda_k=2n-m$ is even as well. In the most nondegenerate case, $W$ is an $m$-dimensional symplectic space. The stabilizer in $Sp_{2n}=Sp(V)$ of this space is $Sp(W)\times Sp(W^\perp).$ By Mackey theory, 
$$\pi|_{Sp(V)}\cong \Ind_{Sp(W)\times Sp(W^\perp)}^{Sp(V)}\sigma|_{Sp(W)}\boxtimes \chi_k|_{Sp(W^\perp)}.$$
Suppose first that $\sigma$ is not one-dimensional. By the induction hypothesis, $\sigma|_{Sp_m}$ has no numbers larger than $\lambda_k$ in its type. Applying the induction hypothesis for induction for induction from $Sp_{2a}\times Sp_{2b}$ to $Sp_{2a+2b}$ the type of $\pi|_{Sp(V)}$ is $\la 1,\dots,1\ra$ if $m\geq \lambda_k$ or 
$$\la \lambda_k-m+1,1,\dots,1\ra=\la 2\lambda_k-2n+1,1,\dots,1\ra$$
otherwise. 

If $\sigma$ is one-dimensional, then this means $k=2$, and we have $\lambda_2\geq \lambda_1$. In that case, applying the induction hypothesis for induction from $Sp_{\lambda_1}\times Sp_{\lambda_2}$ to $Sp_{2n}$, we obtain the type 
$$\la\lambda_2-\lambda_1+1,2,\dots,2\ra$$
as desired.

\textbf{Case 2:} Suppose $m$ is odd. In the most nondegenerate case, $(\cdot,\cdot)$ when restricted to $W$ has a one-dimensional kernel $K$. Choose an isotropic complement $K'$ to $K$ so that $W\oplus K'$ is an $m+1$-dimensional nondegenerate space. Let $W'=(W\oplus K)^\perp$ and $U=(W'\oplus K\oplus K')^\perp\subset W.$ We have the decomposition
$$V=U\oplus K\oplus K'\oplus W'.$$
The set of elements of $Sp(V)$ stabilizing $W$ must stabilize $K$. Let $P_S$ be the parabolic subgroup of $Sp(V)$ of elements stabilizing $K$. We have a Levi subgroup $M_S=GL(K)\times Sp(U\oplus W')$. The elements of $M_S$ which stabilize $W$ are a subgroup $GL(K)\times Sp(U)\times Sp(W')$. The unipotent radical $N_S$ of $P_S$ acts trivially on $K$ and $(U\oplus W')/K$, and hence it stabilizes $W=U\oplus K$. The stabilizer of a point in the open $Sp_{2n}$ orbit of $G/P''$ is then seen to be
$$(GL(K)\times Sp(U)\times Sp(W'))N_S.$$
Here $\chi_k$ is a representation of $GL(K'\oplus W')$, and so it restricts trivially to the $Sp(U)$ factor. Since $\chi_k$ factors through the determinant, it restricts trivially to $N_S$ as well. The representation $\sigma$ of $GL(W)$ restricts trivially to the $Sp(U)$ factor. Then
$$\pi|_{Sp(V)}\cong \Ind_{(GL(K)\times Sp(U)\times Sp(W'))N_S}^{Sp(V)}\sigma|_{(GL(K)\times Sp(U))N_S}\otimes \chi_k|_{GL(K)\times Sp(W')}.$$
The center $Z(N_S)$ of $N_S$ acts trivially on $W$, and so the action of $\sigma|_{N_S}$ factors through the center. We have
$$N_S/Z(N_S)\cong \Hom(U,K)\times \Hom(W',K).$$
Since $\sigma$ acts trivially via $W'$, the action of $\sigma$ is confined to the $\Hom(U,K)$ factor. Note that $(GL(K)\times GL(U))\Hom(U,K)$ is a parabolic subgroup of $GL(W)$. By lemma \ref{restrictiontoparabolic}, $\sigma|_{(GL(K)\times GL(U))\Hom(U,K)}$ weakly contains a representation of the form $\sigma_K\boxtimes\sigma_U\cdot\textbf{1}_{\Hom(U,K)}$ such that $\sigma|_{GL(K)\times GL(U)}\succ \sigma_K\boxtimes\sigma_U.$ It follows that 
$$\pi|_{Sp(V)}\succ \Ind_{(GL(K)\times Sp(U)\times Sp(W'))N_S}^{Sp(V)}\sigma_U|_{Sp(U)}\boxtimes (\sigma_K\otimes \chi_k|_{GL(K)})\boxtimes \chi_k|_{Sp(W')}\cdot \textbf{1}_{N_S}.$$
Let $\pi'=\Ind_{Sp(U)\times Sp(W')}^{Sp(U\oplus W')}\sigma_U|_{Sp(U)}\boxtimes \chi_k|_{Sp(W')}.$ By inducing in stages, we have
$$\pi|_{Sp(V)}\succ \Ind_{P_S} \pi'\boxtimes (\sigma_K\otimes \chi_k|_{GL(K)})\cdot \textbf{1}_{N_S}.$$
This is parabolic induction from which we can read the type of the representation once we compute the type of $\pi'.$ Suppose first that $\sigma$ is not one-dimensional. By the induction hypothesis, $\sigma_U|_{Sp(U)}$ does not contain any numbers in its type larger than $\lambda_k$. Then by the induction hypothesis for induction from $Sp_{2a}\times Sp_{2b}$ to $Sp_{2a+2b}$, the type of $\pi'$ is $\la 1,\dots,1\ra$ if $m-1\geq \lambda_k$ or 
$$\la \lambda_k-m+1,1,\dots,1\ra=\la 2\lambda_k-2n+1,1,\dots,1\ra$$
otherwise. Then $\pi$ has type $\la 2\lambda_k-2n+1,1,\dots,1\ra.$ 

If $\sigma_U$ is one-dimensional, which is only possible if $\sigma$ is one-dimensional, then $k=2$ and since $\lambda_1\leq \lambda_2$, the type of $\pi'$ is $\la \lambda_2-\lambda_1+1,2,\dots,2\ra$. Then the type of $\pi$ is $\la \lambda_2-\lambda_1+1,1,1,2,\dots,2\ra$ as desired. This concludes the proof of lemma \ref{restrictionGLtoSPlemma}.
\end{proof}
\begin{remark}
	In the case of restriction from $GL_{2n}$ to $O_{2n}$, for case 2 above with $m$ odd, $W$ is an $m$-dimensional nondegenerate space. In that case, following the notation of the proof,
	$$\pi\cong \Ind_{O(W)\times O(W^\perp)}^{O(V)}\sigma|_{O(W)}\boxtimes \chi_k|_{O(W^\perp)}.$$
	Knowing the type of this representation necessitates solving the problem of induction from $O(2a+1)\times O(2b+1)$ to $O(2a+2b+2).$ This is solvable using the boundary degeneration techniques of chapter 5 without needing to solve the rest of the odd orthogonal group problems. We omit the argument in this case.
\end{remark}
\subsection{Tensor product with an induced representation}\label{tensortwoinduced}
In this section, we will prove the following sub-case of proposition \ref{weakprop}.

\begin{lemma}\label{tensorlemma}
Let $M=GL_m\times Sp_{2(n-m)}$ be a Levi subgroup for the parabolic subgroup $P=MN$ of $G=Sp_{2n}(\nQ_p)$, and suppose $\pi=\Ind_P^G\sigma\boxtimes\rho\cdot\textbf{1}_N$ is unramified Arthur type representation of the form (\ref{stdrep}), where $\sigma$ is an unramified Arthur type representation of $GL_m$ and $\rho$ is a unipotent unramified Arthur type representations of $Sp_{2b}$. Suppose that $\pi'$ is another unramified Arthur type representation of $Sp_{2n}$ of any type.

If $\la \lambda_i\ra$ and $\la \tau_i\ra$ are the respective types of $\pi$ and $\pi'$, then $\pi\otimes \pi'$ weakly contains an unramified Arthur type representation of type $\la \lambda_i+\tau_j-(2n+1)\ra.$
\end{lemma}
\begin{proof}
This proof will be quite short as we have already done most of the work in section \ref{restrictiontoparabolicsection}. By proposition E.2.5 of \cite{Bekka2008},
$$\pi'\otimes \pi=\pi'\otimes\Ind_P^G\sigma\boxtimes\rho\cdot\textbf{1}_N=\Ind_P^G(\pi'|_P\otimes(\sigma\boxtimes\rho\cdot\textbf{1}_N)).$$
By lemma \ref{restrictiontoparabolic}, $\pi'|_P$ weakly contains a representation of the form $\sigma'\boxtimes \rho'\cdot\textbf{1}_N$, where $\sigma'$ and $\rho'$ are irreducible unramified Arthur type representations of $GL_m$ and $Sp_{2(n-m)}$ respectively. By continuity of induction with respect to the Fell topology, 
$$\pi'\otimes\pi\succ\Ind_P^G(\sigma\otimes\sigma')\boxtimes(\rho\otimes\rho')\cdot\textbf{1}_N.$$
Recall that in the induction hypothesis for proposition \ref{weakprop}, we are assuming the restriction results on $Sp_{2n}$ in order to show the tensor product result on $Sp_{2n}$ (while the restriction results themselves use the tensor product on $Sp_{2m}$ for $m<n$). By continuity of restriction with respect to the Fell topology, $\pi'|_{GL_m}\succ \sigma'$ and $\pi'|_{Sp_{2(n-m)}}\succ \rho'$. Thus the induction hypothesis tells us that the type of $\sigma'$ is $\la \tau_i+m-(2n+1)\ra$ while the type of $\rho'$ is $\la \tau_i+2(n-m)-2n\ra.$ Let $\la s_i\ra$ and $\la r_i\ra$ be the types of $\sigma$ and $\rho$ so that 
$$\la \lambda_i\ra=\la r_i\ra\oplus 2\la s_i\ra.$$ 
Applying the tensor product result of Venkatesh \cite{Venkatesh2005} (alternately Lapid and Rogawski \cite{Lapid2009}) for the tensor product on $GL_m$ and the induction hypothesis for tensor product on $Sp_{2(n-m)}$, we compute the types of $\sigma\otimes\sigma'$ and $\rho\otimes\rho'$ to be 
$$\la s_i+\tau_j-(2n+1)\ra\text{ and }\la r_i+\tau_j-(2n+1)\ra$$
respectively. After parabolically inducing to $G=Sp_{2n}$, we see that the type of $\pi'\otimes \pi$ is $\la \lambda_i+\tau_j-(2n+1)\ra$ as desired. This concludes the proof of lemma \ref{tensortwoinduced}.\end{proof}

\begin{remark}
In the case of $G=O_{2n}$ and $m=n$, there is one subtle case. We have $\pi=\Ind_{GL_nN}^G\sigma$, where $GL_n N$ is a Siegel parabolic subgroup, and in this case the type of $\pi$ is $\la n,n\ra$. Suppose that $\sigma$ is one-dimensional and $\pi'$ has type $\la \tau_1,\tau_2\ra$ with $\tau_1\leq \tau_2$. Then $\sigma'\prec \pi'|_{GL_n}$ has $SL_2$ type given by the exceptional case in theorem \ref{mainthm}. In particular, $\sigma'$ has type $\la \tau_2-(2n-1)+n,2,\dots,2\ra$ if $\tau_1,\tau_2$ are odd or $\la 2,\dots,2\ra$ if $\tau_1=\tau_2=n$ is even.

Since $\sigma$ is one-dimensional, $\sigma\otimes\sigma'$ has the same type as $\sigma'$. Then the type of $\pi\otimes \pi'$ is
$$\la\tau_2+n-(2n-1),\tau_2+n-(2n-1),2,\dots,2\ra$$
if $\tau_1,\tau_2$ are odd or $\la 2,\dots,2\ra$ if $\tau_1=\tau_2=n$ is even. This accounts for the equality case for the exception in part (5) of theorem \ref{mainthm}.
\end{remark}
	
\section{Restriction and tensor product part 2: the unipotent cases}
In the previous chapter, we showed how the restriction and tensor product cases of proposition \ref{weakprop} can be obtained for $Sp_{2n}$ representations under the induction hypotheses that the results of the proposition hold for $Sp_{2m}$ when $m<n$; however, those proofs were only sufficient for non-unipotent representations, i.e. those whose types do not consist solely of distinct odd numbers. Likewise, the same was true for representations of $O_{2n}$ under the induction hypothesis that proposition \ref{weakprop} holds for representations of $O_{2m}$ when $m<n$ for representations that don't have only distinct odd numbers in their types.

In this chapter, we are using a slight modification of that induction hypothesis which does not invalidate the inductive structure of the proof. To solve the $Sp_{2n}$ case, we are assuming that both the $Sp_{2m}$ and $O_{2m}$ cases hold, and we do so likewise for the $O_{2n}$ case. The reason for this is that we will be using theta correspondence to move between symplectic and orthogonal groups.

In the first section, we show how the restriction of a unipotent representation of $Sp_{2n}$ to an $Sp_{2a}\times Sp_{2b}$ subgroup reduces to a tensor product problem for a lower rank orthogonal group. In the second section, we show how the restriction of an $Sp_{2n}$ representation to $GL_n$ is either tempered or reduces to an induction problem from $O_{2m}$ to $GL_{2m}$ for some $m<n$. In the third section we show how the tensor product of two unipotent representations of an orthogonal group reduces to an induction problem for a lower rank symplectic group. 

For each $n$, let $\omega_n$ be a Weil representation for the dual pair $(O_{2m},Sp_{2n})$. Let $n=a+b$. Consider the following see-saw dual pair (\cite{Kudla1983} (2.4))
\begin{center}
	\begin{tikzcd}
		O_{2m}\times O_{2m}\arrow[dash]{dr}&Sp_{2n}\\
		O_{2m}^\Delta\arrow[hookrightarrow]{u}\arrow[dash]{ur}&Sp_{2a}\times Sp_{2b}\arrow[hookrightarrow]{u}
	\end{tikzcd}
\end{center}
where $O_{2m}^\Delta$ is a diagonally embedded copy of $O_{2m}$ in $O_{2m}\times O_{2m}$. Then if we consider $\omega_{n}$ as the Weil representation for the dual pair $(O_{2m}^\Delta,Sp_{2n})$, we have
$$\omega_n|_{O_{2m}^\Delta\times Sp_{2a}\times Sp_{2b}}\cong\omega_a|_{O_{2m}^\Delta\times Sp_{2a}}\otimes \omega_b|_{O_{2m}^\Delta\times Sp_{2b}}.$$
We shall exploit this isomorphism several times in sections \ref{restrictionsphericalunipotent} and \ref{tensortwounipotent}, and we will use an analogous isomorphism for a different see-saw dual pair in section \ref{restrictionGLunipotent}.

The main tools for this chapter come from section \ref{sectiondecayofmatrixcoefficients}. In particular, we reference Li's construction of the theta lift (theorem \ref{convergencetheorem}), the matrix coefficient bounds for unipotent representations (proposition \ref{matrixcoefficientsdecayratethetalifted}), and the theta correspondence for unramified Arthur type representations (proposition \ref{typeofthetalift}). In particular, recall that the theta lift of an unramified representation is nonzero. 
\subsection{Restriction from $Sp_{2n}$ to $Sp_{2a}\times Sp_{2b}$}\label{restrictionsphericalunipotent}
In this section we will prove the following sub-case of proposition \ref{weakprop}.
\begin{lemma}\label{restrictionofunipotentlemma}
Let $\la\lambda_i\ra$ be a partition of $2n+1$ consisting of distinct odd numbers, and suppose $a<n$. There exists an unramified Arthur type representation $\pi$ of $Sp_{2n}(\nQ_p)$ of type $\la \lambda_i\ra$ whose restriction to $Sp_{2a}$ weakly contains an unramfieed Arthur type representation of type $\la \lambda_i+2a-2n\ra.$
\end{lemma}
\begin{proof}The lemma is trivial for one-dimensional representations, so we assume $\la \lambda_i\ra$ is not the trivial type. We write $\la \lambda_i\ra=\la \lambda_1,\dots,\lambda_k,2n-2m+1\ra$ where $\lambda_1<\dots<\lambda_k<2n-2m+1$ and $2m=\lambda_1+\dots+\lambda_k.$ By corollary \ref{iteratedthetalifting}, there exists an unramified Arthur type representation $\theta(\rho)$ of $Sp_{2n}$ of type $\la \lambda_i\ra$ which is a theta lift from an unramified Arthur type representation $\rho$ of $O_{2m}$ of type $\la \lambda_1,\dots,\lambda_{k}\ra$. Let $n=a+b$ for some positive integers $a$ and $b$. We wish to analyze $\theta(\rho)|_{Sp_{2a}\times Sp_{2b}}$ and show that its Arthur type reduces to a problem about $\rho$. We may assume $a\leq b.$ We have three cases to consider:
\begin{enumerate}[(1)]
\item When $m<a\leq b$, we will show that $\theta(\rho)|_{Sp_{2a}}$ weakly contains a representation theta lifted from a tempered representation of $O_{2m}$.
\item When $a\leq m<b$, we will show that $\theta(\rho)|_{Sp_{2a}}$ is tempered (type $\la 1,\dots,1\ra$). We will choose a particular unramified tempered representation $\tau$ of $Sp_{2a}$ and let $\theta_a(\tau)$ be its theta lift to $O_{2m}$. Then for an unramified representation $\sigma$ weakly contained in $\theta(\tau)\otimes \rho$, we will show that its theta lift $\theta_b(\sigma)$ to $Sp_{2b}$ is weakly contained in $\theta(\rho)|_{Sp_{2b}}$. This is the only case that uses the induction hypothesis (for tensor product).
\item When $a\leq b\leq m$, we will show that $\theta(\rho)|_{Sp_{2a}\times Sp_{2b}}$ is tempered.
\end{enumerate}

\textbf{Case 1:} 
 Suppose $m< a\leq b$. According to the result of Sakellaridis \cite{Sakellaridis2017}, representations weakly contained in $\omega_b|_{O_{2m}}$ are tempered, and hence so are those contained in $\omega_b|_{O_{2m}}\otimes \rho.$ Then let $\tau$ be an unramified representation weakly contained in $\omega_b|_{O_{2m}}\otimes\rho$ (recall that tempered representations are automatically of Arthur type by theorem \ref{BSthm}). Since $\tau$ is tempered and unramified, it has nonzero theta lift $\theta_a(\tau)$ by proposition \ref{typeofthetalift} and may be understood using Li's construction \cite{Li1989} of the theta lift as we discussed in section \ref{sectiondecayofmatrixcoefficients}. We claim 
 $$\theta(\rho)|_{Sp_{2a}}\succ \theta_a(\tau).$$
 Note that in this case, the type of $\theta(\rho)$ has the entry $2n-2m+1$, and any other entries are less than or equal to $2m-1$, so the predicted type for $\theta(\rho)|_{Sp_{2a}}$ is $\la 2a-2m+1,1,\dots,1\ra$, which is the type of $\theta_a(\tau).$ The case of restriction to $Sp_{2b}$ will be similar.

 Denote by $\cY_a$ and $\cY_b$ the spaces for $\omega_a$ and $\omega_b$ and $H_\rho$ the space of $\rho$. By Li's construction, the diagonal matrix coefficients of $\theta(\rho)|_{Sp_{2a}\times Sp_{2b}}$ associated to pure tensors $\phi_a\otimes\phi_b\otimes v\in \cY_a\otimes\cY_b\otimes \cH_\rho$ may be expressed in the form 
$$(g_a,g_b)\mapsto \int_{O_{2m}}\la \omega_a(g_ah)\phi_a,\phi_a\ra\la\omega_b(g_bh)\phi_b,\phi_b\ra\la\rho(h)v,v\ra dh.$$
By restricting $\theta(\rho)$ further to $Sp_{2a}$, the diagonal matrix coefficients of $\theta(\rho)|_{Sp_{2a}}$ associated to pure tensors can be then expressed in the form
$$\sF_{\phi_a,\phi_b,v}(g_a)=\int_{O_{2m}}\la \omega_a(g_ah)\phi_a,\phi_a\ra\la\omega_b(h)\phi_b,\phi_b\ra\la\rho(h)v,v\ra dh.$$
The diagonal matrix coefficients of $\theta_a(\tau)$ associated to pure tensors on the other hand are of the form
$$\sG_{\phi_a,v'}(g_a)=\int_{O_{2m}}\la \omega_a(g_ah)\phi_a,\phi_a\ra\la\tau(h)v',v'\ra dh.$$
Let $Q_a$ be a compact subset of $Sp_{2a}$. We wish to approximate the functions $\sG_{\theta_a,v'}$ with finite sums of functions $\sF_{\phi_a,\phi_b,v}$ uniformly on $Q_a.$ We will however first refine our approach. Recall lemma F.1.3 of \cite{Bekka2008} which states that in order to show weak containment, we need only show that the $\sG_{\phi_a,v'}$ for $\phi_a\otimes v'$ in some total subset of $\cY_a\otimes \cH_\tau$ can be approximated, where $\cH_\tau$ is the space of $\tau$. 

In particular, let $K$ be a maximal compact subgroup of $O_{2m}.$ Since $K$-finite vectors are dense in $H_\tau$, we let $\cV$ be the subset of $\cY_a\otimes \cH_\tau$ consisting of pure tensors $\phi_a\otimes v_\mu'$ with $v_\mu'$ a unit vector in some $K$-type $\mu$ of $\tau.$

We may approximate the diagonal matrix coefficients of $\tau$ with those of $\omega_b|_{O_{2m}}\otimes\rho$ on compact subsets of $O_{2m}$. Even better, by proposition F.1.4 of \cite{Bekka2008}, since $\tau$ is irreducible, we may drop the finite sum requirement and approximate the diagonal matrix coefficients of $\tau$ associated to unit vectors with those of $\omega_b|_{O_{2m}}\otimes\rho$ directly uniformly on compacta.

Let $C_1\subset C_2\subset\dots$ be compact subsets of $O_{2m}$ such that $O_{2m}=\bigcup_i C_i.$ Then for each $i$ there exists a vector $\phi_{b_i}\otimes v_i$ of unit norm such that
$$|\la \tau(h)v'_\mu,v'_\mu\ra-\la\omega_b(h)\phi_{b_i},\phi_{b_i}\ra\la\rho(h)v_i,v_i\ra|<\frac{1}{i}$$
for all $h\in C_i$ (for the sake of simplicity, we assume that $\phi_{b_i}\otimes v_i$ is a pure tensor, but the argument works out fine without this assumption; we just end up with a finite sum of matrix coefficients associated to pure tensors instead). It immediately follows that $F_i(h):=\la\omega_b(h)\phi_{b_i},\phi_{b_i}\ra\la\rho(h)v_i,v_i\ra\rightarrow \la \tau(h)v'_\mu,v'_\mu\ra$ as $i\rightarrow \infty$ for all $h\in O_{2m}.$ We may refine this estimate a bit. Let $e_\mu$ be the central idempotent in $L^1(K)$ corresponding to the $K$-type $\mu.$ As Howe notes in Chapter 6 of his paper \cite{Howe2010}, 
$$(e_\mu\ast F_i\ast e_\mu^\ast)(h)=\la\omega_b(h)\psi_{b_i},\psi_{b_i}\ra\la\rho(h)w_i,w_i\ra,$$
where the $\ast$ operation denotes convolution over $K$, $e_\mu^\ast(h)=\overline{e_\mu(h^{-1})}$, and  $\psi_{b_i}=\omega_b(e_\mu)\phi_{b_i}$ $w_i=\rho(e_\mu)v_i$ are the projections of $\phi_{b_i}$ and $v_i$ to the $\mu$-isotypic subspaces of $\omega_b|_{O_{2m}}$ and $\rho$ respecitvely. Since $v'_\mu$ is already in $\cH_\mu$, if we denote $F(h)=\la \tau(h)v'_\mu,v'_\mu\ra$, then we see that $e_\mu\ast F\ast e_\mu^\ast = F.$ Therefore we see that $\tilde{F}_i(h):=\la \omega_b(h)\psi_{b_i},\psi_{b_i}\ra\la\rho(h)w_i,w_i\ra\rightarrow F(h)$ for all $h\in O_{2m}.$

Next, recall theorem 7.1 and lemma 6.2 (b) of \cite{Howe2010}, which together imply that since $\tau$ is tempered,
$$|F(h)|\leq (\dim\mu)^4\cdot\Xi(h),$$
 where $\Xi$ is Harish-Chandra's $\Xi$ function, which is an $L^{2+\epsilon}$ function on $O_{2m}$ for all $\epsilon>0$. Lemma 6.2b of \cite{Howe2010} tells us that because this estimate will hold for all representations in the support of $\omega_b|_{O_{2m}}\otimes\rho$, it also holds for the functions $\tilde{F_i}$. That is, for each $i$ we have the estimate
$$|\tilde{F_i}(h)|\leq \|\psi_{b_i}\|^2\|w_i\|^2(\dim\mu)^4\Xi(h)\leq (\dim\mu)^4\Xi(h)$$
since the $\psi_{b_i}$ and $w_i$ are projections of unit vectors. Hence, the sequence $\tilde{F_i}$ converges pointwise to $F$ and is uniformly bounded by an $L^{2+\epsilon}$ function. By Lebesgue's dominated convergence theorem,
$$\|\tilde{F_i}-F\|_{L^{2+\epsilon}(O_{2m})}\rightarrow 0$$
as $i\rightarrow\infty$ for all $\epsilon>0.$

Before we proceed with the final estimation for the matrix coefficients on the compact subset $Q_a$ of $Sp_{2a}$, we must take care of the $\phi_a$ part of the integral. Note that by Theorem 3.2 of \cite{Li1989}, for each fixed $g_a\in Q_a$, the function $G_{g_a}(h)= \la \omega_a(g_ah)\phi_a,\phi_a\ra$ is in $L^p(O_{2m})$ for $p>\frac{2m}{a}$. Since $m<a$, we may choose such a $p<2.$ Let $C=\max_{g_a\in Q_a}\|G_{g_a}\|_{L^p(O_{2m})}$, which exists because $Q_a$ is compact and matrix coefficients are continuous. 

Finally, let $\frac{1}{p}+\frac{1}{q}=1$ so that $q>2.$ we apply H{\"o}lder's inequality to see that for all $g_a\in Q_a$,
\begin{align*}
|(\sG_{\phi_a,v'_\mu}&-\sF_{\phi_a,\psi_{b_i},w_i})(g_a)|\\
&=\left|\int_{O_{2m}}\la\omega(g_ah)\phi_a,\phi_a\ra(\la\tau(h)v'_\mu,v'_\mu\ra-\la\omega_b(h)\psi_{b_i},\psi_{b_i}\ra\la\rho(h)w_i,w_i\ra) dh\right|\\
&\leq \int_{O_{2m}}|G_{g_a}(h)(F(h)-\tilde{F_i}(h))|dh\\
&\leq \|G_{g_a}\|_{p}\|F-\tilde{F_i}\|_q\\
&\leq C\|F-\tilde{F_i}\|_q.
\end{align*}
This approaches $0$ as $i\rightarrow \infty$; thus we have approximated the matrix coefficients of $\theta_a(\tau)$ in the total subset $\cV$ with those of $\theta(\rho)|_{Sp_{2a}}$ uniformly on $Q_a.$ Therefore $\theta_a(\tau)$ is weakly contained in $\theta(\rho)|_{Sp_{2a}}.$ This proves lemma \ref{restrictionofunipotentlemma} in this case.
\begin{remark}\label{alternateproofremark}
	There is an alternate way to prove case 1 of lemma \ref{restrictionofunipotentlemma}. It is possible using the argument of section \ref{restrictionGLunipotent} to show that if $\sigma_a\boxtimes\sigma_b$ is a tempered representation of $O_{2m}\times O_{2m}$ weakly contained in $\Ind_{O_{2m}}^{O_{2m}\times O_{2m}}\rho$, then 
	$$\theta(\rho)|_{Sp_{2a}\times Sp_{2b}}\succ\theta_a(\sigma_a)\boxtimes\theta_b(\sigma_b).$$
\end{remark}

\textbf{Case 2:} Suppose $a\leq m< b$. By the result of Sakellaridis \cite{Sakellaridis2017}, $\omega_a$ is supported on representations of the form $\tau\otimes\theta_a(\tau)$, where $\tau$ is a tempered representation of $Sp_{2a}$ and $\theta_a(\tau)$ is its theta lift to $O_{2m}.$ Let $\tau$ be one unramified such representation (which again is automatically of Arthur type). By proposition \ref{typeofthetalift}, $\theta_a(\tau)$ is an Arthur type unramified representation of type $\la 2m-2a-1,1,\dots,1\ra.$

Let $\sigma$ be an unramified (Arthur type) representation weakly contained in $\theta_a(\tau)\otimes\rho$. By the induction hypothesis for tensor product on $O_{2m}$, $\sigma$ has type 
$$\la \lambda_1-2a,\dots,\lambda_k-2a\ra.$$
Let $\theta_b(\sigma)$ be the theta lift of $\sigma$ to $Sp_{2b}$, which is nonzero and of type 
$$\la \lambda_1-2a,\dots,\lambda_k-2a,2n-2m-2a+1\ra$$
by proposition \ref{typeofthetalift}. We claim 
$$\theta(\rho)|_{Sp_{2a}\times Sp_{2b}}\succ\tau\otimes\theta_b(\sigma).$$
Observe that $\theta_b(\sigma)$ has the type predicted by lemma \ref{restrictionofunipotentlemma}. Also observe that $\theta_b(\sigma)$ may be constructed using Li's construction by corollary \ref{convergencecorollary}.

As in the first case, the diagonal matrix coefficients of $\theta(\rho)|_{Sp_{2a}\times Sp_{2b}}$ associated to pure tensors are of the form
$$\sF_{\phi_a,\phi_b,v}(g_a,g_b)=\int_{O_{2m}}\la \omega_a(g_ah)\phi_a,\phi_a\ra\la \omega_b(g_bh)\phi_b,\phi_b\ra\la\rho(h)v,v\ra dh.$$
The diagonal matrix coefficients of $\tau\otimes \theta_b(\sigma)$ associated to pure tensors on the other hand are of the form
$$\sG_{v',\phi_b,v}(g_a,g_b)=\la\tau(g_a)v',v'\ra\int_{O_{2m}}\la\omega_b(g_bh)\phi_b,\phi_b\ra\la\sigma(h)v,v\ra dh.$$
We wish to approximate the latter with finite sums of the former on a given compact subset $Q_a\times Q_b$ of $Sp_{2a}\times Sp_{2b}$. 

Once again, let $C_1\subset C_2\subset\dots$ be compact subsets of $O_{2m}$ such that $O_{2m}=\bigcup_i C_i$, and let $\cV$ be the subset of unit pure tensors $v'\otimes \phi_a\otimes v\in \cH_{\tau}\otimes \cY_a\otimes \cH_\rho$ such that $v'=v'_\mu$ is an element of a $K$-type $\mu$ of $\tau$ and $v=v_\nu$ is an element of a $K$-type $\nu$ of $\rho.$ As we noted in the first case, $\cV$ is a total subset, and so it suffices to show that we can approximate $\sG_{v'_\mu,\phi_b,v_\nu}.$

The diagonal matrix coefficients of $\theta_a(\tau)\otimes \rho$ are of the form
$$h\mapsto \la\rho(h)v,v\ra\la\theta_a(\tau)(h)v',v'\ra.$$
Since $\sigma$ is irreducible and weakly contained in $\theta_a(\tau)\otimes\rho$, this means that for each $i$, there exists $v_i\otimes v_i'\in \cH_{\rho}\otimes\cH_{\theta_a(\tau)}$ (again we use a pure tensor for simplicity but it does not affect the argument) such that
$$\left|\la \sigma(h)v_\nu,v_\nu\ra-\la\rho(h)v_i,v_i\ra\la\theta_a(\tau)(h)v'_i,v'_i\ra\right|<\frac{1}{2i}$$
for all $h\in C_i$. 

Now consider the matrix coefficients of $\tau\otimes \theta_a(\tau)$, which are of the form
$$\sH_{v',v}(g_a,h)= \la\tau(g_a)v',v'\ra\la\theta_a(\tau)(h)v,v\ra.$$
We may approximate these by the matrix coefficients of $\omega_a.$ In particular, for each $i$, we choose $\psi_{a_i}$ for each $i$ so that
$$\left|\la\tau(g_a)v'_\mu,v'_\mu\ra\la\theta_a(\tau)(h)v'_i,v'_i\ra -\la\omega_a(g_ah)\psi_{a_i},\psi_{a_i}\ra\right|<\frac{1}{2i}$$
uniformly for $g_a\in Q_a$ and $h\in C_i.$ Let $F_i(g_a,h)=\la\omega_a(g_ah)\psi_{a_i},\psi_{a_i}\ra\la\rho(h)v_i,v_i\ra.$ By the triangle inequality and exploiting the fact that all vectors involved are unit vectors, we have
$$\left|\la\sigma(h)v_\nu,v_\nu\ra\la\tau(g_a)v_\mu',v_\mu'\ra-F_i(g_a,h)\right|<\frac{1}{i}$$
uniformly for $g_a\in Q_a$ and $h\in C_i$ , and thus $F_i\rightarrow F:=\la\sigma(h)v_\nu,v_\nu\ra\la\tau(g_a)v_\mu',v_\mu'\ra$ for $g_a\in Q_a$ (uniformly) and $h\in O_{2m}$ (pointwise)

We now need to make a dominated convergence argument again in the same way as in the previous case. Treat $g_a$ as fixed for now. Let $e_\nu\in L^1(K)$ be the central idempotent of $\nu$. Then, convolving over $K$, $e_\nu\ast F_i(g_a,\cdot)\ast e_\nu^\ast\rightarrow e_\nu \ast F(g_a,\cdot)\ast e_\nu^\ast=F(g_a,\cdot)$ for each $g_a\in Q_a.$ For each $g_a\in Q_a$, let 
$$\tilde{F}_i(g_a,h)=(e_\nu\ast F_i(g_a,\cdot)\ast e_\nu^\ast)(h)=\la\omega_a(g_ah)\xi_{i},\xi_{i}\ra\la\rho(h)w_i,w_i\ra,$$
where $w_i=\omega_b(e_\nu)v_i$ and $\xi_{i}=\omega_a(e_\nu)(\psi_{a_i})$ (note $\omega_a(e_\nu)$ commutes with $\omega_a(g_a)$) so that $\tilde{F_i}\rightarrow F$ pointwise in $O_{2m}$ and uniformly on $Q_a.$

Keeping $g_a$ fixed, it's clear that $\tilde{F_i}$ is a matrix coefficient of $\omega_a|_{O_{2m}}\otimes \rho.$ Let us consider the $L^p$ decay of its matrix coefficients so that we may apply Lebesgue's dominated convergence theorem to find the values of $p$ for which $\|\tilde{F_i}-F\|_{L^p(O_{2m})}\rightarrow 0$ as $i\rightarrow \infty.$

Recall that the result of Sakellaridis \cite{Sakellaridis2017} says that all of the representations in the spectrum of $\omega_a|_{O_{2m}}$ are theta lifts of tempered representations of $Sp_{2a}$. By corollary \ref{matrixcoefficientsdecayratethetalifted}, such a representation  has $K$-finite matrix coefficients bounded by a multiple of $\Xi_{O_{2m}}^\frac{a}{m-1}$, and lemma 6.2 of \cite{Howe2010} says that then $\omega_a|_{O_{2m}}$ has this bound as well. Meanwhile, from corollary \ref{matrixcoefficientsdecayratethetalifted}, $\rho$ has $K$-finite matrix coefficients bounded by a multiple of $\Xi_{O_{2m}}^\frac{n'}{m-1}$, where $2n'=2m-\lambda_k-1$ (recall $\lambda_k$ is the largest number in the type of $\rho$). Then $\tilde{F_i}$ is bounded by a multiple of $\Xi_\frac{a+n'}{m-1}$, where the multiple in question depends on $\|\xi_i\|\leq 1$ and $\|w_i\|\leq 1$ for all $i$. 

Let $p_0=\frac{2(m-1)}{n'+a}$. Then we've shown all $\tilde{F_i}$ are uniformly bounded by a function which is in $L^p(O_{2m})$ for all $p>p_0.$ By Lebesgue's dominated convergence theorem, we conclude $\|\tilde{F_i}-F\|_{L^p(O_{2m})}\rightarrow 0$ as $i\rightarrow \infty$ for all $p>p_0$. Since the convergence $\tilde{F_i}\rightarrow F_i$ was already uniform for $g_a\in Q_a$, the convergence $\|\tilde{F_i}-F\|_{L^p(O_{2m})}\rightarrow 0$ for $p>p_0$ is still uniform on $Q_a.$

Recall from the proof of theorem \ref{convergencetheorem} that $\omega_b|_{O_{2m}}$ has matrix coefficients in $L^q$ for all $q>q_0$, where $q_0=\frac{2(m-1)}{b}.$ Observe that
$$\frac{1}{p_0}+\frac{1}{q_0}=\frac{a+b+n'}{2m-2}=\frac{2n+2m-\lambda_k-1}{4m-4}> \frac{4m-2}{4m-4}>1$$
since $2n-2m+1>\lambda_k$. Then we may choose $p>p_0$ and $q>q_0$ such that $\frac{1}{p}+\frac{1}{q}=1.$ By H{\"o}lder's inequality,
\begin{align*}
|\sG_{v_\mu',\phi_b,v_\nu}(g_a,g_b)&-\sF_{\xi_{i},\phi_b,w_i}(g_a,g_b)|\\
&=\Big|\int_{O_{2m}}\la\omega_b(g_bh)\phi_b,\phi_b\ra\la\tau(g_a)v_\mu',v_\mu'\ra\la\sigma(h)v_\nu,v_\nu\ra\\
&\quad- \la \omega_a(g_ah)\xi_{i},\xi_{i}\ra\la \omega_b(g_bh)\phi_b,\phi_b\ra\la\rho(h)w_i,w_i\ra dh\Big|\\
&\leq\int_{O_{2m}}|\la\omega_b(g_bh)\phi_b,\phi_b\ra(F(g_a,h)-\tilde{F_i}(g_a,h))|dh\\
&\leq\|\la\omega_b(g_b \cdot)\phi_b,\phi_b\ra\|_{q}\cdot\|F(g_a,\cdot)-\tilde{F_i}(g_a,\cdot)\|_p.
\end{align*}
In the first factor, we may take the maximum over all $g_b\in Q_b$ to obtain a constant multiple. In the second factor, since the $L^p$ convergence was uniform in $Q_a$, that factor will tend to zero. Therefore $\theta(\rho)|_{Sp_{2a}\times Sp_{2b}}$ weakly contains $\tau\otimes \theta(\sigma)$. This proves lemma \ref{restrictionofunipotentlemma} in this case.\\

\begin{remark}In the case $m=a$, we can't use the construction of Li for $\theta_a(\tau)$ since the integral (\ref{Liintegral}) does not converge, which is why we refrain from writing out its matrix coefficients in terms of that construction and why we prove the $m=a$ case in this part rather than in case 1, which relies on Li's construction throughout.\end{remark}

\textbf{Case 3:} Suppose $a\leq b\leq 
m$. In this case, $\theta(\rho)|_{Sp_{2a}\times Sp_{2b}}$ contains only tempered representations. Recall that the type of $\theta(\rho)$ is $\la\lambda_1,\dots,\lambda_k,2n-2m+1\ra$ with $\lambda_1<\dots<\lambda_k<2n-2m+1.$

We calculate just the restriction to $Sp_{2b}$ since the restriction to $Sp_{2a}$ follows by restricting in stages. Restricting first through the subgroup $Sp_{2a}\times Sp_{2b}$, the matrix coefficients of $\theta(\rho)|_{Sp_{2b}}$ look like finite sums of functions of the form
$$F(g)=\int_{O_{2m}}\la \omega_b(gh)\phi,\phi\ra\la\omega_a(h)\psi,\psi\ra\la\rho(h)v,v\ra dh.$$
We will assume that we are working with smooth matrix coefficients only since that is sufficient to show a representation is tempered. Recall from the proof of proposition \ref{matrixcoefficientsdecayrate} that 
\begin{align*}
\la\omega_b(gh)\phi,\phi\ra&\leq C(\phi)\cdot\Xi_{O_{2m}}(h)^\frac{b}{m-1}\cdot\Xi_{Sp_{2b}}(g)^\frac{m}{b-1},\\
\la\omega_a(h)\psi,\psi\ra&\leq D(\psi)\cdot\Xi_{O_{2m}}(h)^\frac{a}{m-1}
\end{align*}
for some constants $C(\phi)$ and $D(\psi)$.

Suppose first that $\rho$ is one-dimensional. Its matrix coefficients $\la\rho(h)v,v\ra$ are bounded by $\|v\|^2.$ Then $2m=\lambda_k+1<2n-2m+1$, so $a+b=n>2m$ we have
\begin{align*}
|F(g)|&\leq C(\phi)D(\psi)\|v\|^2 \Xi_{Sp_{2b}}(g)^\frac{m}{b}\int_{O_{2m}}\Xi_{O_{2m}}(h)^\frac{b}{m-1}\cdot\Xi_{O_{2m}}(h)^\frac{a}{m-1}dh\\
&= C(\phi)D(\psi)\|v\|^2 \Xi_{Sp_{2b}}(g)^\frac{m}{b}\int_{O_{2m}}\Xi_{O_{2m}}(h)^\frac{n}{m-1}dh.
\end{align*}
Since $\frac{n}{m-1}>2$, the $O_{2m}$ integral converges. Then $\|F\|_{L^p(Sp_{2b})}^p$ is bounded by a multiple of $\int_{Sp_{2b}}\Xi_{Sp_{2b}}(g)^\frac{pm}{b}dg.$ Since $m\geq b$, this converges for all $p>2$. We conclude $\theta(\rho)|_{Sp_{2b}}$ is strongly $L^p$ for $p>2$, and hence it is supported on tempered representations by lemma 6.2 of \cite{Howe2010}.

Now suppose $\rho$ is nontrivial so that $\rho$ is an iterated theta lift from a one-dimensional character. In that case by corollary \ref{matrixcoefficientsdecayratethetalifted}, $\rho$ has matrix coefficients in $L^p$ for $p>\frac{2m-2}{n'}$, where 
$$2m-\lambda_k-1=2n'=\lambda_1+\dots+\lambda_{k-1}$$ 
(here $\rho$ is theta lifted from a representation of $Sp_{2n'}$), and we have the bound
$$\la \rho(h)v,v\ra\leq B(v)\cdot\Xi_{O_{2m}}(h)^\frac{n'}{m-1}$$
for some constant $B(v)$. Then 
\begin{align*}
|F(g)|&\leq B(v)C(\phi)D(\psi)\Xi_{Sp_{2b}}(g)^\frac{m}{b}\int_{O_{2m}}\Xi_{O_{2m}}(h)^\frac{b}{m-1}\Xi_{O_{2m}}(h)^\frac{a}{m-1}\Xi_{O_{2m}}(h)^\frac{n'}{m-1}dh\\
&=B(v)C(\phi)D(\psi)\Xi_{Sp_{2b}}(g)^\frac{m}{b}\int_{O_{2m}}\Xi_{O_{2m}}(h)^\frac{n+n'}{m-1}dh.
\end{align*}
Since $2n-2m+1>\lambda_k$, we have
$$2n+2n'=2n+2m-\lambda_k-1>4m-2>4(m-1).$$
Then $\frac{n+n'}{m-1}>2$, so the $O_{2m}$ integral above converges. Similarly to the previous case, we then conclude $F\in L^p(Sp_{2b})$ when $p>2$, and so $\theta(\rho)|_{Sp_{2b}}$ is supported on tempered representations by lemma 6.2 of \cite{Howe2010}. We conclude $\theta(\rho)|_{Sp_{2b}}$ has type $\la 1,\dots,1\ra$ in this case. This concludes the proof of lemma \ref{restrictionofunipotentlemma}.\end{proof}

\subsection{Restriction from $Sp_{2n}$ to $GL_{n}$}\label{restrictionGLunipotent}
In this section, we will prove the following sub-case of proposition \ref{weakprop}.
\begin{lemma}\label{restrictionofunipotenttoGLlemma}
	Let $\la\lambda_i\ra$ be a partition of $2n+1$ consisting of distinct odd numbers. There exists an unramified Arthur type representation $\pi$ of $Sp_{2n}(\nQ_p)$ of type $\la \lambda_i\ra$ whose restriction to $GL_n$ weakly contains an unramfied Arthur type representation of type $\la \lambda_i-n-1\ra.$
\end{lemma}
\begin{proof}The lemma is trivial in the one-dimensional case, so we may assume $\la\lambda_i\ra$ is not the trivial type. We write $\la \lambda_i\ra=\la\lambda_1,\dots,\lambda_k,2n-2m+1\ra$, where $\lambda_1<\dots<\lambda_k<2n-2m+1$ and $2m=\lambda_1+\dots+\lambda_k$. By corollary \ref{iteratedthetalifting}, there exists an unramified Arthur type representation $\theta(\rho)$ of $Sp_{2n}$ of type $\la \lambda_i\ra$ which is a theta lift from an unramified Arthur type representation $\rho$ of $O_{2m}$ of type $\la\lambda_1,\dots,\lambda_k\ra$. We have two cases to consider.
\begin{enumerate}[(1)]
	\item When $2m<n$, we will show $\theta(\rho)|_{GL_n}$ weakly contains a representation theta lifted from a tempered representation of $GL_{2m}$. This will produce a representation of type $\la n-2m,1,\dots,1\ra=\la 2n-2m+1+n-(2n+1),1,\dots,1\ra$ as required.
	\item When $2m\geq n$, we will show  $\theta(\rho)|_{GL_n}$ is tempered.
\end{enumerate}

\textbf{Case 1:} Suppose $2m<n$. In this case, $\theta(\rho)$ is a low rank representation and can be understood using Howe's theory of low rank unitary representations \cite{Howe2010}. That is, if $P=(GL_n)N\subset Sp_{2n}$ is the Siegel parabolic, then $\theta(\rho)|_{N}$ is supported on rank $2m$ characters. The $GL_n$-stabilizer of rank $2m$ character of $N$ is a subgroup of the form 
$$(O_{2m}\times GL_{n-2m})N'\subset P'=M'N',$$ 
where $P'$ is a parabolic subgroup of $GL_n$ with Levi subgroup $M'=GL_{2m}\times GL_{n-2m}.$ From the proof of lemma \ref{restrictiontoparabolic}, we have
$$\theta(\rho)|_{GL_n}\succ\Ind_{(O_{2m}\times GL_{n-2m})N'}^{GL_n}\rho\boxtimes\textbf{1}_{GL_{n-2m}}\cdot \textbf{1}_{N'}.$$
By inducing in stages,
$$\theta(\rho)|_{GL_n}\succ\Ind_{P'}^{GL_n}\left(\Ind_{O_{2m}}^{GL_{2m}}\rho\right)\boxtimes\textbf{1}_{GL_{n-2m}}\cdot \textbf{1}_{N'}.$$
From the induction from $O_{2m}$ to $GL_{2m}$ case, we see that the type of this representation is $\la n-2m,1,\dots,1\ra$ as desired.
\begin{remark}
	In the case of restriction from $O_{2n}$ to $GL_n$, where $\rho$ is an unramified unipotent Arthur type representation of $Sp_{2m}$ and $\theta(\rho)$ is its lift to $O_{2n}$, the argument is largely the same but with one exceptional case. If $\rho$ is one dimensional, then $\Ind_{Sp_{2m}}^{GL_{2m}}\rho$ has type $\la 2,\dots,2\ra$, so the final type is for $\theta(\rho)|_{GL_n}$ is $\la n-2m,2,\dots,2\ra.$
\end{remark}

\textbf{Case 2:} Suppose $2m\geq n$. We shall exploit the see-saw dual pair (\cite{Kudla1983} (2.23))
\begin{center}
	\begin{tikzcd}
		GL_{2m}\arrow[dash]{dr}&Sp_{2n}\\
		O_{2m}\arrow[hookrightarrow]{u}\arrow[dash]{ur}&GL_n\arrow[hookrightarrow]{u}
	\end{tikzcd}.
\end{center}
Thus if $\omega_H$ is a Weil representation for the dual pair $(O_{2m},Sp_{2n})$ and $\omega_G$ is a Weil representation for the dual pair $(GL_{2m},GL_n)$, then we have 
$$\omega_H|_{O_{2m}\times GL_n}=\omega_G|_{O_{2m}\times GL_n}.$$
Let $\rho$ be a unipotent unramified Arthur type representation of $O_{2m}$, and let $\theta(\rho)$ be its theta lift to $Sp_{2n}$. This case is much easier than the previous case and is similar to case 3 of section \ref{restrictionsphericalunipotent}, as we will show $\theta(\rho)|_{GL_n}$ weakly contains only tempered representations.

Consider the diagonal matrix coefficients of $\theta(\rho).$ They are of the form
$$g\mapsto \int_{O_{2m}}\la \omega_H(gh)\phi,\phi\ra\la\rho(h)v,v\ra dh,$$
but when restricted to $GL_n$, using the see-saw dual pair identity, the matrix coefficients of $\theta(\rho)|_{GL_n}$ are of the form
$$F(g)=\int_{O_{2m}}\la \omega_{G}(gh)\phi,\phi\ra\la\rho(h)v,v\ra dh.$$
Recall from section 2.3 the bounds on the smooth matrix coefficients of the Weil representation:
\begin{align*}
|\la\omega_{G}(gh)\phi,\phi\ra|&\leq A(\phi)\Xi_{GL_{2m}}(h)^\frac{n}{2m-1}\Xi_{GL_n}(g)^\frac{2m}{n-1}\\
|\la\omega_G(gh)\phi,\phi\ra|&\leq B(\phi)\Xi_{O_{2m}}(h)^\frac{n}{m-1}\Xi_{Sp_{2n}}(g)^\frac{2m}{n}
\end{align*}
for all $h\in GL_{2m}$ $g\in GL_n$ in the first line or $h\in O_{2m}$ and $g\in Sp_{2n}$ in the second line and for some constants $A(\phi)$ and $B(\phi).$ It follows from (22) of \cite{Li1989}, (7.7) of \cite{Howe2010}, and the fact that $\omega_G|_{O_{2m}\times GL_n}=\omega_H|_{O_{2m}\times GL_n}$ that we have the combined estimate
\begin{equation}\label{eqn425}|\la\omega_{G}(gh)\phi,\phi\ra|\leq C(\phi)\Xi_{O_{2m}}(h)^\frac{n}{m-1}\Xi_{GL_n}(g)^\frac{2m}{n-1}\end{equation}
when $h\in O_{2m}$ and $g\in GL_n$ for some constant $C(\phi)$.

Suppose for now that $\rho$ is one-dimensional. Note $|\la \rho(h)v,v\ra|\leq \|v\|^2$ for all $h\in O_{2m}$. Also our assumption that $2n-2m+1$ is the largest number in the type of $\rho$ means $2n-2m+1>2m-1$, so $n>2m-1$. Since in this case we are assuming $2m\geq n$, the only possibility is $n=2m.$ Then
$$|F(g)|\leq \|v\|^2\cdot C(\phi)\cdot \Xi_{GL_{n}}(g)^\frac{2m}{2m-1}\int_{O_{2m}}\Xi_{O_{2m}}(h)^\frac{2m}{m-1}dh.$$
The integral converges since $m>m-1$, and so $|F(g)|$ is bounded by a multiple of $\Xi_{GL_{n}}^\frac{2m}{2m-1}.$ We conclude $\theta(\rho)|_{GL_{n}}$ is tempered in this case.

Now suppose $\rho$ is not one-dimensional so that $\rho$ is itself a theta lift from a representation of $Sp_{2n'}$, where $2n-2m+1>2m-2n'-1$ (hence $n'>2m-n-1$). Recall that by corollary \ref{matrixcoefficientsdecayratethetalifted}, we have the following bound for the $K$-finite matrix coefficients of $\rho$:
$$|\la\rho(h)v,v\ra|\leq D(v)\cdot \Xi_{O_{2m}}(h)^\frac{n'}{m-1},$$
where $D(v)$ is some constant. Hence we have
$$|F(g)|\leq C(\phi)D(v)\cdot \Xi_{GL_n}(g)^\frac{2m}{n-1}\int_{O_{2m}}\Xi_{O_{2m}}(h)^\frac{n}{m-1}\cdot\Xi_{O_{2m}}(h)^\frac{n'}{m-1} dh.$$
Note that $n+n'>2m-1>2(m-1)$, and so the integral converges. Since $2m\geq n-1$, the right side of the above inequality is in $L^p(GL_n)$ for all $p>2.$ We conclude once again $\theta(\rho)|_{GL_n}$ is tempered. We have then shown $\theta(\rho)|_{GL_n}$ has type $\la 1,\dots,1\ra$. This concludes the proof of lemma \ref{restrictionofunipotenttoGLlemma}.\end{proof}

\subsection{Tensor product of two unipotent representations}\label{tensortwounipotent}
While most of the proofs have been done for $Sp_{2n}$ representations, here we do it for $O_{2m}$ as a stylistic choice since the restriction of unipotent case reduced to a tensor product on the orthogonal side. It is also more interesting in that the exceptional case of theorem \ref{mainthm} arises on the orthogonal side. We will prove the following sub-case of proposition \ref{weakprop}.
\begin{lemma}\label{tensorproductofunipotentslemma}
Let $\la\lambda_i\ra$ and $\la\tau_i\ra$ be two partitions of $2m$ such that each partition consists only of distinct odd numbers. There exist two unramified Arthur type representations $\pi_a$ and $\pi_b$ of $O_{2m}(\nQ_p)$ of respective types $\la\lambda_i\ra$ and $\la\tau_i\ra$ such that $\pi_a\otimes\pi_b$  weakly contains an unramified Arthur type representation $\sigma$ which is of type $\la \lambda_i+\tau_i-(2m-1)\ra$, except in the case where the given partitions are $\la \lambda_1,\lambda_2\ra$ and $\la\tau_1,\tau_2\ra$ in which case, assuming $\lambda_1\leq \tau_1$ (so that $\tau_1+\lambda_2-(2n-1)>0$), $\sigma$ has type
$$\la\lambda_2+\tau_2-(2n-1),\tau_1+\lambda_2-(2n-1),2,\dots,2\ra.$$
\end{lemma}

\begin{proof}For this proof, we return to the see-saw dual pair 
\begin{center}
	\begin{tikzcd}
		Sp_{2a+2b}\arrow[dash]{dr}&O_{2m}\times O_{2m}\\
		Sp_{2a}\times Sp_{2b}\arrow[hookrightarrow]{u}\arrow[dash]{ur}&O_{2m}\arrow[hookrightarrow]{u}{\Delta}
	\end{tikzcd}.
\end{center}
Assume $\lambda_1<\dots<\lambda_k$, $\tau_1<\dots<\tau_l$, and let $2a=2m-\lambda_k-1$ and $2b=2m-\tau_l-1.$ By corollary \ref{convergencecorollary}, we may take $\pi_a$ and $\pi_b$ to be theta lifts $\theta(\rho_a)$ and $\theta(\rho_b)$ from unramified Arthur type representations $\rho_a$ and $\rho_b$ of $Sp_{2a}$ and $Sp_{2b}$ of types $\la\lambda_1,\dots,\lambda_{k-1}\ra$ and $\la\tau_1,\dots,\tau_{l-1}\ra$ respectively. We may assume $a\leq b.$ We have two cases to consider:
\begin{enumerate}[(1)]
\item Suppose $a+b<m$. In this case, we will show $\theta(\rho_a)\otimes\theta(\rho_b)$ weakly contains a representation theta lifted from an unramified Arthur type representation of $Sp_{2a+2b}$. We will show this representation is of the correct type using the induction hypothesis for induction from $Sp_{2a}\times Sp_{2b}$ to $Sp_{2a+2b}$.
\item Suppose $a+b\geq m.$ In this case, we show $\theta_a(\rho)\otimes\theta_b(\rho)$ weakly contains only tempered representations.
\end{enumerate}

\textbf{Case 1:} Suppose $2a+2b<2m.$ Let $\sigma$ be an unramified Arthur type representation which is weakly contained in $\pi=\Ind_{Sp_{2a}\times Sp_{2b}}^{Sp_{2a+2b}}\rho_a\boxtimes\rho_b$, and let $\theta(\sigma)$ be its (necessarily nonzero) theta lift to $O_{2m}.$ We claim 
$$\theta(\rho)\otimes\theta(\tau)\succ \theta(\sigma).$$

Let $G=Sp_{2a+2b}$ and $H=Sp_{2a}\times Sp_{2b}$. For this proof, when we restrict a representation of $O_{2m}\times O_{2m}$ to $O_{2m}$, it will be understood that we are restricting to a diagonally embedded copy of $O_{2m}$ as in the see-saw dual pair. Let $\rho=\rho_a\boxtimes \rho_b$ and $\theta(\rho)=\theta(\rho_a)\boxtimes\theta(\rho_b)$ so that $\sigma\prec \Ind_H^G \rho$ and $\theta(\rho)|_{O_{2m}}=\theta(\rho_a)\otimes\theta(\rho_b);$ hence we wish to show $\theta(\sigma)\prec \theta(\rho)|_{O_{2m}}$. Denote by $\omega_G$ the Weil representation for the dual pair $(G,O_{2m})$ and $\omega_H$ the Weil representation for $(H,O_{2m}\times O_{2m})$ so that $\omega_G|_{H\times O_{2m}}=\omega_H|_{H\times O_{2m}}.$ for each representation $\eta$ in this section, $\cH_\eta$ will denote the space on which $\eta$ acts.

Our proof will proceed as follows. We first verify that we may use Li's construction of the theta lift for $\theta(\sigma)$ and verify that it has the type predicted in lemma \ref{tensorproductofunipotentslemma}. Next, we will construct a sequence of functions on $G\times O_{2m}$ which pointwise approaches a diagonal matrix coefficient of $\omega_G\otimes\sigma$. We will then integrate this sequence over $G$ to obtain a sequence of functions on $O_{2m}$ which pointwise approaches (by Lebesgue's dominated convergence theorem) the integral of the limit, which is a matrix coefficient for $\theta(\sigma)$. It will then turn out that each of these functions of $O_{2m}$ is a diagonal matrix coefficient for $\theta(\rho)$. Finally, we will show uniform convergence on an arbitrary compact subset of $O_{2m}$. Hence we will have shown weak containment by approximation of diagonal matrix coefficients uniformly on compacta. 

Let us begin. By the induction hypothesis for induction from $Sp_{2a}\times Sp_{2b}$ to $Sp_{2a+2b}$, which we may apply since $2a+2b<2m$, the type of $\sigma$ is $\la\tau_1-2a,\dots,\tau_{l-1}-2a\ra$ (or $\la 2b-2a+1,2,\dots,2\ra$ if $k=l=2$). The largest number in the type of $\sigma$ then is $\max\{\tau_{l-1}-2a,1\}$ (or $2$ in the case $k=l=2$ and $a=b$), and since $2a+2b+2\leq 2m$ and 
$$2a+2b+\tau_{l-1}-2a=2m-\tau_l+\tau_{l-1}-1<2m$$
since $\tau_{l-1}<\tau_l$, we may use Li's construction for the theta lift of $\sigma$ by theorem \ref{convergencetheorem}.

By proposition \ref{typeofthetalift}, the theta lift $\theta(\sigma)$ of $\sigma$ to $O_{2m}$ is nonzero and an Arthur type unramifed representation of type 
$$\la \tau_1-2a,\dots,\tau_{l-1}-2a,2m-(2a+2b)-1\ra$$ 
(or $\la 2m-(2a+2b)-1,2b-2a+1,2,\dots,2\ra$ if $k=l=2$). Since $2a=2m-\lambda_k-1$ and $2b=2m-\tau_k-1$, this produces a final type of
$$\la \tau_1+\lambda_k-(2m-1),\dots,\tau_k+\lambda_k-(2m-1)\ra$$
(or $\la \lambda_2+\tau_2-(2m-1),\lambda_2+\tau_1-(2m-1),2,\dots,2\ra$ if $k=l=2$). Note that the terms $\lambda_i+\tau_j$ do not appear for $i<k$ because 
$$\lambda_i+\tau_j<2a+1+\tau_l=2a+2m-2b\leq 2m$$
since $a\leq b$, and so $\theta(\sigma)$ has the type predicted in lemma \ref{tensorproductofunipotentslemma}. 

By corollary \ref{convergencecorollary}, we may also use Li's construction of the theta lift for $\theta(\rho)$. A diagonal matrix coefficient for $\theta(\sigma)$ associated to a pure tensor is a function of the form
\begin{equation}\label{eqn422}x\mapsto\int_{G}\la \omega_G(xg)\phi,\phi\ra\la \sigma(g)v,v\ra dg\end{equation}
for $x\in GL_n$, $\phi\in\cH_{\omega_G}$, and $v\in \cH_\sigma.$ We wish to approximate this with diagonal matrix coefficients of $\theta(\rho)$ uniformly on compact subsets of $O_{2m}$.

Recall by lemma F.1.3 of \cite{Bekka2008} that in order to show weak containment, we only need to approximate the diagonal matrix coefficients of $\theta(\sigma)$ associated to vectors in some total subset of $\cH_{\theta(\sigma)}.$ Recall from section 2.3 that $\cH_{\theta(\sigma)}=(\cH_{\omega_G}\otimes \cH_{\sigma})/R$ for some subspace $R\subset \cH_{\omega_G}\otimes \cH_{\sigma}$. Let $K$ be a maximal compact subgroup of $G$. The set of pure tensors $\phi\otimes v$ with $\phi\in \cH_{\omega_G}$ a unit $K$-finite vector and $v\in \cH_{\sigma}$ $v$ a unit vector in some $K$-type forms a total set in $\cH_{\omega_G}\otimes \cH_{\sigma}$, and hence the set of $(\phi\otimes v)R$ with $\phi,v$ as above is total in $\cH_{\theta(\sigma)}$. Thus we assume $\phi$ is a $K$-finite unit vector and $v$ is a unit vector in some $K$-type $\mu$.

Let $C_1\subset C_2\subset\dots$ be compact subsets of $G=Sp_{2a+2b}$ such that $G=\bigcup_i C_i.$ By proposition F.1.4 of \cite{Bekka2008}, since $\sigma\prec \pi=\Ind_H^G\rho$ is irreducible and $v$ is a unit vector, for each $i$ there exists $\zeta_i\in \cH_\pi$ such that
\begin{equation}\label{eqn421}|\la \sigma(g)v,v\ra-\la \pi(g)\zeta_i,\zeta_i\ra|<\frac{1}{2i}\end{equation}
for all $g\in C_i.$ We wish to refine this estimate.

For $f\in C_c^\infty(G)$, $w\in \cH_\rho$, and $g\in G$, let
$$\xi_{f,w}(g)=\int_{H} f(gh)\rho(h)wdh.$$

By the proof of lemma E.1.3 in \cite{Bekka2008}, the set
$$\sH_\pi=\{\xi_{f,w}:f\in C_c(G),w\in \cH_\rho\}$$
is dense in $\cH_\pi$. Then for each $i$, we may find $f_i\in C_c(G)$ and $w_i\in \cH_\rho$ such that 
$$\|\zeta_i-\xi_{f_i,w_i}\|<\frac{1}{8i}.$$
Observe that since $\zeta_i$ is a unit vector, it follows that $|1-\|\xi_{f_i,w_i}\||<\frac{1}{8i}$. This then implies
$$\|\zeta_i-\|\xi_{f_i,w_i}\|^{-1}\xi_{f_i,w_i}\|<\frac{1}{4i}.$$
Note that if we let $v_i=\|\xi_{f_i,w_i}\|^{-1}w_i$, then 
$$\|\xi_{f_i,w_i}\|^{-1}\xi_{f_i,w_i}=\xi_{f_i,\|\xi_{f_i,w_i}\|^{-1}w_i}=\xi_{f_i,v_i}$$
is a unit vector. Since $\pi$ is unitary, we then have
\begin{align*}
|\la\pi(g)\zeta_i,\zeta_i\ra-\la\pi(g)\xi_{f_i,v_i},\xi_{f_i,v_i}\ra|&\leq|\la\pi(g)\zeta_i,\zeta_i-\xi_{f_i,v_i}\ra|+|\la\pi(g)(\zeta_i-\xi_{f_i,v_i}),\xi_{f_i,v_i}\ra|\\
&\leq \|\zeta_i\|\cdot\|\zeta_i-\xi_{f_i,v_i}\|+\|\zeta_i-\xi_{f_i,v_i}\|\cdot\|\xi_{f_i,v_i}\|\\
&<\frac{1}{2i}.
\end{align*}
It follows from (\ref{eqn421}) then that
\begin{equation}\label{eqn423}|\la\sigma(g)v,v\ra-\la\pi(g)\xi_{f_i,v_i},\xi_{f_i,v_i}\ra|<\frac{1}{i}\end{equation}
for all $g\in C_i$. Recall that $v$ belongs to a $K$-type $\mu$. As in section \ref{restrictionsphericalunipotent}, we may project $\xi_{f_i,v_i}$ to the $\mu$-isotypic subspace of $\cH_\pi$ by convolving with the central idempotent $e_\mu$ over $K$, and the projected vectors will still have the convergence (\ref{eqn423}). It's easy to show by the Fubini-Tonelli theorem that such a projection is still a vector in the set $\sH_\pi$. Without changing notation then, we assume $\xi_{f_i,v_i}$ to belong to the $K$-type $\mu$, and since projecting a vector does not increase its norm, we have $\|\xi_{f_i,v_i}\|\leq 1$.

For each $g\in G$, $x\in O_{2m}$, and each $i$, let
$$F_i(g,x)=\la \omega_G(xg)\phi,\phi\ra\la\pi(g)\xi_{f_i,v_i},\xi_{f_i,v_i}\ra.$$
We have so far shown 
$$F_i(g,x)\rightarrow F(g,x)=\la\omega_G(xg)\phi,\phi\ra\la\sigma(g)v,v\ra$$
as $i\rightarrow \infty$ pointwise in $G\times O_{2m}.$ 

We now seek to establish dominated convergence for the sequence of integrals of the $F_i$ over $G$. From the proof of theorem \ref{convergencetheorem}, the induction from $Sp_{2a}\times Sp_{2b}$ to $Sp_{2a+2b}$ case, and lemma 6.2 of \cite{Howe2010}, we know that $\pi=\Ind_{H}^G \rho$ has $K$-finite matrix coefficients bounded in absolute value by a multiple of $\Xi_{G}^q$ for some $q$ and that $\omega_G|_{G}$ has matrix coefficients bounded in absolute value by a multiple of $\Xi_{G}^p$ for some $p$ such that $p+q>2$. It follows for fixed $x\in O_{2m}$ that $F_i(g,x)$ has matrix coefficients bounded by a constant multiple of $\Xi_{G}^{p+q}$, which is an $L^1(G)$ function. The constant multiple in question depends on $\|\xi_i\|\leq 1$, so it is uniformly bounded in $i$. By Lebesgue's dominated convergence theorem, we conclude
$$\int_G F_i(g,x)dg\rightarrow \int_G F(g,x)dg$$
as $i\rightarrow \infty$ for each $x\in O_{2m}.$

Recall that $x\mapsto \int_G F(g,x)dg$ is a matrix coefficient for $\theta(\sigma)$. We will show that $x\mapsto \int_G F_i(g,x)dg$ is in fact a matrix coefficient for $\theta(\rho)|_{O_{2m}}.$ First, observe that $F_i(g,x)$ is a matrix coefficient for the $G\times O_{2m}$-representation $\omega_G\otimes \Ind_{H}^G\rho.$ By proposition E.2.5 of \cite{Bekka2008} and using the see-saw dual pair identity,
$$\omega_G\otimes\Ind_{H}^G\rho\cong \pi':=\Ind_H^G(\omega_G|_{H\times O_{2m}}\otimes \rho)=\Ind_H^G(\omega_H|_{H\times O_{2m}}\otimes \rho).$$
The isomorphism $\omega_G\otimes\pi=\pi'$ is given by the unitary map $U:\cH_{\omega_G}\otimes\cH_{\pi}\rightarrow \cH_{\pi'}$ defined by
$$U(\theta\otimes \xi)(g)=\omega_G(g^{-1})\theta\otimes \xi(g).$$
We have
$$\int_G F_i(g,x)dg=\int_G\la \pi'(g,x)U(\phi\otimes \xi_{f_i,v_i}),U(\phi\otimes\xi_{f_i,v_i})\ra dg.$$

Observe for each $x\in O_{2m}$ and $g\in G$ (cf. proof of F.3.5 in \cite{Bekka2008}),
\begin{align*}
&\la \pi'(g,x)U(\phi\otimes \xi_{f_i,v_i}),U(\phi\otimes\xi_{f_i,v_i})\ra\\
&=\int_{G/H}\la \pi'(g,x)U(\phi\otimes \xi_{f_i,v_i})(y),U(\phi\otimes \xi_{f_i,v_i})(y)\ra dy\\
&=\int_{G/H}\la \omega_G(y^{-1}gx)\phi,\omega_G(y^{-1})\phi\ra\cdot \la\xi_{f_i,v_i}(g^{-1}y),\xi_{f_i,v_i}(y)\ra dy\\
&=\int_{G/H}\int_H\int_Hf_i(g^{-1}yh)\overline{f_i(yk)}\la\omega_G(gx)\phi,\phi\ra\la\rho(k^{-1}h)v_i,v_i\ra dh\,dk\,dy\\
&=\int_G\int_H f_i(g^{-1}yh)\overline{f_i(y)}\la\omega_G(gx)\phi,\phi\ra\la\rho(h)v_i,v_i\ra dh\,dy.
\end{align*}
Integrating over $G$ and using the Fubini-Tonelli theorem, for each $x\in O_{2m}$ we have 
\begin{align*}
\int_GF_i(g,x)dg&=\int_G\int_G\int_H f_i(g^{-1}yh)\overline{f_i(y)}\la\omega_G(gx)\phi,\phi\ra\la\rho(h)v_i,v_i\ra dh\,dy\,dg\\
&=\int_H\int_G\int_G f_i(g^{-1})\overline{f_i(y)}\la\omega_G(yhgx)\phi,\phi\ra\la\rho(h)v_i,v_i\ra dg\,dy\,dh\\
&=\int_H\left\la \int_G f_i(g^{-1})\omega_G(hgx)\phi\,dg,\int_G f_i(y)\omega_G(y^{-1})\phi\,dy\right\ra\la\rho(h)v_i,v_i\ra dh\\
&=\int_H \la \omega_H(xh)\Phi_i,\Phi_i\ra\la\rho(h)v_i,v_i\ra dh
\end{align*}
where $\Phi_i=\int_G f_i(g^{-1})\omega_G(g)\phi \,dg$, and we have used the facts that $x$ commutes with $G$, $\omega_G$ is unitary, and $\omega_G|_{H\times O_{2m}}=\omega_H|_{H\times O_{2m}}.$ We recognize the last line as a diagonal matrix coefficient for $\theta(\rho)|_{O_{2m}}.$ Therefore we have constructed a sequence of matrix coefficients of $\theta(\rho)|_{O_{2m}}$ converging pointwise to the matrix coefficient (\ref{eqn422}) of $\theta(\sigma)$. 

We must show that the convergence $\int_G F_i(g,x)dg\rightarrow \int_G F(g,x)dg$ is uniform on compact subsets of $O_{2m}$. Let $Q\subset O_{2m}$ be compact. Recall that $\la\pi(g)\xi_{f_i,v_i},\xi_{f_i,v_i}\ra$ is bounded by a multiple of $\Xi_{G}^p$ and converges pointwise to $\la \sigma(g)v,v\ra$. By Lebesgue's dominated convergence theorem, for any $r>\frac{2}{p}$
$$\|\la \pi(\cdot)\xi_{f_i,v_i},\xi_{f_i,v_i}\ra-\la\sigma(\cdot)v,v\ra\|_{L^r(G)}\rightarrow 0$$
as $i\rightarrow \infty.$ Also recall that $\omega_G|_{G}$ has matrix coefficients bounded by $\Xi_G^q$ where $p+q>2$, and so it has matrix coefficients in $L^r(G)$ for $r>\frac{2}{q}$. Then $\frac{p}{2}+\frac{q}{2}>1$, and hence there exists $p_0>\frac{2}{p}$ and $q_0>\frac{2}{q}$ with $\frac{1}{p_0}+\frac{1}{q_0}=1$. By H\"{o}lder's inequality,
\begin{align*}
&\left|\int_G F_i(g,x)-F(g,x)dg\right|\leq \int_G |\la\omega_G(xg)\phi,\phi\ra(\la \pi(g)\xi_{f_i,v_i},\xi_{f_i,v_i}\ra-\la \sigma(g)v,v\ra)|dg\\
&\hspace{3.5cm}\leq \|\la\omega_G(x\cdot)\phi,\phi\ra\|_{L^{p_0}(G)}\cdot\|\la\pi(\cdot)\xi_{f_i,v_i},\xi_{f_i,v_i}\ra-\la\sigma(\cdot)v,v\ra\|_{L^{q_0}(G)}.
\end{align*}
In the last line above, we may take the supremum over all $x\in Q$ to obtain a constant multiple in the first factor, and then the second factor tends to $0$. This shows $\int_G F_i(g,x)dg\rightarrow \int_G F(g,x)dg$ uniformly for $x\in Q$ as $i\rightarrow \infty.$ Therefore we have approximated the matrix coefficients of $\theta(\sigma)$ in a total subset by those of $\theta(\rho)|_{O_{2m}}$ uniformly on compact subsets of $O_{2m}$, and we conclude
$$\theta(\rho)|_{O_{2m}}\succ \theta(\sigma).$$
This proves lemma \ref{tensorproductofunipotentslemma} in this case.

\textbf{Case 2:} Suppose $a+b\geq m$. We need to show $\theta(\rho_a)\otimes\theta(\rho_b)$ has matrix coefficients in $L^r$ for any $r>2$. By corollary \ref{matrixcoefficientsdecayratethetalifted}, $\theta(\rho_a)$ has matrix coefficients in $L^p$ for $p>\frac{4m-4}{2m-(2m-2a-1)-1}=\frac{2m-2}{a}.$ Likewise, $\theta(\rho_b)$ has matrix coefficients in $L^q$ for $q>\frac{2m-2}{b}.$ Note that 
$$\left(\frac{2m-2}{a}\right)^{-1}+\left(\frac{2m-2}{b}\right)^{-1}=\frac{a+b}{2m-2}\geq\frac{1}{2}\cdot\frac{m}{m-1}>\frac{1}{2}.$$
Let $r>2$. Then there exists $p>\frac{2m-2}{a}$ and $q>\frac{2m-2}{b}$ such that $\frac{1}{p}+\frac{1}{q}=\frac{1}{r}.$ Let $f(g)$ and $f'(g)$ be matrix coefficients of $\theta(\rho_a)$ and $\theta(\rho_b)$ respectively. By H\"{o}lder's inequality,
$$\|ff'\|_{L^r(O_{2m})}\leq \|f\|_{L^p}\|f'\|_{L^q}\|.$$
Thus the matrix coefficients of $\theta(\rho)\otimes\theta(\tau)$ are in $L^r$ for all $r>2$, which makes the tensor product tempered. This concludes the proof of lemma \ref{tensorproductofunipotentslemma} \end{proof}

\section{The induction problems}
This chapter involves the most explicit calculations in this paper in that we realize all of the groups involved explicitly by their actions on a vector space with a skew-symmetric bilinear form. This is necessary as the general framework for the sort of analysis involved in this section has not yet been fleshed out enough for our purposes. 
\subsection{Induction from $Sp_{2a}\times Sp_{2b}$ to $Sp_{2a+2b}$}\label{inductionspherical}
In this section, we will prove the following sub-case of proposition \ref{weakprop}.
\begin{lemma}\label{inductionsphericallemma}
Let $\rho_a$ and $\rho_b$ be unramified Arthur type representations of $Sp_{2a}$ and $Sp_{2b}$ respectively, and assume $a\leq b$. If $\la\lambda_i\ra$ is the type of $\rho_b$, then $\pi=\Ind_{Sp_{2a}\times Sp_{2b}}^{Sp_{2a+2b}}\rho_a\boxtimes\rho_b$ weakly contains an unramified Arthur type representation of type $\la \lambda_i-2a\ra$, except in the case where $\rho_a$ and $\rho_b$ are both one-dimensional, in which case the type is $\la 2b-2a+1,2,\dots,2\ra.$
\end{lemma}

Before we start into the proof properly, let us give some motivation for how we study this induction problem. Let $G=Sp_{2a+2b}$ and $H=Sp_{2a}\times Sp_{2b}$. The space $G/H$ is a spherical variety, and the space $L^2(G/H)$ forms a representation of $G$ via left translation isomorphic to $\Ind_H^G\textbf{1}_H$. In \cite{Sakellaridis2012}, Sakellaridis and Venkatesh show how the spectrum of $L^2(G/H)$ is controlled by so called ``boundary degenerations." This is accomplished by constructing a ``wonderful compactification" of $G/H$, i.e. a compact space $\overline{X}$ into which $G/H$ embeds as an open dense sub-variety with some desirable properties. By studying the boundary points of $G/H$ in this compactification, they are able to decompose a part of $L^2(G/H)$ in terms of parabolically induced representations.

Their theory is not easy to apply to our setting where instead of $\Ind_H^G\textbf{1}_H$, we are considering $\Ind_H^G\rho_a\boxtimes\rho_b$. However, philosophically we can take the same approach. In this case, we may identify $G/H$ with the Grassmannian of $2a$-dimensional nondegenerate spaces in a $2a+2b$-dimensional symplectic space $V$. The full Grassmannian of $2a$-dimensional subspaces of $V$ is a compact space into which $G/H$ embeds as an open dense sub-variety.

Our proof will use the following strategy. We may take a sequence of points $\{V_i\}$ in $G/H$ realized as $2a$-dimensional nondegenerate subspaces in $V$ and have it converge to a $2a$-dimensional isotropic subspace $W$ of $V$, the $G$-stabilizer of which is a parabolic subgroup $P$ of $G$. The representation $\rho_a\boxtimes \rho_b$ will be realized on $Sp(V_i)\times Sp(V_i^{\perp})$ for each $i$ in a natural way. We shall study the behavior of the matrix coefficients of the representation $\rho_a\boxtimes\rho_b$ as $i\rightarrow \infty$. 

In the end, we obtain a description of the following form. There is a fixed group $M'=(Sp(V_i)\times Sp(V_i^\perp))\cap P$ which does not depend on $i$ and is contained in a Levi subgroup $M$ of $P$. We will show that if a sequence of elements $\{x_i\}\subset Sp(V)/M'$ with $x_i\in Sp(V_i)/M'$ for each $i$ converges in $Sp(V)/M'$ as $i\rightarrow \infty$, then $\lim_{i\rightarrow \infty} x_i$, in some sense, lies in the unipotent radical $N$ of $P$, and the matrix coefficients of $\rho_a\boxtimes\rho_b$ approach the matrix coefficients of the trivial representation of $N$ for such a sequence as $i\rightarrow \infty$. By integrating these matrix coefficients and addressing some analytical difficulties, we will eventually obtain a result of the form
$$\Ind_H^G\rho_a\boxtimes\rho_b\succ \Ind_{M'N}^G(\rho_a\boxtimes\rho_b)|_{M'}\cdot\textbf{1}_N.$$
By inducing in stages through the parabolic $P$, we obtain a representation of the form (\ref{stdrep}), from which we are able to read the $SL_2$-type. This same technique will be used with appropriate modification in section \ref{inductionexceptionalprop}. Let us begin.
\begin{proof}Let $n=a+b$, and let $V$ be a $2n$-dimensional vector space over $\nQ_p$ with a non-degenerate skew-symmetric bilinear form $\la \cdot,\cdot \ra$. We have a left action of $GL(V)\cong GL_{2n}$ on this space, and the set of elements which stabilize $\la \cdot,\cdot \ra$ is $Sp(V)\cong Sp_{2n}.$ Choose a symplectic basis $\{e_1,\dots,e_n,f_1,\dots,f_n\}$ with $\la e_i,e_j\ra= \la f_i,f_j\ra=0$ and $\la e_i,f_j\ra=\delta_{ij}.$ For $t\in \nQ_p^\times,$ let 
$$v_{2i-1}^\pm(t)=\frac{1}{2}(e_{2i-1}\pm t^2f_{2i})\text{ and } v_{2i}^\pm(t)=\frac{1}{2}(e_{2i}\mp t^2f_{2i-1})$$
for $i=1,\dots,a$ so that for all $i,j$,
\begin{align*}
\la v_i^\pm(t),v_j^\mp(t)\ra&=0,\\
\la v_{2i-1}^\pm(t),v_{2j-1}^\pm(t)\ra&=0,\\
\la v_{2i}^\pm(t),v_{2j}^\pm(t)\ra&=0,\\
\la v_{2i-1}^\pm(t),v_{2j}^\pm(t)\ra&=\begin{cases}t^2&\text{if }i=j\\0&\text{if }i\not=j\end{cases}.
\end{align*}

Let $U$ be the $2b-2a$-dimensional nondegenerate space
$$U=\begin{cases}\{0\}&\text{if }a=b\\\spn\{e_{2a+1},\dots,e_n,f_{2a+1},\dots,f_n\}&\text{if }a<b\end{cases},$$
and for each $t\not=0$, let 
\[V_t^\pm=\spn\{v_1^\pm(t),\dots,v_{2a}^\pm(t)\}\]
The spaces $V_t^+$ and $V_t^-$ are nondegenerate symplectic spaces of dimension $2a$, and $(V_t^+)^\perp=V_t^-\oplus U$. For each $t\not=0$, we have the symplectic decomposition
$$V=V_t^+\oplus V_t^-\oplus U.$$
Let 
$$H_t=Sp(V_t^+)\times Sp(V_t^-\oplus U)\cong Sp_{2a}\times Sp_{2b}$$ 
be the set of elements of $Sp(V)$ which stabilize $V_t^+$ (and hence $V_t^-\oplus U$).

The groups $H_t$ are all isomorphic and conjugate to $H_1$ via the element $g_t\in Sp(V)$ defined by 
\begin{align*}
g_te_i&=\begin{cases}\frac{1}{t}e_i&\text{if }i\leq 2a\\e_i&\text{if }i>2a\end{cases},\\
g_tf_i&=\begin{cases}tf_i&\text{if }i\leq 2a\\f_i&\text{if }i>2a\end{cases}.
\end{align*}
Denote $\rho=\rho_a\boxtimes \rho_b$, and assume that $\rho$ is a representation of $H_1$. Then define $\rho_t(h_t)=\rho(g_t^{-1}h_tg_t)$ for $h_t\in H_t.$ Then $\rho_t$ for $t\not=0$ is a representation of $H_t$ isomorphic to $\rho=\rho_1.$

Let 
$$W=\spn\{e_1,\dots,e_{2a}\}\text{ and }W'=\spn\{f_1,\dots,f_{2a}\}.$$ 
These spaces are isotropic. Let $P$ be the set of elements of $Sp(V)$ which stabilize $W$. This is a parabolic subgroup with Levi subgroup $M=GL(W)\times Sp(U)$ and unipotent radical $N$ which acts by the identity endomorphism on $W$, $(W\oplus U)/W$, and $V/(W\oplus U).$ Observe that we have the isotropic splitting
$$V_t^+\oplus V_t^-=W\oplus W'$$
for each $t\not=0$ and that as $t\rightarrow 0$, $v_i^\pm(t)\rightarrow e_i$. 

\textbf{Claim:} For every $t\not=0$, the subgroup of elements of $H_t$ which stabilize $W$ is the same subgroup for every $t$, and it is of the form $Sp(W)\times Sp(U)$, where $Sp(W)$ is the stabilizer of some symplectic form on $W$. 

Observe that we have a natural surjection $\phi_t:V_t^-\oplus U\rightarrow V_t^+$ given by $\phi_t(v_i^-(t))=v_i^+(t)$, $\phi(U)=0$. Note $\phi_t|_{V_t^-}$ is an isomorphism with graph $W$ with the property $\la w,w'\ra=\la \phi(w'),\phi(w)\ra$ for all $w,w'\in V_t^-.$ The graph of $\phi_t$ is $W\oplus U.$ The isomorphism $\phi_t|_{V_t^-}$ defines a symplectic form on $W$ via
$$(w+\phi_t(w),w'+\phi_t(w'))_t=\la w,w'\ra.$$
The set of elements of $H_t$ which stabilize $\phi_t$ are precisely the elements of $H_t$ which stabilize $W$. Such elements must stabilize $V_t^-$ and hence $U$, and so this stabilizer is contained in $Sp(V_t^+)\times Sp(V_t^-)\times Sp(U).$ Clearly elements of $Sp(U)$ stabilize $W$ since they act by the identity, and the elements of $Sp(V_t^+)\times Sp(V_t^-)$ stabilizing the isomorphism $\phi|_{V_t^-}$ are those that stabilize $(\cdot,\cdot)_t$, which we denote $Sp(W)$. Note $Sp(W)$ is diagonally embedded in $Sp(V_t^+)\times Sp(V_t^-).$

A priori, it seems like $Sp(W)$ should depend on $t$, but in fact it does not. Every element of $H_t$ which stabilizes $W$ is of the form $g_thg_t^{-1}$, where $h$ is an element of $H_1$ that stabilizes $W$. Since $g_t$ acts by a scalar on $W$ and $h$ stabilizes $W$, $g_thg_t^{-1}W=hW$. On the other hand, $Sp(W)\subset GL(W)$ since the only elements of $H_t$ which act trivially on $W$ are elements of $Sp(U)$ (so $N$ and $H_t$ have trivial intersection), and $GL(W)$ also stabilizes $W'$, so we have 
$$g_thg_t^{-1}W'=t h \frac{1}{t}W'=hW'.$$ 
Thus $Sp(W)$ is point-wise invariant under conjugation by $g_t$ and does not depend on $t$.

The stabilizer of $W$ in $H_t$ is seen to be $Sp(W)\times Sp(U)$, and moreover we have 
$$H_t\cap H_{t'}=H_t\cap P=Sp(W)\times Sp(U)\subset M=GL(W)\times Sp(U)$$ 
for $t\not= t'.$ This proves the claim.

Let $H_0=(Sp(W)\times Sp(U))\ltimes N.$ We now claim
\begin{equation}\label{eqn510}\Ind_{H_1}^{Sp(V)}\rho\succ \Ind_{H_0}^{Sp(V)}\rho|_{Sp(W)\times Sp(U)}\cdot \textbf{1}_N.\end{equation}

In order to show this, we will prove the following. 

\begin{lemma}\label{deformationlemma}For every $t\in \nQ_p^\times$, there exists a continuous bijection $\psi_t$ from an open dense subset of the center $Z(N)$ of $N$ onto an open dense subset of $(Sp(V_t^+)\times Sp(V_t^-))/Sp(W)$ and $\Psi_t$ from an open dense subset of $N/Z(N)$ onto an open dense subset of 
$$H_t/(Sp(V_t^+)\times Sp(V_t^-)\times Sp(U))=Sp(V_t^-\oplus U)/(Sp(V_t^-)\times Sp(U))$$
with the following properties:
\begin{enumerate}[(1)]
	\item For every $z\in Z(N)$, $z$ is in the domain of $\psi_t$ for all but finitely many $t$. Likewise, for any $n\in N/Z(N)$, $n$ is in the domain of $\Psi_t$ for all but finitely many $t$.
	\item There is a natural continuous bijection $B:(Sp(V_t^+)\times Sp(V_t^-))/Sp(W)\rightarrow Sp(V_t^+)$ given by 
	$$B((x,y)Sp(W))=B((x\cdot\phi(y^{-1}),1)Sp(W))=x\cdot\phi(y^{-1}).$$
	 Under this bijection, $\lim_{t\rightarrow 0}B(\psi_t(z))=z$ for all $z\in Z(N)$, where the limit is taken in $Sp(V).$
	\item There is a natural identification between open dense subsets
	$$H_t/(Sp(V_t^+)\times Sp(V_t^-)\times Sp(U))\leftrightarrow \Hom(V_t^-,U)\subset \Hom(W\oplus W',U),$$
	and also $N/Z(N)\cong \Hom(W,U)\subset \Hom(W\oplus W',U)$. Likewise there is an identification between open dense subsets
	$$G/(Sp(W\oplus W')\times Sp(U))\leftrightarrow \Hom(W\oplus W',U),$$
	Under this identification, $\lim_{t\rightarrow 0}\Psi_t(n)=n$ for all $n\in N/Z(N)$ (limit in $G/(Sp(W\oplus W')\times Sp(U))$).
	\item Given a Haar measure $dp$ on $H_0$, there exist Haar measures $dh$ on $Sp(W)\times Sp(U)$, $dz$ on $Z(N)$, and $d\mu_0$ on $N/Z(N)$ such that
	$$\int_{H_0}f(p)dp=\int_{N/Z(N)}\int_{Z(N)}\int_{Sp(W)\times Sp(U)}f(nzh)dz\, d\mu_0(n)\, dh$$
	for any $f\in C_c(G)$. There then exists for each $t\not=0$ invariant measures $dz_t$ on $(Sp(V_t^+)\times Sp(V_t^-))/Sp(W)$ and $d\mu_t$ on $H_t/(Sp(V_t^+)\times Sp(V_t^-)\times Sp(U))$ which define a Haar measure $dh_t$ on $H_t$ in the following way. Let $K\subset H_t$ be compact, and let $\chi_K$ be the characteristic function of $K$. Then we define
	\begin{align*}dh_t(K)=\int_{H_t/(Sp(V_t^+)\times Sp(V_t^-)\times Sp(U))}&\int_{(Sp(V_t^+)\times Sp(V_t^-))/Sp(W)}\\&\int_{Sp(W)\times Sp(U)}\chi_K(nzh)dh\,dz_t(z)\,d\mu_t(n).\end{align*}
	This $dh_t$ has the property that,
	\begin{equation}\label{eqn512}\lim_{t\rightarrow 0}\int_{H_t}f(x)dh_t(x)=\int_{H_0} f(p)dp.\end{equation}
	 for any $f\in C_c(G)$
\end{enumerate}\end{lemma}
The measures in the last part will be constructed using $\psi_t$ and $\Psi_t$. A dense subset of matrix coefficients for the right side of (\ref{eqn510}) will be expressed as an integral over $G$ of the right side of (\ref{eqn512}), while the left side of (\ref{eqn512}) will be matrix coefficients expressed as an integral over $G$ of the left side of (\ref{eqn510}) for an appropriate $f\in C_c(G)$.

\textbf{Constructing $\psi_t$:} We begin now the proof of lemma \ref{deformationlemma}. Let us begin by constructing $\psi_t$, which is a function from an open dense subset of $Z(N)\cong \Sym_{2a}(\nQ_p)$ to an open dense subset of $(Sp(V_t^+)\times Sp(V_t^-))/Sp(W).$ Each element of $(Sp(V_t^+)\times Sp(V_t^-))/Sp(W)$ has a unique representative of the form $(h,1)$, where $h\in Sp(V_t^+)$, which produces the natural continuous bijection
$$B:Sp(V_t^+)\leftrightarrow (Sp(V_t^+)\times Sp(V_t^-))/Sp(W).$$
 It then suffices to give a representative of $\psi_t(z)$ in $Sp(V_t^+)$ for each $z$ in an open dense subset of $Z(N)$. 
 
 Consider first the elements $\alpha n_{i,i}\in Z(N)$, $\alpha \in \nQ_p$ ($1\leq i\leq 2a$), which act trivially each $e_j$ and on each $f_j$ except for $f_i$, which is mapped to $f_i+\alpha e_i.$ Define $\psi_t(\alpha n_{2i,2i})$ ($i\leq a$) to act trivially on $V_t^-\oplus U$, and trivially on each $v_j^+(t)$ except $v_{2i-1}^+(t)$ and let 
 $$\psi_t(\alpha n_{2i,2i})v_{2i-1}^+(t)=v_{2i-1}^+(t)+2\alpha t^2v_{2i}^+(t).$$ 
 This is a lower-triangular element of 
 $$Sp(\spn\{v_{2i-1}^+(t),v_{2i}^+(t)\})\subset Sp(V_t^+).$$ 
Then given the identities
$$e_{i}=v_i^+(t)+v_i^-(t),\quad f_{2i-1}=\frac{-1}{t^2}(v_{2i}^+(t)-v_{2i}^-(t)),\quad f_{2i}=\frac{1}{t^2}(v_{2i-1}^+(t)-v_{2i-1}^-(t)),$$
we see that $\psi_t(\alpha n_{2i,2i})$ fixes all $e_j$ and $f_j$ except $e_{2i-1}$ and $f_{2i}$, and 
$$\psi_t(\alpha n_{2i,2i})e_{2i-1}=e_{2i-1}+2\alpha t^2v_{2i}^+(t),\quad \psi_t(\alpha n_{2i,2i})f_{2i}=f_{2i}+\alpha e_{2i}+\alpha t^2f_{2i-1}.$$
Observe that as $t\rightarrow 0$, $\psi_t(\alpha n_{2i,2i})\rightarrow \alpha n_{2i,2i}.$ We analagously define $\psi_t(\alpha n_{2i-1,2i-1})$ by switching the role of $2i$ and $2i-1$ in the formula above and changing the sign on the $2at^2$ to negative and also get $\psi_t(\alpha n_{2i-1,2i-1})\rightarrow \alpha n_{2i-1,2i-1}$ as $t\rightarrow 0.$

Next, consider the maps $\alpha n_{i,j}$, $i\not=j$, $i,j\leq 2a$, which fix all $e_k$ and all $f_k$ except $f_i,f_j$, which are mapped respectively to $f_i+\alpha e_j$ and $f_j+\alpha e_i$. Define $\psi_t(\alpha n_{2i,2j})$ to fix $V_t^-\oplus U$ and all $v_i^+(t)$, except for $v_{2i-1}^+(t)$ and $v_{2j-1}^+(t)$, which are defined as
\begin{align*}\psi_t(\alpha n_{2i,2j})v_{2i-1}^+(t)&=v_{2i-1}^+(t)+2\alpha t^2v_{2j}^+(t)\\
\psi_t(\alpha n_{2i,2j})v_{2j-1}^+(t)&=v_{2j-1}^+(t)+2\alpha t^2v_{2i}^+(t).\end{align*}
This is an element of the unipotent radical of a Borel subgroup of
$$Sp(\spn\{v_{2i-1}^+(t),v_{2i}^+(t),v_{2j-1}^+(t),v_{2j}^+(t)\})\subset Sp(V_t^+).$$ 
The actions on the non-fixed $e_k,f_k$ are
\begin{align*}
\psi_t(\alpha n_{2i,2j})e_{2i-1}&=e_{2i-1}+2\alpha t^2v_{2j}^+(t),\\
\psi_t(\alpha n_{2i,2j})f_{2i}&=f_{2i}+\alpha e_{2j}+\alpha t^2f_{2j},
\end{align*}
and similar for $e_{2j-1}$ and $f_{2j}$. We analagously define $\psi_t(\alpha n_{2i-1,2j-1})$, $\psi_t(\alpha n_{2i,2j-1})$, and $\psi_t(\alpha n_{2i-1,2j})$, exchanging $2i$ with $2i-1$, exchanging $2j$ with $2j-1$, and switching signs of terms in the formulas as needed, except for the case $\psi_t(\alpha n_{2i-1,2i})$, for which the above construction does not produce an element of determinant $1$. 

We define the map $\psi_t(\alpha n_{2i-1,2i})$ to fix all $v_{k}^\pm(t)$ except 
\begin{align*}
\psi_t(\alpha n_{2i-1,2i})v_{2i-1}^+(t)&=(1+2\alpha t^2)v_{2i-1}^+(t),\\
\psi_t(\alpha n_{2i-1,2i})v_{2i}^+(t)&= (1+2\alpha t^2)^{-1}v_{2i}^+(t).
\end{align*} 
This is an element of the standard maximal torus of $Sp(V_t^+)$, and it is defined when $\alpha\not= -\frac{1}{2t^2}.$ Then the actions on the non-fixed $e_j,f_j$ are
\begin{align*}
\psi_t(\alpha n_{2i-1,2i})e_{2i-1}&=e_{2i-1}+2\alpha t^2v_{2i-1}^+(t),\\
\psi_t(\alpha n_{2i-1,2i})e_{2i}&=e_{2i}-\frac{2\alpha t^2}{1+2\alpha t^2}v_{2i}^+(t),\\
\psi_t(\alpha n_{2i-1,2i})f_{2i-1}&=f_{2i-1}+\frac{\alpha }{1+2\alpha t^2}e_{2i}-\frac{\alpha t^2}{1+2\alpha t^2}f_{2i-1},\\
\psi_t(\alpha n_{2i-1,2i})f_{2i}&=f_{2i}+\alpha e_{2i-1}+\alpha t^2f_{2i}.
\end{align*}
Once again as $t\rightarrow 0$, $\psi_t(\alpha n_{2i-1,2i})\rightarrow \alpha n_{2i-1,2i}.$	

Finally, all of the elements of $Z(N)$ are of the form $\sum_{i\geq j}\alpha _{i,j}n_{i,j}.$ Fix an ordering $(\iota_k)$, $k=1,\dots, a(2a+1)$ of the set of indices $\cI$ of pairs $(i,j)$ with $i\leq j$, $i,j=1,\dots, 2a.$ Then we define 
$$\psi_t\left(\sum_{i\geq j}\alpha _{i,j}n_{i,j}\right)=\prod_k \psi_t(\alpha _{\iota_k}n_{\iota_k}).$$
The ordering of the indices is necessary because $Sp(V_t^+)$ is not abelian whereas $Z(N)$ is. Then $\psi_t$ is defined on a dense subset of $Z(N)$ (the set of elements whose $\alpha_{2i-1,2i}n_{2i-1,2i}$ terms satisfy $\alpha_{2i-1,2i}\not=-\frac{1}{2t^2}$), and we claim its image is a dense subset of $Sp(V_t^+)$ for an appropriate ordering $(\iota_k)$. 

Observe that the product of the elements $\psi_t(\alpha n_{2i-1,2i})$ form a maximal torus $T_t^+$ of $Sp(V_t^+)$. Such a torus is contained in a Borel subgroup $B_t^+$. The rest of the elements are upper triangular or lower triangular, so they lie in the unipotent radical $N_t^+$ of $B_t^+$ (for an appropriate choice of a Borel subgroup) or the unipotent radical $\overline{N_t^+}$ of the opposite Borel subgroup $\overline{B_t^+}$. It's not hard to see that as those elements range over all $i$ and $j$, the entire unipotent radicals are covered. After choosing an appropriate ordering $(\iota_k)$, we see that the image of $\psi_t$ is $\overline{N_t^+}T_t^+N_t^+$, which is open and dense in $Sp(V_t^+)$. We have now established parts (1) and (2) of lemma \ref{deformationlemma} for $\psi_t$.

\textbf{Constructing $\Psi_t$:} We now turn to $N/Z(N),$ which we may identify with $\Hom(W',U),$ assuming $a<b.$ We may extend any linear map $W'\rightarrow U$ to $W\oplus W'$ by setting it to be zero on $W$. Thus we have
$$N/Z(N)\cong \Hom(W',U)\subset \Hom(W\oplus W',U).$$
We identify $Sp(V_t^-\oplus U)/(Sp(V_t^-)\times Sp(U))$ with the symplectic Grassmannian of $2a$-dimensional nondegenerate subspaces of $V_t^-\oplus U$. The set of linear maps $V_t^-\rightarrow U$ with nondegenerate graph is identified, via their graphs, with a dense subset of $Sp(V_t^-\oplus U)/(Sp(V_t^-)\times Sp(U))$. A linear map $V_t^-\rightarrow U$ can be extended to $W\oplus W'=V_t^+\oplus V_t^-$ as zero on $V_t^+.$ Thus we are identifying both $N/Z(N)$ and a dense subset of $Sp(V_t^-\oplus U)/(Sp(V_t^-)\times Sp(U))$ with subsets of $\Hom(W\oplus W',U).$

For $\alpha \in \nQ_p$, let $\alpha E_{i,j}\in \Hom(W',U)$ be the linear map that sends each $f_k$ to zero except $f_i$, which is sent to $\alpha e_j$, where $i\leq 2a$ and $j> 2a.$ Let $\alpha F_{i,j}$ be similarly defined: sending $f_i$ to $\alpha f_j$. Any element of $N/Z(N)$ is a sum of such maps. Define $\Psi_t(\alpha E_{2i,j})\in \Hom(V_t^-,U)$ to send all $v_{k}^-(t)$ (and all $v_k^+(t)$) to zero except $v_{2i-1}^-(t)$, which is mapped to $-\alpha t^2 e_{j}$. Then all the $e_k$ and $f_k$ are fixed by this map except
\begin{align*}\Psi_t(\alpha E_{2i,j})e_{2i-1}&=-t^2\alpha e_{j},\\
\Psi_t(\alpha E_{2i,j})f_{2i}&=\frac{-1}{t^2}\Psi_t(\alpha E_{2i,j})v_{2i-1}^-(t)=\alpha e_{j}.\end{align*}
Note that the graph of this map is 
$$\spn_{k\not=2i-1,k\leq 2a}\{v_i^-(t)\}\cup \{v_{2i-1}^-(t)-\alpha t^2e_j\},$$
which is nondegenerate. We analagously define $$\Psi_t(\alpha E_{2i-1,j}),\Psi_t(\alpha F_{i,j})\in\Hom(V_i^-(t),U)\subset \Hom(W\oplus W',U).$$ We see that as $t\rightarrow 0$, 
\begin{equation}\label{eqn513}\Psi_t(\alpha E_{i,j})\rightarrow \alpha E_{i,j}\text{ and }\Psi_t(\alpha F_{i,j})\rightarrow \alpha F_{i,j}\end{equation} in $\Hom(W\oplus W',U)$. After fixing an ordering of the $i,j$ indices, we once again extend this to a map defined on all of $N/Z(N)$ whose image is all of $\Hom(V_t^-,U)$. As we mentioned, the set of elements of $\Hom(V_t^-,U)$ with nondegenerate graphs is a dense subset of $Sp(V_t^-\oplus U)/(Sp(V_t^-)\times Sp(U))$, so by restricting the domain of $\Psi_t$, we obtain a continuous bijection between dense subsets of $N/Z(N)$ and $Sp(V_t^-\oplus U)/(Sp(V_t^-)\times Sp(U))$. Moreover, for a fixed $n\in N/Z(N)$, it is easy to check that $\Psi_t(n)$ has nondegenerate graph for all but finitely many $t$.

Thus the convergence $(\ref{eqn513})$ in $\Hom(W\oplus W',U)$ establishes part (3) of lemma \ref{deformationlemma}. We have now established parts (1), (2), and (3) of lemma \ref{deformationlemma}.

\textbf{Invariant measures:} 
Fix a Haar measure $dp$ on $(Sp(W)\times Sp(U))N$, and choose Haar measures $dh$, $dz$, and $d\mu_0$ on $Sp(W)\times Sp(U)$, $Z(N)$, and $N/Z(N)$ respectively so that
$$\int_{H_0}f(p)dp=\int_{Sp(W)\times Sp(U)}\int_{N/Z(N)}\int_{Z(N)}f(hnz)dz\, d\mu_0(n)\,dh.$$
The continuous bijection of dense open subsets 
$$\psi_t:Z(N)\rightarrow (Sp(V_t^+)\times Sp(V_t^-))/Sp(W)\leftrightarrow Sp(V_t^+)$$ 
gives rise to a continuous bijection between tangent spaces
$$T_{z}Z(N)\leftrightarrow T_{\psi_t(z)}(Sp(V_t^+)\times Sp(V_t^-))/Sp(W).$$
In particular, since $\psi_t$ maps the identity in $Z(N)$ to the identity coset in $(Sp(V_t^+)\times Sp(V_t^-))/(Sp(W))$, we have a bijection between dense open subsets
$$\Lie Z(N)\leftrightarrow T_{e}(Sp(V_t^+)\times Sp(V_t^-))/Sp(W)=\Lie Sp(V_t^+).$$
This produces a bijection between the dual spaces and from there a bijection between the sets of volume forms
$$\Lambda^{\dim Z(N)}\Lie Z(N)^\ast\leftrightarrow \Lambda^{\dim Z(N)}T_e^\ast (Sp(V_t^+)\times Sp(V_t^-))/Sp(W).$$
The measure $dz$ is induced by such an element $\omega\in \Lambda^{\dim Z(N)}\Lie Z(N)^\ast$ by left translation. From this bijection, we obtain for every $t\not=0$, a volume form $\omega_t$ on $(Sp(V_t^+)\times Sp(V_t^-))/Sp(W)$. Such a volume form induces a left-invariant measure $dz_t$ by left translation. This volume form has the following property. We may embed 
$$Z(N)=(Sp(W)Z(N))/Sp(W)\hookrightarrow Sp(V)/Sp(W).$$
Since we have $\phi_t(z)\rightarrow z$ in $Sp(V)/Sp(W)$ as $t\rightarrow 0$, we have $\omega_t\rightarrow \omega$ in $\Lambda^{\dim N}T_e^\ast Sp(V)/Sp(W)$ as $t\rightarrow 0.$ As a consequence, by left translating we have 
$$L^\ast_{\psi_t(z)}\omega_{z_t}\rightarrow L^\ast_{z}\omega_z$$
in $\Lambda^{\dim N}	T^\ast Sp(V)/Sp(W)$ as $t\rightarrow 0$. For any compact subset $A$ of $Z(N)$, this means that $dz_t(\psi_t(A))\rightarrow dz(A)$ as $t\rightarrow 0.$ 

Similarly, the Haar measure $d\mu_0$ on $N/Z(N)$ induces via
$$\Lambda^{\dim N/Z(N)}\Lie (N/Z(N))^\ast\xleftrightarrow{\Phi_t} \Lambda^{\dim N/Z(N)}T_e^\ast (H_t/(Sp(V_t^+)\times Sp(V_t^-)\times Sp(U)))$$
a left-invariant measure $d\mu_t$ for each $t\not=0$ with the property that $d\mu_t(\Phi_t(A))\rightarrow d\mu_0(A)$ for any compact subset $A\subset N/Z(N).$ Note that $d\mu_t$ and $dz_t$ are not the pushforward measures of $d\mu_0$ and $dz$ but are instead constructed via the bijection on volume forms; thus their left-invariance property is ensured.

The measures $d\mu_t$, $dz_t$, and $dh$ induce a left-invariant measure $dh_t$ on $H_t$ via 
\begin{align*}dh_t(K)=\int_{H_t/(Sp(V_t^+)\times Sp(V_t^-)\times Sp(U))}&\int_{(Sp(V_t^+)\times Sp(V_t^-))/Sp(W)}\\&\int_{Sp(W)\times Sp(U)}\chi_K(nzh)dh\,dz_t\,d\mu_t(n),\end{align*}
where for a compact subset $K\subset H_t$, $\chi_K$ is the characteristic function of $K$. If $dz_t\circ \psi_t$ and $d\mu_t \circ \Psi_t$ denote the pushforward measures of $dz_t$ and $d\mu_t$ to $Z(N)$ and $N/Z(N)$ respectively, then we have
$$\int_A 1dz_t\circ\psi_t\rightarrow\int_A1 dz\text{ and }\int_B1 d\mu_t\circ \Psi_t\rightarrow \int_B 1d\mu_0$$
for any measurable subsets $A\subset Z(N)$ and $B\subset N/Z(N)$ as $t\rightarrow 0.$ It follows that for any functions $f\in C_c(Z(N))$ and $g\in C_c(N/Z(N))$ that
$$\int_{A}f(z)dz_t\circ \psi_t(z)\rightarrow \int_A f(z)dz\text{ and }\int_B g(n)d\mu_t\circ\Psi_t(n)\rightarrow \int_B g(n)d\mu_0(n)$$
as $i\rightarrow 0.$ Since $dz_t\circ \psi_t$ and $d\mu_t\circ \Psi_t$ are pushforward measures, this implies that if $F\in C_c(G)$ (recall $G=Sp(V)$) then
\begin{align*}
&\lim_{t\rightarrow 0}\int_{H_t}F(h)dh_t=\lim_{t\rightarrow 0}\int_{H_t/(Sp(V_t^+)\times Sp(V_t^-)\times Sp(U))}\int_{(Sp(V_t^+)\times Sp(V_t^-))/Sp(W)}\\
&\hspace{7cm}\int_{Sp(W)\times Sp(U)}F(nzh)dh\,dz_t(z)\,d\mu_t(n)\\
&=\lim_{t\rightarrow 0}\int_{N/Z(N)}\int_{Z(N)}\int_{Sp(W)\times Sp(U)}F(\Psi_t(n)\psi_t(z)h)dh\,dz_t\circ \psi_t(z)\,d\mu_t\circ \Psi_t(n)\\
&=\int_{N/Z(N)}\int_{Z(N)}\int_{Sp(W)\times Sp(U)}F(nzh)dz\,d\mu_0(n)\,dh\\
&=\int_{H_0}F(p)dp.
\end{align*}
We have just shown (\ref{eqn512}). This concludes the proof of lemma 5.1.3. 

\textbf{Approximation of matrix coefficients:} We are now ready to complete the proof of lemma \ref{inductionsphericallemma}. Recall that we wish to show
$$\Ind_{H_1}^G \rho\succ \Ind_{H_0}^G \rho|_{Sp(W)\times Sp(U)}\cdot\textbf{1}_N.$$
Recall following from lemma E.1.3 of \cite{Bekka2008} that the set of functions 
$$\xi_{f,v}(x)=\int_{N/Z(N)}\int_{Z(N)}\int_{Sp(W)\times Sp(U)}f(xhnz)\rho(h)vdz\,d\mu_0(n)\,dh$$
with unit norm for $f\in C_c(G)$ and $v\in \cV_\rho$ is total in the space of the representation $\Ind_{H_0}^G \rho|_{Sp(W)\times Sp(U)}\cdot\textbf{1}_N$ and that to show weak containment, we need only approximate matrix coefficients for vectors on this total subset (\cite{Bekka2008}, {E.1.3 and F.1.3}). We have
\begin{align}\label{eqn515}
&\la g^{-1}\cdot \xi_{f,v},\xi_{f,v}\ra\\
&\quad=\int_{Sp(V)/H_0}\int_{H_0}\int_{H_0}f(gxp)\overline{f(xq)}\la\rho(q^{-1}p)v,v\ra dp\, dq\,dx\nonumber\\
&\quad=\int_{Sp(V)}\int_{N/Z(N)}\int_{Z(N)}\int_{Sp(W)\times Sp(U)}f(gxnzh)\overline{f(x)}\la \rho(h)v,v\ra dh\,dz\,d\mu_0(n)dx.\nonumber
\end{align}
for $g\in Sp(V)$. We will handle this expression one integral at a time.

\textbf{Claim 1:} For each $y\in G$ and $z\in Z(N)$,
\begin{align*}\lim_{t\rightarrow 0}&\int_{Sp(W)\times Sp(U)}f(y\psi_t(z)h)\la\rho_t(\psi_t(z)h)v,v\ra dh\\
&=\int_{Sp(W)\times Sp(U)}f(yzh)\la\rho(h)v,v\ra dh.\end{align*}
We may consider $\psi_t(z)$ an element of $Sp(V_t^+)$. Recall that we defined $\rho_t(\psi_t(z)h)=\rho(g_t^{-1}\psi_t(z)hg_t)$, where $g_t$ acts on $W$ by scaling by $\frac{1}{t}$, on $W'$ by scaling by $t$, and $g_t$ fixes $U$. Then $g_t$ commutes with $h\in Sp(W)\times Sp(U).$  We will show $\lim_{t\rightarrow 0} g_t^{-1}\psi_t(z)g_t$ in $Sp(V)$ is the identity element.

Note that $\psi_t(z)$ always acts trivially on elements of $U$, so we need only consider elements of $W$ and $W'$. Observe that $g_t(v_i^\pm(1))=\frac{1}{t}(v_i^\pm(t))$ for each $i$. From the construction of $\psi_t(z)$, we see that for any $w\in W$, $\psi_t(z)w$ is of the form 
$$w+t^3g_tw_z(t),$$ 
where $w_z(t):\nQ_p\rightarrow V_1^+$ is continuous at $t=0$. Then
$$g_t^{-1}\psi_t(z)g_t w=\frac{1}{t}g_t^{-1}\psi_t(z)w=\frac{1}{t}g_t^{-1}(w+t^3g_tw_z(t))=w+t^2w_z(t)\rightarrow w$$
as $t\rightarrow 0$. On the other hand, for $w'\in W'$, $\psi_t(z)w'$ is of the form 
$$w(z)+w'+t^2 w_z'(t),$$ 
where $w_z'(t):\nQ_p\rightarrow W'$ is continuous at $t=0$ and $w(z)\in W$ does not depend on $t$. Then because $g_tw'=tw'$, we have
$$g_t^{-1}\psi_t(z)g_t w'=tg_t^{-1}\psi_t(z)w'=tg_t^{-1}(w(z)+w'+t^2w_z'(t))=t^2w(z)+w'+t^2w_z'(t).$$
This converges to $w'$ as $t\rightarrow 0.$ Hence $\lim_{t\rightarrow 0} g_t^{-1}\psi_t(z)g_t$ in $Sp(V)$ is the identity element.

Because $\rho$ is continuous, we conclude 
$$\la\rho(g_t^{-1}\psi_t(z)g_th)v,v\ra\rightarrow \la  \rho(h)v,v\ra$$
as $t\rightarrow 0.$ Since $f\in C_c(G)$, by Lebesgue's dominated convergence theorem we have
\begin{align*}\lim_{t\rightarrow 0}&\int_{Sp(W)\times Sp(U)}f(y\psi_t(z)h)\la\rho_t(\psi_t(z)h)v,v\ra dh\\
&=\lim_{t\rightarrow 0}\int_{Sp(W)\times Sp(U)}f(y\psi_t(z)h)\la\rho(g_t^{-1}\psi_t(z)g_th)v,v\ra dh\\
&=\int_{Sp(W)\times Sp(U)}f(yzh)\la\rho(h)v,v\ra dh.\end{align*}
This proves claim 1

\textbf{Claim 2:} For each $y\in G$ and $n\in N/Z(N)$,
\begin{align*}\lim_{t\rightarrow 0}&\int_{(Sp(V_t^+)\times Sp(V_t^-))/Sp(W)}\int_{Sp(W)\times Sp(U)}f(y\Psi_t(n)zh)\la\rho_t(\Psi_t(n)zh)v,v\ra dh\,dz\\
&=\int_{Z(N)}\int_{Sp(W)\times Sp(U)}f(ynzh)\la\rho(h)v,v\ra dh dz.\end{align*}
Recall $\Psi_t(n)\in H_t/(Sp(V_t^+)\times Sp(V_t^-)\times Sp(U))$, and so the first integral above is well defined. Using the pushforward measure $dz_t\circ \psi_t$, we have 
\begin{align*}\lim_{t\rightarrow 0}&\int_{(Sp(V_t^+)\times Sp(V_t^-))/Sp(W)}\int_{Sp(W)\times Sp(U)}f(ynzh)\la\rho_t(\Psi_t(n)zh)v,v\ra dh\,dz\\
&=\lim_{t\rightarrow 0}\int_{Z(N)}\int_{Sp(W)\times Sp(U)}f(y\Psi_t(n)\psi_t(z)h)\la\rho_t(\Psi_t(n)\psi_t(z)h)v,v\ra dh\,dz_t\circ \psi_t(z)
\end{align*}

Let us consider 
$$\rho_t(\Psi_t(n)\psi_t(n)h)=\rho(g_t^{-1}\Psi_t(n)g_t\cdot g_t^{-1}\psi_t(n)g_t h).$$
Of course, this is only well-defined in the integral. We need only describe the behavior of $g_t^{-1}\Psi_t(n)g_t$ modulo $Sp(V_t^+)\times Sp(V_t^-)\times Sp(U)$ as $t\rightarrow 0$. Recall that we have described this already as a linear map $V_t^-\rightarrow U$ which by lemma \ref{deformationlemma} (3) converges in $\Hom(W\oplus W',U)$ to $n\in N/Z(N)=\Hom(W',U)$ (recall $\Psi_t(n)$ and $n$ extend to $W'\oplus W$ by acting trivially on $V_t^-$ and $W$ respectively). 

Recall from the construction of $\Psi_t(n)$ that it maps elements of $W$ to elements of the form $t^2u_n$, where $u_n\in U$ does not depend on $t$. Then
$$g_t^{-1}\Psi_t(n)g_t(w)=\frac{1}{t}g_t^{-1}\Psi_t(n)=\frac{1}{t}g_t^{-1}(t^2u_n)=tu_n\rightarrow 0$$
as $t\rightarrow 0$. On the other hand, $\Psi_t(n)$ maps elements of $W'$ to elements of the form $u_n'$, where $u_n'\in U$ does not depend on $U$. Then
$$g_t^{-1}\Psi_t(n)g_t(w')=tg_t^{-1}\Psi_t(n)(w')=tg_t^{-1}u_n'=tu_n'\rightarrow 0$$
as $t\rightarrow 0$. Thus $g_t^{-1}\Psi_t(n)g_t$ approaches the zero map as $t\rightarrow 0$. As elements of $Sp(V_1^-\oplus U)/(Sp(V_1^-)\times Sp(U))$, this means $g_t^{-1}\Psi_t(n)g_t^{-1}$ approaches the identity coset as $t\rightarrow 0.$ By the continuity of $\rho$ and since $f\in C_c(G)$, we conclude by Lebesgue's dominated convergence theorem, claim 1, and the proof of part (4) of lemma \ref{deformationlemma} that
\begin{align*}\lim_{t\rightarrow 0}&\int_{(Sp(V_t^+)\times Sp(V_t^-))/Sp(W)}\int_{Sp(W)\times Sp(U)}f(ynzh)\la\rho_t(\Psi_t(n)zh)v,v\ra dh\,dz\\
&=\lim_{t\rightarrow 0}\int_{Z(N)}\int_{Sp(W)\times Sp(U)}f(y\Psi_t(n)\psi_t(z)h)\\&\hspace{4.5cm}\la\rho(g_t^{-1}\Psi_t(n)g_t\cdot g_t^{-1}\psi_t(z)g_th)v,v\ra dh\,dz_t\circ \psi_t(z)\\
&=\int_{Z(N)}\int_{Sp(W)\times Sp(U)}f(ynzh)\la\rho(h)v,v\ra dh\,dz.
\end{align*}
This proves claim 2. 

The final two integrals are easier. By lemma \ref{deformationlemma} (4), and Lebesgue's dominated convergence theorem, since $f\in C_c(G)$ we can now say for each $g,x\in G$ that
\begin{align*}
&\lim_{t\rightarrow 0}\overline{f(x)}\int_{H_t/(Sp(V_t^+)\times Sp(V_t^-)\times Sp(U))}\int_{(Sp(V_t^+)\times Sp(V_t^-))/Sp(W)}\int_{Sp(W)\times Sp(U)}\\&\hspace{6cm}f(gxnzh)\la\rho_t(nzh)v,v\ra dh\,dz_t(z)\,d\mu_t(n)\\
&=\lim_{t\rightarrow 0}\overline{f(x)}\int_{N/Z(N)}\int_{Z(N)}\int_{Sp(W)\times Sp(U)}\\&\hspace{2cm}f\big(gx\Psi_t(n)\psi_t(z)h\big)\la\rho_t\big(\Psi_t(n)\psi_t(z)h\big)v,v\ra dh\,dz_t\circ\psi_t(z)\,d\mu_t\circ\Psi_t(n)\\\
&=\overline{f(x)}\int_{N/Z(N)}\int_{Z(N)}\int_{Sp(W)\times Sp(U)}f(gxnzh)\la \rho(h)v,v\ra dh\,dz\,d\mu_0(n).
\end{align*}
Finally since this is compactly supported as a function of $x$, we have for each $g\in G=Sp(V)$, by Lebesgue's dominated convergence theorem
\begin{align*}
&\lim_{t\rightarrow 0}\int_G\int_{H_t}f(gxh)\overline{f(x)}\la\rho_t(h)v,v\ra dh_t(h)\,dx\\
&=\lim_{t\rightarrow 0}\int_{G}\int_{H_t/(Sp(V_t^+)\times Sp(V_t^-)\times Sp(U))}\int_{(Sp(V_t^+)\times Sp(V_t^-))/Sp(W)}\int_{Sp(W)\times Sp(U)}\\&\hspace{5cm}f(gxnzh)\overline{f(x)}\la\rho_t(nzh)v,v\ra dh\,dz_t(z)\,d\mu_t(n)\,dx\\
&=\int_G\int_{N/Z(N)}\int_{Z(N)}\int_{Sp(W)\times Sp(U)}f(gxnh)\overline{f(x)}\la \rho(h)v,v\ra dh\,dz\,d\mu_0(n)\,dx.
\end{align*}

The first line we recognize as a matrix coefficient for $\Ind_{H_t}^G\rho_t\cong \Ind_{H_1}^G\rho_1$, and the last line is the matrix coefficient (\ref{eqn515}) for $\Ind_{H_0}^G\rho|_{Sp(W)\times Sp(U)}\cdot\textbf{1}_N.$

\textbf{Uniform convergence on compacta:} We have nearly shown the weak containment (\ref{eqn510}). The convergence of matrix coefficients we have established is only pointwise for $g\in G$. Let $Q\subset G$ be a compact subset. Let $K$ be the support of $f$. Note that $gx\in K$ when $x\in Q^{-1}K$. Let $L=Q^{-1}K\cup Q$. It suffices to show uniform convergence on $L$. 

For each $t\in \nQ_p$, let
$$M_t(g)=\begin{cases}\int_{G}\int_{H_0}f(gxp)\overline{f(x)}\la(\rho|_{Sp(W)\times Sp(U)}\cdot\textbf{1}_N)(p)v,v\ra dp\,dx&\text{if }t=0\\
\int_{G}\int_{H_t}f(gxh)\overline{f(x)}\la\rho(h)v,v\ra dh_t(h)\,dx&\text{if }t\not=0\end{cases}.$$
The pointwise convergence $M_t(g)\rightarrow M_0(g)$ as $t\rightarrow 0$ is what we have just established. Define
$$F_t(x)=\begin{cases}\int_{H_0}\overline{f(xp^{-1})}\la(\rho|_{Sp(W)\times Sp(U)}\cdot\textbf{1}_N)(p)v,v\ra dp\,dx&\text{if }t=0\\
\int_{H_1}\overline{f(xh^{-1})}\la\rho(h)v,v\ra dh_t(h)\,dx&\text{if }t\not=0\end{cases}.$$
Recall that we have already shown $F_t\rightarrow F_0$ pointwise as $t\rightarrow 0.$ Then by change of variables, we have
$$M_t(g)=\int_G f(gx)F_t(x)dx=\int_L f(gx)F_t(x)dx.$$
Hence
\begin{align*}
\sup_{g\in L}|M_t(g)-M_0(g)|\leq (\sup_{g\in G} |f(g)|)\left(\int_L |F_t(x)-F_0(x)|dx\right).
\end{align*}
Observe that $(t,x)\mapsto |F_t(x)-F_0(x)|$ is continuous as a function on $\nQ_p\times L$, and in particular this means it is bounded on the compact set
$$\{t\in \nQ_p:|t|<1\}\times L.$$
Since $|F_t(x)-F_0(x)|\rightarrow 0$ pointwise as $t\rightarrow 0$ and $L$ is compact, by Lebesgue's dominated convergence theorem, $\int_L |F_t(x)-F_0(x)|dx\rightarrow 0$. Thus $M_t(x)\rightarrow M_0(g)$ uniformly on $L$ and hence on $Q\subset L$. We have therefore shown (\ref{eqn510}):
$$\Ind_{H_1}^G\rho\succ\Ind_{H_0}^G\rho|_{Sp(W)\times Sp(U)}\cdot \textbf{1}_N.$$

\textbf{The $SL_2$ type:} At long last, let us see that this representation has the type predicted by lemma \ref{inductionsphericallemma}. Recall that $\rho=\rho_a\boxtimes\rho_b$ is a representation of 
$$H_1=Sp(V_1^+)\times Sp(V_1^-\oplus U)=Sp_{2a}\times Sp_{2b},$$ 
and $Sp(W)=Sp_{2a}$ is diagonally embedded in this product, while $Sp(U)=Sp_{2b-2a}$ is contained in the $Sp_{2b}$ factor of $H_1.$ If $\la\lambda_i\ra$ is the type of $\rho_b$, then by the induction hypothesis for restriction on $Sp_{2b}$, $\rho_b|_{Sp(V_1^-)\times Sp(U)}$ weakly contains a representation of the form $\gamma_a\boxtimes\gamma_b$, where $\gamma_a$ is an unramified Arthur type representation of $Sp(V_1^-)$, and $\gamma_b$ is an unramified Arthur type representation of $Sp(U)$ of type $\la \lambda_i-2a\ra.$ Then by inducing in stages,
$$\Ind_{H_1}^G\rho_a\boxtimes\rho_b\succ\Ind_{H_0}^G (\rho_a\otimes \gamma_a)\boxtimes\gamma_b\cdot\textbf{1}_N=\Ind_{P}^G\left(\Ind_{Sp(W)}^{GL(W)}\rho_a\otimes\gamma_a\right)\boxtimes\gamma_b\cdot\textbf{1}_N.$$
This is parabolically induced, and so this is of the form (\ref{stdrep}), and we may read the type of this representation. If $\rho_a\otimes\gamma_a$ is not one-dimensional, then $\Ind_{Sp(W)}^{GL(W)}\rho_a\otimes\gamma_a$ is of type $\la 1,\dots,1\ra$. The type of $\pi$ is then $\la \lambda_i-2a\ra$ as expected. 

However, if $\rho_a\otimes\gamma_a$ is one-dimensional (which is only possible if $\rho_a$ and $\rho_b$ are one-dimensional), then $\Ind_{Sp(W)}^{GL(W)}\rho_a\otimes\gamma_a$ has type $\la 2,\dots,2\ra$, and consequently $\pi$ has type $\la 2b-2a+1,2,\dots,2\ra$ since the type of $\rho_b$ in that case is $\la 2b+1\ra.$ This completes the proof of lemma \ref{inductionsphericallemma}.
\end{proof}
\begin{remark}
Since we are working in the $p$-adic case, the map 
$$t\mapsto \begin{cases}\la g^{-1}\cdot \xi_{f,v},\xi_{f,v}\ra&t=0\\\int_{Sp(V)}\int_{H_t}f(gxh)\overline{f(x)}\la \rho_t(h)v,v\ra dh_t(h)dx&t\not=0\end{cases}$$
for locally constant $f$ is a locally constant function $\nQ_p\rightarrow \nC.$ That means for $t$ sufficiently small, we actually get equality of matrix coefficients
$$\int_{Sp(V)}\int_{H_t}f(gxh)\overline{f(x)}\la \rho_t(h)v,v\ra dh_t(h)dx=\la g^{-1}\cdot \xi_{f,v},\xi_{f,v}\ra$$
for all $g$ in a given compact subset.
\end{remark}
\begin{remark}
In the case where $\rho$ is trivial, this result coincides with lemma 11.1.1 of \cite{Sakellaridis2012}.
\end{remark}
\begin{remark}
The above argument is an explicit model of what should be a much more general phenomenon. Given a spherical variety $G/H$ and a representation $\rho$ of $H$, we might expect that for an appropriate choice of a parabolic subgroup $P=MN$ (such as the standard parabolic associated to that variety) that there's a subgroup $M_0$ of $M$ such that 
$$\Ind_H^G \rho\succ \Ind_{M_0N}^G \rho|_{M_0}\cdot \textbf{1}_N.$$
This could be computed using the theory of deformations to the normal cone, but we don't consider this general theory here. It may be possible to obtain a full decomposition of $\Ind_H^G \rho$ in terms of such deformations. 
\end{remark}

\subsection{Induction from $Sp_{2n}$ to $GL_{2n}$}\label{inductionexceptionalprop}
In this section, we give a proof of the following sub-case of proposition \ref{weakprop}. 
\begin{lemma}\label{inductionexceptionallemma}
Let $\rho$ be an unramified Arthur type representation of $Sp_{2n}$, and let $\pi=\Ind_{Sp_{2n}}^{GL_{2n}}\rho$.
\begin{enumerate}[(1)]
	\item $\pi$ weakly contains a tempered representation (of type $\la 1,\dots,1\ra$) if $\rho$ is not one dimensional and $n>1$.
	\item $\pi$ weakly contains a representation of type $\la 2,\dots,2\ra$ if $\rho$ is one dimensional.
	\item If $n=1$, $\pi=\Ind_{SL_2}^{GL_2}\rho$ has the same type as $\rho$ (as an $SL_2$ representation), and $\pi$ is tempered if $\rho$ is.
\end{enumerate}
\end{lemma}
The proof for this section does not rely on any of the cases of theorem \ref{mainthm}. Thus once we have shown lemma \ref{inductionexceptionallemma}, we may immediately conclude proposition \ref{weakprop} holds by the global argument given after proposition \ref{unramifiedspectrumofarthurtype}.

We will use a group deformation argument similar to the last section and argue by induction on $n$. For the base case $Sp_2$ to $GL_2$, since $Sp_2=SL_2$, we know that inducing $SL_2$ to $GL_2$ does not change the type of the representation (see section 3.5 of \cite{Clozel2007}). Thus a one-dimensional character of $Sp_2$ has the trivial type when induced to $GL_2$, which is $\la 2\ra$, and any other Arthur type representation of $GL_2$ has type $\la 1,1\ra.$ 

Given the structure of the inductive argument of theorem \ref{mainthm} and the way all of the exceptional cases have arisen so far, it is worth taking special note of the following.
\begin{remark}\label{exceptionalremark}
	All of the exceptional cases of theorem \ref{mainthm} are ultimately a consequence of the fact that $Sp_2=SL_2$ and inducing from $SL_2$ to $GL_2$ does not change the $SL_2$-type.
\end{remark}

\begin{proof}[Proof of lemma \ref{inductionexceptionallemma}] To begin our inductive argument, suppose lemma \ref{inductionexceptionallemma} holds for the induction from $Sp_{2m}$ to $GL_{2m}$ when $m<n$. We realize $G=GL_{2n}$ as $GL(V)$ for some $2n$-dimensional vector space $V$, and fix a basis $e_1,\dots,e_n,f_1,\dots,f_n$ of $V$. For each $t\not=0$, define a symplectic form $\la\cdot,\cdot\ra_t$ according to the relations
\begin{align*}
\la e_i,e_j\ra_t=\la f_i,f_j\ra_t&=0, \\
\left\la \frac{1}{t}e_1,\frac{1}{t}f_j\right\ra_t&=\delta_{1,j},\\
\left\la \frac{1}{t}e_i,\frac{1}{t}f_1\right\ra_t&=\delta_{i,1},\\
\left\la \frac{1}{t^2}e_i,\frac{1}{t^2}f_j\right\ra_t&=\delta_{i,j}\text{ for }i,j>1.
\end{align*}
Then $\la e_1,f_1\ra_t=t^2$ and $\la e_i,f_i\ra=t^4$, and so as $t\rightarrow 0$, this form becomes identically zero. Let $H_t$ be the set of all elements of $GL(V)$ which preserve the form $\la\cdot,\cdot\ra_t$ for $t\not=0.$ Let $W=\spn\{e_1,f_1\}$ and $W'=\spn\{ e_2,\dots,e_n,f_2,\dots,f_n\}$. For each $t\not=0$, $W'$ is the $\la\cdot,\cdot\ra_t$-complement of $W$.

Let $P$ be the parabolic subgroup of $GL(V)$ stabilizing $W'$. We may decompose $P=MN$ where $M=GL(W)\times GL(W')$ and $N\cong \Hom(W,W').$ The forms $\la\cdot,\cdot\ra_t$ when restricted to $W$ or to $W'$ are all multiples of each other and as such have the same stabilizer. That is, $\{g\in GL(W):\la gv,gw\ra_t=\la v,w\ra_t\forall v,w\in W\}$ does not change as $t$ varies since for $s\not=t$, 
$$\la gv,gw\ra_t=\la v,w\ra_t\text{ implies } \la gv,gw\ra_s=\frac{s^2}{t^2}\la gv,gw\ra_t=\frac{s^2}{t^2}\la v,w\ra_t=\la v,w\ra_s,$$
and this is similarly true for $W'$. As such, there exists a subgroup $Sp(W)\times Sp(W')$ contained in the intersection of all $H_t$. In fact, this is equal to the intersection of all $H_t$ since the scaling factors on the restrictions of the forms to $W$ and $W'$ are different.

Define $H_0=(Sp(W)\times Sp(W'))N.$ We wish to show for a representation $\rho$ of $H_1$ that
\begin{equation}\label{eqn521}\Ind_{H_1}^{GL(V)}\rho\succ \Ind_{H_0}^{GL(V)}\rho|_{Sp(W)\times Sp(W')}\cdot \textbf{1}_N.\end{equation}
Let $g_t$ be the map $e_1\mapsto \frac{1}{t} e_1,$ $f_1\mapsto \frac{1}{t}f_1$, $e_i\mapsto \frac{1}{t^2}e_i$, and $f_i\mapsto \frac{1}{t^2}f_i$ for $i>1.$ Then $\la g_tv,g_tw\ra_t=\la v,w\ra_1$ for every $t\not=0$ and $v,w\in V.$ Then for $h\in H_1$, 
$$\la g_thg_t^{-1}v,g_thg_t^{-1}w\ra_t=\la hg_t^{-1}v,hg_t^{-1}w\ra_1=\la g_t^{-1}v,g_t^{-1}w\ra_1=\la v,w\ra_t,$$
and so $g_tH_1g_t^{-1}=H_t.$ We may then define the representation $\rho_t(h)=\rho_1(g_t^{-1}hg_t)$ of $H_t$. 

In order to show (\ref{eqn521}), we will prove the following (cf. lemma \ref{deformationlemma})
\begin{lemma}\label{deformationlemma2}
For every $t\not=0$, there exists a continuous bijection $\phi_t$ from an open dense subset of $N$ to an open dense subset of $H_t/(Sp(W)\times Sp(W'))$ with the following properties.
\begin{enumerate}[(1)]
	\item For every $n\in N$, $n$ is in the domain of $\phi_t$ for all but finitely many $t\in \nQ_p^\times.$
	\item We have natural identification $N\cong \Hom(W,W')$, and we identify a dense open subset of $\Hom(W,W')$ with $H_t/(Sp(W)\times Sp(W')).$
	In $\Hom(W,W')$, $\lim_{t\rightarrow 0} \phi_t(n)=n$ for all $n\in N$.
	\item Given a Haar measure $dp$ on $H_0$, there exist Haar measures $dh$ and $dn$ on $Sp(W)\times Sp(W')$ and $N$ respectively such that
	$$\int_{H_0} f(p)dp=\int_{N}\int_{Sp(W)\times Sp(W')}f(nh)dh\,dn$$
	for all $f\in C_c(G)$. There then exists for each $t\not=0$ an invariant measure $dn_t$ on $H_t/(Sp(W)\times Sp(W'))$ which defines a Haar measure $dh_t$ on $H_t$ in the following way. Let $K\subset H_t$ be compact and let $\chi_K$ be the characteristic function of $K$. Then we define
	$$dh_t(K)=\int_{H_t/(Sp(W)\times Sp(W'))}\int_{Sp(W)\times Sp(W')}\chi_K(nh)dh\,dn_t(n).$$
	This $dh_t$ has the property that
	$$\lim_{t\rightarrow 0}\int_{H_t}f(x)dh_t(x)=\int_{H_0}f(p)dp$$
	for any $f\in C_c(G).$
\end{enumerate}
\end{lemma}

\textbf{Constructing $\phi_t$:} We begin the proof of lemma \ref{deformationlemma2} by constructing $\phi_t$. Observe that $N=\Hom(W,W')$ and we may identify $H_t/(Sp(W)\times Sp(W'))$ with the Grassmannian of $2$-dimensional nondegenerate spaces, almost all of which can be realized as the graph of a homomorphism $W\rightarrow W'.$ The set of $2$-dimensional nondegenerate spaces is open and dense in the Grassmannian of $2$-dimensional subspaces of $V$, and so the set of homomorphisms $W\rightarrow W'$ with nondegenerate graph is open and dense in $N$. Then we immediately have a bijection between open dense subsets of $N$ and $H_t/(Sp(W)\times Sp(W')).$ We must show that this identification produces the required properties.

Consider the homomorphism $n_{e_1,ae_i}$ which maps $e_1$ to $ae_i$, $i>1$, $a\in \nQ_p$, and sends $f_1$ to $0$. This can be written as $\phi_t(n_{e_1,ae_i})\in H_t/(Sp(W)\times Sp(W'))$ which maps $\frac{1}{t}e_i\mapsto at\left(\frac{1}{t^2}e_i\right)$; this is the expression of the map in the symplectic basis for $\la\cdot,\cdot\ra_t$. This map has nondegenerate graph, as a basis for the graph is $\{e_1+ae_i,f_1\}.$ We have $g_t^{-1}\phi_t(n_{e_1,ae_i})g_te_i= ate_i$, which as a sequence of elements of $H_1/(Sp(W)\times Sp(W'))$ approaches the identity element. Thus by continuity of $\rho$, $\la \rho_t(\phi_t(n_{e_1,ae_i})\eta,\eta\ra\rightarrow \la \eta,\eta\ra$ as $t\rightarrow 0.$

We now argue that $\phi_t(n_{e_1,ae_i})\rightarrow n_{e_1,ae_i}$ as $t\rightarrow 0.$ Observe that the symplectic complement of  $\graph \phi_t(ne_{e_1,ae_i})$ defines the graph of a map $\phi_t(n_{e_1,ae_i})^\perp:W'\rightarrow W$ which sends fixes all $f_j$ and $e_j$ except $f_i$, which is mapped to $-at^2f_1$, since $\la e_1+ae_i,-at^2f_1+f_i\ra_t=0.$ Both of these maps describe the same element of $H_t/(Sp(W)\times Sp(W'))$. From the description of $\phi_t(n_{e_1,ae_i})^\perp$ it's clear that the action of this element on $W'$ approaches the identity as $t\rightarrow 0$, and from the description of $\phi_t(n_{e_1,ae_i})$ it's clear that the action on $W$ approaches the $n_{e_1,a_ei}$ action on $W$. It must be the case then that $\phi_t(n_{e_1,ae_i})\rightarrow n_{e_1,ae_i}$ as $t\rightarrow 0.$ 

We may similarly construct maps $\phi_t(n_{e_1,af_i})$, $\phi_t(n_{f_1,ae_i})$, and $\phi_t(n_{f_1,af_i})$, which then fully describes $\phi_t$, and we have $\phi_t(n)\rightarrow n$ as $t\rightarrow 0$. This shows part (2) of lemma \ref{deformationlemma}.

We now show part (1) of lemma \ref{deformationlemma2}. Suppose $n\in N=\Hom(W,W')$. Recall that for any $t$, the form $\la \cdot,\cdot\ra_t$ only changes by scaling when restricted to $W$ or $W'$. If for all $t\not=0$ the range of $n$ in $W'$ is an isotropic space, then 
$$\la e_1+n(e_1),f_1+n(f_1)\ra_t=\la e_1,f_1\ra_t\not=0.$$
Otherwise, the range of $n$ is non-isotropic for all $t\not=0$. Suppose $x,y\in W$ are such that $\la n(x),n(y)\ra_t\not=0$ for all $t\not=0$ (this again is possible because $\la\cdot,\cdot\ra_t$ on $W'$ just changes by scalar multiple as $t$ varies). Then 
$$\la x+n(x),y+n(y)\ra_t=\la x,y\ra_t+\la n(x),n(y)\ra_t=t^2\la x,y\ra_1+t^4\la n(x),n(y)\ra_1.$$
This can be zero for at most two nonzero values of $t$. Thus $\graph n$ is nondegenerate for all but finitely many $t$, which proves part (1) of lemma \ref{deformationlemma2}.

\textbf{Invariant measures and matrix coefficients:} From here, the proof is very similar to the last section, and so we will omit most of the details. Given $dp$, $dh$, and $dn$ as in lemma \ref{deformationlemma2} (3), we construct $dn_t$ via the correspondence between volume forms
$$\Lambda^{\dim N}\Lie N^\ast\leftrightarrow \Lambda^{\dim N}T_e^\ast H_t/(Sp(W)\times Sp(W')).$$
We then define $dh_t$ as in lemma \ref{deformationlemma2} (3), and by the arguments of the previous section we have
$$\int_{H_t}f(x)dh_t(x)\rightarrow \int_{H_0}f(p)dp$$
as $t\rightarrow 0$ for $f\in C_c(G).$ This completes the proof of lemma \ref{deformationlemma2}.

Continuing with the proof of lemma \ref{inductionexceptionallemma}, we will now approximate matrix coefficients to show (\ref{eqn521}). Since $\phi_t(n)\rightarrow n$ as $t\rightarrow 0$, we obtain for each $y\in G=GL(V)$, $n\in N$, $v\in \cV_\rho$ and $f\in C_c(G)$
$$\int_{Sp(W)\times Sp(W')}f(y\phi_t(n)h)\la\rho_t(\phi_t(n)h)v,v\ra dh\rightarrow \int_{Sp(W)\times Sp(U)}f(ynh)\la \rho(h)v,v\ra dh$$
as $t\rightarrow 0$ by the arguments of the previous section. It follows from Lebesgue's dominated convergence theorem that we get convergence of matrix coefficients
\begin{align*}\int_{G}\int_{H_t}&f(gxh)\overline{f(x)}\la\rho_t(h)v,v\ra dh_t(h)\\&\rightarrow \int_{G}\int_{N}\int_{Sp(W)\times Sp(W')}f(gxnh)\overline{f(x)}\la\rho(h)v,v\ra dhdndx\end{align*}
as $t\rightarrow 0$ for each $g\in G$. After showing uniform convergence on a given compact subset of $G$, we obtain 
$$\pi=\Ind_{H_1}^G \rho\succ\Ind_{H_0}^G \rho|_{Sp(W)\times Sp(W')}\cdot\textbf{1}_N.$$

\textbf{The $SL_2$ type:} Suppose $\rho_a\boxtimes\rho_b$ is an unramified Arthur type representation weakly contained in $\rho|_{Sp(W)\times Sp(W')}$, where $\rho_a$ and $\rho_b$ are unramified Arthur type representations of $Sp(W)=Sp_2=SL_2$ and $Sp(W')=Sp_{2n-2}$ respectively. By inducing in stages,
$$\pi\succ \Ind_{P}^G\left(\Ind_{Sp(W)}^{GL(W)}\rho_a\right)\boxtimes\left(\Ind_{Sp(W')}^{GL(W')}\rho_b\right)\cdot\textbf{1}_N.$$
This is parabolic induction, and so we may read the type of this representation. If $\rho$ is not one-dimensional, then $\rho_a$ must be an $SL_2$ representation of type $\la 1,1\ra$ (in fact, it is tempered, which is easy to show by Mackey theorey as in section \ref{restrictionsphericalinduced}). On the other hand, if $\rho$ is not one-dimensional, then neither is $\rho_b$, and so the induction hypothesis for this section says that $\Ind_{Sp(W')}^{GL(W)}\rho_b$ is tempered. Hence $\pi$ is tempered in this case.

On the other hand, if $\rho$ is one-dimensional, then $\rho_a$ is a one dimensional representation of $SL_2$, so its induction to $GL(W)$ is of type $\la 2\ra$. Likewise, $\rho_b$ is one dimensional, and the induction hypothesis for this section says $\Ind_{Sp(W')}^{GL(W')}\rho_b$ is of type $\la 2,\dots,2\ra$. This yields a type of $\la 2,\dots,2\ra$ for $\pi$. This completes the proof of lemma \ref{inductionexceptionallemma}.\end{proof}
\begin{remark}\label{gurevichremark}
	The proof technique for this section is not really necessary in the case of induction from $O_n$ to $GL_n$ since this result is already known. By corollary 5.10 of \cite{Gurevich2016}, the spherical variety $GL_n/O_n$ is strongly tempered in the sense of Sakellaridis and Venkatesh \cite{Sakellaridis2012}. By theorem 6.2.1 of \cite{Sakellaridis2012} this implies $L^2(GL_n/O_n)=\Ind_{O_n}^{GL_n}\textbf{1}_{O_n}$ is a tempered representaton of $GL_n$. Therefore by lemma 1 of \cite{Venkatesh2005}, $\Ind_{O_n}^{GL_n}\rho$ is tempered for any unitary representation $\rho$.  
\end{remark}
\subsection{Induction from $GL_n$ to $O_{2n}$}\label{inductionGL}
In this section, we will prove the following sub-case of proposition \ref{weakprop}.
\begin{lemma}\label{inductionlevilemma}
	Let $\sigma$ be an unramified Arthur type representation of $GL_n$, where $GL_n$ is embedded in $O_{2n}$ as a Siegel Levi subgroup. Then
	$$\pi=\Ind_{GL_n}^{O_{2n}}\sigma$$
	weakly contains a representation of type $\la 1,\dots,1\ra$, except in the case where $\sigma$ is one-dimensional, in which case $\pi$ weakly contains a representation of type $\la 2,\dots,2\ra$ if $n$ is even or $\la 1,1,2\dots,2\ra$ if $n$ is odd.
\end{lemma}
There are two reasons we prove the $O_{2n}$ case rather than the $Sp_{2n}$ case here. The first is that it is already known that inducing from $GL_{n}$ to $Sp_{2n}$ is tempered (see remark 5.14 of \cite{Gurevich2016} and lemma 1 of \cite{Venkatesh2005}). The second is that the orthogonal case produces an exceptional case for theorem \ref{mainthm}. 

There are two ways we could approach the proof of lemma \ref{inductionlevilemma}. The more difficult way is to use the boundary degeneration technique of section \ref{inductionspherical}, which we briefly sketch. If $\{e_1,\dots,e_n,f_1,\dots f_n\}$ is an orthogonal basis for a $2n$-dimensional orthogonal vector space $V$, then setting $v_i^\pm(t)=e_i\pm tf_i$ and $W^\pm_t=\spn\{v_1^\pm(t),\dots,v_n^\pm(t)\}$, we can realize $GL_n$ as the stabilizer of the isotropic splitting $V=W_t^+\oplus W_t^-$ for each $t\not=0$. 
	
As $t\rightarrow 0$, this approaches a subgroup of the stabilizer of $W=\spn\{e_1,\dots,e_n\}$, which is the Siegel parabolic $P$. The set of elements of $GL(W_t^+)$ which already stabilize $W$ are the stabilizer of a skew-symmetric form on $W$ which is nondegenerate if $\dim W=n$ is even or has one-dimensional kernel if $\dim W=n$ is odd as we will see when we prove lemma \ref{inductionlevilemma} below. Using the arguments of section \ref{inductionspherical}, we could eventually obtain the weak containment (\ref{eqn531}). 

Instead of going through all of that, we will induce through the Siegel parabolic and perform harmonic analysis on the abelian unipotent radical. Let us begin.

\begin{proof}[Proof of lemma \ref{inductionlevilemma}] Let $P$ be a Siegel parabolic subgroup of $G=O_{2n}$ with Levi subgroup $M=GL_n$. We will induce $\sigma$ in stages through $P$. This proof is essentially the same as the argument for induction from a Levi subgroup in \cite{Venkatesh2005} but with appropriate modification for the split orthogonal case. We have
$$\Ind_{M}^{G}\sigma=\Ind_P^{G}(\sigma\cdot \textbf{1}_N)\otimes L^2(P/M).$$

Let us analyze $L^2(P/M)$, which as a $P$-representation we may realize on $L^2(N).$ The group $M$ acts on the abelian group $N$ by conjugation and hence acts on $\hat{N}$ via
$$(m\cdot \chi)(n)=\chi(m^{-1}nm)$$
for $\chi\in \hat{N}.$ The Fourier transform is an isometry $L^2(N)\rightarrow L^2(\hat{N})$ which intertwines the action of $M$ on $L^2(N)$ with an $M$-action on $L^2(\hat{N})$.

After identifying $N$ with the set of skew-symmetric $n\times n$ matrices over $\nQ_p$, we can identify the $M$ orbits on $L^2(\hat{N})$ with the possible dimensions of kernels of skew-symmetric bilinear forms on an $n$-dimensional space $V$. There is one open orbit; it corresponds to the most non-degenerate skew-symmetric forms on $V$. Suppose $\zeta\in \hat{N}$ belongs to the open $M$-orbit, and let $M_\zeta$ be its stabilizer. Then we have
$$\Ind_{M}^{MN}\textbf{1}_M\cong \Ind_{M_\zeta N}^{MN}\textbf{1}_{M_\zeta}\cdot\zeta.$$
If we let $\zeta_t(x)=\zeta(tx)$, then as $t\rightarrow 0$, $\zeta_t$ approaches the trivial character in the Fell topology on $\hat{N}$. The stabilizer $M_\zeta$ does not change under this scaling by $t$. By continuity of induction with respect to the Fell topology, we have 
$$\Ind_{M}^{MN}\textbf{1}_M\succ \Ind_{M_\zeta N}^{MN}\textbf{1}_{M_\zeta N}=L^2(GL_n/M_\zeta)\cdot \textbf{1}_N.$$
We therefore have
\begin{equation}\label{eqn531}\Ind_M^{G}\sigma\succ\Ind_P^{G}\left(\Ind_{M_\zeta}^G\sigma|_{M_\zeta}\right)\cdot \textbf{1}_N.\end{equation}

Suppose that $n$ is even. Then $M_\zeta=Sp_{n}$, and $\Ind_{M_\zeta}^M\sigma|_{M_\zeta}$ is of type $\la 1,\dots,1\ra$ unless $\sigma$ is one-dimensional, in which case the type is $\la 2,\dots,2\ra$. Since we have parabolic induction to $G$, the type of $\Ind_M^G\sigma$ is therefore $\la 1,\dots,1\ra$ unless $\sigma$ is one-dimensional, in which case the type is $\la 2,\dots,2\ra$ in the case when $n$ is even.

On the other hand, suppose $n$ is odd. Then the a skew-symmetric form $(\cdot,\cdot)$ on the $n$-dimensional space $V$ which is as non-degenerate as possible has one-dimensional kernel $K$. Let $P'=M'N'$ be the parabolic subgroup of $GL(V)$ stabilizing $K$. For any $n-1$-dimensional complement $V'$ to $K$, we have a Levi subgroup $M'=GL(K)\times GL(V')$ of $P'$. The form $(\cdot,\cdot)$ is nondegenerate when restricted to $V'$, and so the stabilizer in $M'$ of such a form is $GL(K)\times Sp(V')$. The unipotent radical $N'$ acts trivially on $K$ and $V/K$, and so it will stabilize $(\cdot,\cdot)$ as well. From this description, we see
$$M_\zeta=(GL(K)\times Sp(V'))N'.$$
By proposition \ref{restrictiontoparabolic}, $\sigma|_{P'}$ weakly contains a representation of the form $\sigma_K\boxtimes \sigma_{S}\cdot \textbf{1}_{N'}$, where $\sigma_K\boxtimes \sigma_S$ is an irreducible unramified Arthur type representation weakly contained in $\sigma|_{M'}$. Let
$$\pi=\Ind_{P'}^{GL(V)}\sigma_K\boxtimes\left(\Ind_{Sp(V')}^{GL(V')}\sigma_S|_{Sp(V')}\right)\cdot\textbf{1}_{N'}.$$
Then from (\ref{eqn531}), by inducing in stages we have
$$\Ind_M^G\sigma\succ\Ind_P^G\pi\cdot\textbf{1}_N.$$
The type of $\pi$ is $\la 1,\dots,1\ra$ if $\sigma$ is not one-dimensional or $\la 1,2,\dots,2\ra$ if $\sigma$ is one-dimensional. Therefore we have shown that the type of $\Ind_M^G\sigma$ is $\la 1,\dots,1\ra$ unless $\sigma$ is one-dimensional, in which case it is of type $\la 1,1,2,\dots,2\ra$ when $n$ is odd. This proves lemma \ref{inductionlevilemma}.
\end{proof}
We have now given a proof for every case of proposition \ref{weakprop}. This therefore concludes the proofs of proposition \ref{weakprop} and theorem \ref{mainthm}.
\newpage
\addcontentsline{toc}{section}{Bibliography}
\bibliographystyle{plain}
\bibliography{thesisbibliography}

\begin{thebibliography}{10}

\bibitem{Arthur2013}
J.~Arthur.
\newblock {\em {The Endoscopic Classification of Representations: Orthogonal
  and Symplectic Groups}}, volume~61 of {\em Colloquium Publications}.
\newblock American Mathematical Society, 2013.

\bibitem{Bekka2008}
B.~Bekka, P.~de~la Harpe, and A.~Valette.
\newblock {\em Kazhdan{\textquotesingle}s Property (T)}.
\newblock Cambridge University Press, 2008.

\bibitem{Burger1992}
M.~Burger, J.~S. Li, and P.~Sarnak.
\newblock Ramanujan duals and automorphic spectrum.
\newblock {\em Bulletin of the American Mathematical Society}, 26(2):253--258,
  April 1992.

\bibitem{Burger1991}
M.~Burger and P.~Sarnak.
\newblock Ramanujan duals {II}.
\newblock {\em Inventiones Mathematicae}, 106(1):1--11, December 1991.

\bibitem{Casselman1995}
W.~Casselman.
\newblock Introduction to the theory of admissible representations of $p$-adic
  reductive groups.
\newblock Available on Casselman's website:
  \href{https://www.math.ubc.ca/~cass/research/publications.html}{https://www.math.ubc.ca/\url{~}cass/research/publications.html},
  May 1995.

\bibitem{Clozel2004}
L.~Clozel.
\newblock {Combinatorial consequences of Arthur's conjectures and the
  Burger-Sarnak method}.
\newblock {\em International Mathematics Research Notices}, 2004(11):511, 2004.

\bibitem{Clozel2007}
L.~Clozel.
\newblock Spectral theory of automorphic forms.
\newblock In Peter Sarnak and Freydoon Shahidi, editors, {\em Automorphic Forms
  and Applications}, volume~12 of {\em Ias/ Park City Mathematics Series},
  pages 43--93. American Mathematical Society, 2007.

\bibitem{Clozel2003}
L.~Clozel and E.~Ullmo.
\newblock Equidistribution des points de {H}ecke.
\newblock In {\em Contributions to automorphic forms, geometry, and number
  theory}, pages 193--254. Johns Hopkins University Press, 2004.

\bibitem{Cogdell2004}
J.~W. Cogdell, H.~H. Kim, I.~I. Piatetski-Shapiro, and F.~Shahidi.
\newblock Functoriality for the classical groups.
\newblock {\em Publications math{\'{e}}matiques de
  l{\textquotesingle}{IH}{\'{E}}S}, 99(1):163--233, June 2004.

\bibitem{Collingwood2017}
D.~Collingwood and W.~McGovern.
\newblock {\em {Nilpotent Orbits in Semisimple Lie Algebras}}.
\newblock New York: Van Nostrand Reinhold, October 1993.

\bibitem{Gan2014}
W.~T. Gan and R.~Gomez.
\newblock A conjecture of {S}akellaridis{\textendash}{V}enkatesh on the unitary
  spectrum of spherical varieties.
\newblock In {\em Symmetry: Representation Theory and Its Applications}, pages
  185--226. Springer New York, 2014.

\bibitem{Gan2014a}
W.~T. Gan, Y.~Qiu, and S.~Takeda.
\newblock The regularized {S}iegel{\textendash}{W}eil formula (the second term
  identity) and the {R}allis inner product formula.
\newblock {\em Inventiones mathematicae}, 198(3):739--831, March 2014.

\bibitem{Gurevich2016}
M.~Gurevich and O.~Offen.
\newblock A criterion for integrability of matrix coefficients with respect to
  a symmetric space.
\newblock {\em Journal of Functional Analysis}, 270(12):4478--4512, June 2016.

\bibitem{Harris2011}
M.~Harris, J.~S. Li, and B.~Sun.
\newblock Theta correspondence for close unitary groups.
\newblock {\em Advanced Lectures in Mathematics}, 19:265--308, January 2011.

\bibitem{Hongyu2000}
H.~He.
\newblock Theta correspondence {I} {\textemdash} semistable range: Construction
  and irreducibility.
\newblock {\em Communications in Contemporary Mathematics}, 02(02):255--283,
  May 2000.

\bibitem{Howe2010}
R.~Howe.
\newblock On a notion of rank for unitary representations of the classical
  groups.
\newblock In {\em Harmonic Analysis and Group Representation}, pages 224--331.
  Springer Berlin Heidelberg, 2010.

\bibitem{Knapp1986}
A.~Knapp.
\newblock {\em Representation Theory of Semisimple Groups: An Overview Based on
  Examples}.
\newblock Princeton University Press, revised edition, 1986.

\bibitem{Kudla1983}
S.~Kudla.
\newblock Seesaw dual reductive pairs.
\newblock In {\em Automorphic Forms of Several Variables}, pages 244--268,
  1983.

\bibitem{Kudla1994}
S.~Kudla and S.~Rallis.
\newblock {A Regularized Siegel-Weil Formula: The First Term Identity}.
\newblock {\em The Annals of Mathematics}, 140(1):1--80, July 1994.

\bibitem{Lapid2009}
E.~Lapid and J.~Rogawski.
\newblock On a result of {V}enkatesh on {C}lozel's conjecture.
\newblock In {\em Automorphic Forms and L-functions II: Local Aspects
  (Contemporary Mathematics: Israel Mathematical Conference Proceedings)},
  pages 173--178. Amer Mathematical Society, 2009.

\bibitem{Li1989}
J.~S. Li.
\newblock Singular unitary representations of classical groups.
\newblock {\em Inventiones Mathematicae}, 97(2):237--255, June 1989.

\bibitem{Ma2017}
J.~Ma, B.~Sun, and C.~B. Zhu.
\newblock Unipotent representations of real classical groups.
\newblock Preprint,
  \href{https://arxiv.org/abs/1712.05552}{https://arxiv.org/abs/1712.05552},
  2017.

\bibitem{Moeglin2009}
C.~M{\oe}glin.
\newblock {Comparaison des param\`{e}tres de Langlands et des exposants \'{a}
  l'int\'{e}rieur d'un paquet d'Arthur}.
\newblock {\em {Journal of Lie Theory}}, 1(4):797--840, 2009.

\bibitem{Moeglin2011}
C.~M{\oe}glin.
\newblock {Multiplicit\'{e} 1 dans les paquets d'Arthur aux places p-adiques}.
\newblock In {\em On Certain L-Functions}, volume~13 of {\em Clay Mathematics
  Proceedings}, pages 333--376. American Mathematical Society, 2011.

\bibitem{Sakellaridis2017}
Y.~Sakellaridis.
\newblock Plancherel decomposition of {Howe} duality and {Euler} factorization
  of automorphic functionals.
\newblock In {\em Representation Theory, Number Theory, and Invariant Theory},
  pages 545--585. Springer International Publishing, 2017.

\bibitem{Sakellaridis2012}
Y.~Sakellaridis and A.~Venkatesh.
\newblock {\em Periods and Harmonic Analysis on Spherical Varieties}, volume
  396.
\newblock Ast\'{e}risque, 2017.

\bibitem{Venkatesh2005}
A.~Venkatesh.
\newblock {The Burger-Sarnak method and operations on the unitary dual of
  GL(n)}.
\newblock {\em Representation Theory of the American Mathematical Society},
  009(08):268--287, August 2005.

\bibitem{Wallach1988}
N.~Wallach.
\newblock {\em Real Reductive Groups I}.
\newblock Academic Pr Inc, 1988.

\end{thebibliography}
\end{document}